\documentclass[12pt,oneside]{amsart}

\usepackage{setspace}

\usepackage{amssymb,fontenc}   

\usepackage{dsfont}

\usepackage{amsmath}

\usepackage{mathrsfs}

\usepackage{stmaryrd}

\usepackage{amsthm}
\usepackage{amscd}
\usepackage{mathtools}
\usepackage{tikz-cd}
\usepackage{graphicx}

\usepackage{xcolor}
\usepackage{hyperref}
\hypersetup{
    colorlinks=true,
    citecolor=blue}

\usepackage{framed}
\usepackage{comment}

\usepackage{enumitem}

\usepackage[all]{xy}

\usepackage{textcomp}

\usepackage{amsbsy}

\usepackage{cancel}

    \DeclareFontFamily{U}{wncy}{}
    \DeclareFontShape{U}{wncy}{m}{n}{<->wncyr10}{}
    \DeclareSymbolFont{mcy}{U}{wncy}{m}{n}
    \DeclareMathSymbol{\Sha}{\mathord}{mcy}{"58}

\newtheorem{thm}{Theorem}[section]

\newtheorem{thm-defn}[thm]{Theorem/Definition}
\newtheorem{lem}[thm]{Lemma}
\newtheorem{prop}[thm]{Proposition}
\newtheorem{cor}[thm]{Corollary}

\theoremstyle{definition}
\newtheorem{defn}[thm]{Definition}
\newtheorem{notn}[thm]{Notation}
\newtheorem{eg}[thm]{Example}

\newtheorem{set-up}[thm]{Set-up}

\theoremstyle{remark}
\newtheorem{rem}[thm]{Remark}

\numberwithin{equation}{section}

\usepackage{relsize}
\usepackage[bbgreekl]{mathbbol}
\usepackage{amsfonts}

\DeclareSymbolFontAlphabet{\mathbb}{AMSb}
\DeclareSymbolFontAlphabet{\mathbbl}{bbold}
\newcommand{\Prism}{{\mathlarger{\mathbbl{\Delta}}}}
\DeclareMathOperator{\Gal}{Gal}

\DeclareMathOperator{\Spf}{Spf}

\DeclareMathOperator{\Ker}{Ker}

\newcommand{\rank}{\mathrm{rank}}

\newcommand{\cris}{\mathrm{cris}}

\newcommand{\st}{\mathrm{st}}

\newcommand{\dR}{\mathrm{dR}}

\newcommand{\N}{\mathbf{N}}

\newcommand{\Z}{\mathbf{Z}}
\newcommand{\Q}{\mathbf{Q}}

\newcommand{\B}{\mathbf{B}}

\renewcommand{\log}{\mathrm{log}}

\newcommand{\fkM}{{\mathfrak M}}

\newcommand{\fkS}{{\mathfrak S}}

\newcommand{\calC}{\mathcal{C}}

\newcommand{\calE}{\mathcal{E}}

\newcommand{\calO}{\mathcal{O}}

\theoremstyle{definition}
\newtheorem*{theoremA}{Theorem A}
\theoremstyle{plain}

\theoremstyle{definition}
\newtheorem*{theoremB}{Theorem B}
\theoremstyle{plain}%

\theoremstyle{definition}
\newtheorem*{theoremC}{Theorem C}
\theoremstyle{plain}%

\theoremstyle{definition}
\newtheorem*{ConstructionA}{Construction A}
\theoremstyle{plain}%

\theoremstyle{definition}
\newtheorem*{ConstructionB}{Construction B}
\theoremstyle{plain}

\begin{document}

\pagenumbering{arabic}

\title{A prismatic Herr complex for Bloch--Kato Selmer groups}

\author{Heng Du} 
\address{Yau Mathematical Sciences Center, Tsinghua University, Beijing 100084, China}
\email{hengdu@mail.tsinghua.edu.cn}

\author{Luming Zhao} 
\address{School of Mathematical Sciences, Dalian University of Technology, Dalian 116024, China}
\email{lumingzhao@dlut.edu.cn}

\begin{abstract} Let $p$ be an odd prime and $K$ a $p$-adic field. We introduce crystalline and semistable prismatic
$(\varphi,\tau)$-modules $(\mathfrak M,\varphi_\mathfrak M, \tau_\mathfrak M)$.
For a $\mathbf Z_p$-lattice $T$ of either type, its associated module $\mathfrak M$ gives an explicit two-term complex for the cohomology of
the corresponding $F$-crystal on the absolute prismatic or logarithmic
prismatic site of $\operatorname{Spf}(\mathcal O_K)$, respectively.
For effective $\mathfrak M$, the Frobenius fibre is an explicit
three-term Herr(-style) complex admitting a natural integral two-term model
with terms $\mathfrak M/u\mathfrak M$, a finite free $W(k)$-module.
For arbitrary Hodge--Tate weights, we also construct a localized version Herr complex.
In both constructions, first cohomology injects into the local
Bloch--Kato Selmer group $H_f^1(K,T)$ in the crystalline case and the
semistable local condition $H_{\st}^1(K,T)$ in the semistable case.
All these comparison maps become isomorphisms after inverting $p$,
while for the localized Herr complex, the comparison is also an integral
isomorphism for non-positive Hodge--Tate weights, with
convention $\operatorname{HT}(\mathbf Q_p(1))=-1$.
For $p$-divisible groups over $\mathcal O_K$, the good truncation of the localized Herr
complex recovers Nekov\'a\v r's complex in \cite{Nek18}.
In the effective crystalline case, we identify the Herr complex
with the syntomic cohomology of the associated reflexive $F$-gauge in \cite{bhatt-F-gauge},
giving explicit integral models over ramified as well as
unramified bases.
\end{abstract}

\maketitle

\tableofcontents

\section{Introduction}

\noindent
\textbf{Bloch--Kato Selmer groups and local conditions.}
For an elliptic curve $E$ over a number field $F$, the Birch and
Swinnerton--Dyer conjecture predicts that $\rank_\Z E(F)$ equals
the order of vanishing of $L(E,s)$ at $s=1$
(see \cite[\S1, conjecture~(A)]{TateBSD1966}).
For an abelian variety $A/F$, put
$\widehat{\mathrm{Sel}}_p(A/F)
:=\allowbreak\varprojlim_n\mathrm{Sel}_{p^n}(A/F)$, with transition maps
induced by multiplication by $p$. This compact $p$-adic Selmer group
fits into the exact sequence
\[
0\longrightarrow A(F)\otimes_\Z\Z_p
\longrightarrow\widehat{\mathrm{Sel}}_p(A/F)
\longrightarrow T_p\Sha(A/F)\longrightarrow0.
\]
Here $\Sha(A/F)$ is the Tate--Shafarevich group and
$T_p\Sha(A/F)=\varprojlim_n\Sha(A/F)[p^n]$.
The conjectured finiteness of $\Sha(A/F)$ would imply that
$\widehat{\mathrm{Sel}}_p(A/F)$ has $\Z_p$-rank
$\rank_\Z A(F)$
(see \cite[\S1, p.~420]{TateBSD1966} and
\cite[Chapter~I, Remark~6.7]{MilneADT2006}).

Bloch and Kato generalized this construction to geometric
$p$-adic Galois representations
(see \cite[\S3.7 and Definition~5.1]
{Bloch-Kato90L-functions-Tamagawa-number}).
For such a representation $V$ of $G_F$, their Selmer group is
\[
H_f^1(F,V)
=\Ker\!\big(
H^1(F,V)\longrightarrow
\prod_v H^1(F_v,V)/H_f^1(F_v,V)\big),
\]
where the product runs over the finite places of $F$ and
$F_v$ denotes the completion at $v$.
For $V=T_pA\otimes_{\Z_p}\Q_p$, the Kummer maps identify
$H_f^1(F,V)$ with
$\widehat{\mathrm{Sel}}_p(A/F)\otimes_{\Z_p}\Q_p$.
For $v\nmid p$, the Bloch--Kato local condition
is the unramified condition.
The main difficulty in this construction lies at the places above
$p$, where the local conditions are defined using Fontaine's
period rings and their integral structure is particularly subtle.
These local conditions and their semistable analogues are the
focus of this paper.
Let $p>2$, and let $K$ be a complete
discretely valued field of mixed characteristic $(0,p)$, with ring
of integers $\calO_K$ and perfect residue field $k$.
For $\star\in\{\cris,\st\}$ and a $p$-adic representation $V$ of
$G_K$ of type $\star$, put
\[
H_\star^1(K,V)
=\Ker\!\left(
H^1(K,V)\longrightarrow
H^1(K,\B_\star\otimes_{\Q_p}V)\right).
\]
Here $H_\cris^1(K,V)=H_f^1(K,V)$ is the local Bloch--Kato
Selmer group, and $H_\st^1(K,V)$ is the semistable local condition.
For a $G_K$-stable $\Z_p$-lattice $T\subset V$, define
$H_\star^1(K,T)$ as the inverse image of $H_\star^1(K,V)$
under $H^1(K,T)\to H^1(K,V)$
(see \eqref{s6:eq-integral-local-condition}).

We construct explicit complexes describing these local conditions.
We introduce \emph{crystalline and semistable prismatic
$(\varphi,\tau)$-modules}\footnote{See
Subsection~\ref{subsection finite height} for the general definition
and Definition~\ref{def:prismatic-phi-tau} for the effective case.
These modules can be viewed as variants of Du--Liu's prismatic
$(\varphi,\widehat G)$-modules in
\cite{du-liu-prismaticphiGhatmodule}, as well as refinements of
Caruso's $(\varphi,\tau)$-modules in \cite{Caruso_phitau}.}.
For a $G_K$-stable $\Z_p$-lattice $T$ in a representation of type
$\star\in\{\cris,\st\}$, write
$(\fkM,\varphi_{\fkM},\tau_{\fkM})$ for its associated module
of the corresponding type.
Here $\fkM$ is finite free over $\fkS=W(k)[[u]]$, and Frobenius
is defined with values in $\fkM[1/E(u)]$, where $E(u)$ is the
Eisenstein polynomial of a fixed uniformizer of $K$.
We first present the complexes and compare their cohomology
with the local conditions above, before explaining their
prismatic origin.

\begin{ConstructionA}
[see Definitions~\ref{def ambient-finite-height-Herr}
 and~\ref{def effective prismatic Herr complex}]
Let $T$ and $\fkM$ be as above, and put $\fkM_E=\fkM[1/E(u)]$.
Using the coefficient rings of
\cite[\S\S2.2--2.3]{du-liu-prismaticphiGhatmodule}
and their first structural maps, put
\[
\begin{aligned}
N_{\cris}(\fkM)
 &=\left(A^{(2)}\otimes_{p_0,\fkS}\fkM\right)_0,
 &N_{\st}(\fkM)
 &=\left(A_{\st}^{(2)}\otimes_{i_0,\fkS}\fkM\right)_0.
\end{aligned}
\]
The subscript-zero condition is defined in
Definition~\ref{def of sub-zero} and discussed below.
Let $N_\star^{\mathrm{fh}}(\fkM)$ be the subgroup defined by the
same condition in the corresponding coefficient module after
inverting $E(u_0)$ (see Subsection~\ref{s6 finite-height}).
The \emph{$E$-localized prismatic Herr complex} is
\[
\calC_{\varphi,\tau}^{\star,\mathrm{fh}}(T)
:=\big[
\fkM\xrightarrow{d^0_{\mathrm{fh}}}
\fkM_E\oplus N_\star(\fkM)
\xrightarrow{d^1_{\mathrm{fh}}}N_\star^{\mathrm{fh}}(\fkM)
\big]
\]
in degrees $0,1,2$, where
\[
d^0_{\mathrm{fh}}(m)
=\bigl((\varphi_{\fkM}-1)m,(\tau_{\fkM}-1)m\bigr),
\qquad
d^1_{\mathrm{fh}}(a,b)
=(\tau_{\fkM}-1)a-(\varphi-1)b.
\]
Here $(\tau_{\fkM}-1)x$ means
$\tau_{\fkM}(1\otimes x)-1\otimes x$, using the localized
extension of $\tau_{\fkM}$ when $x\in\fkM_E$.
If $\fkM$ is effective, the same differential formulas define the
\emph{effective prismatic Herr complex}
\[
\calC_{\varphi,\tau}^{\star}(T)
:=\big[
\fkM\xrightarrow{d^0}
\fkM\oplus N_\star(\fkM)
\xrightarrow{d^1}N_\star(\fkM)
\big]
\]
in degrees $0,1,2$.
This is the mapping fibre of $\varphi-1$ on the two-term complex
$[\fkM\xrightarrow{\tau_{\fkM}-1}N_\star(\fkM)]$
in degrees $0,1$.
\end{ConstructionA}

The superscript $\mathrm{fh}$ retains the finite-height notation
of Subsection~\ref{s6 finite-height}. The localization enlarges
only the Frobenius target, not the whole complex.

\begin{theoremA}
[{\normalfont = Theorems~\ref{s6 thm-fh-main}
 and~\ref{thm H0 H1}}]
Let $T$ be a $G_K$-stable $\Z_p$-lattice in a $p$-adic
representation of type $\star\in\{\cris,\st\}$.
There are natural maps
\[
\begin{aligned}
H^0\!\left(\calC_{\varphi,\tau}^{\star,\mathrm{fh}}(T)\right)
 &\xrightarrow{\sim} H^0(K,T),\\
H^1\!\left(\calC_{\varphi,\tau}^{\star,\mathrm{fh}}(T)\right)
 &\lhook\joinrel\longrightarrow H_\star^1(K,T).
\end{aligned}
\]
The second map is always an isomorphism after inverting $p$,
and is an isomorphism integrally when
$\operatorname{HT}(T)\subset[-r,0]$ for some $r\geq0$.
If $\fkM(T)$ is effective, the natural inclusion
$\calC_{\varphi,\tau}^{\star}(T)
\to\calC_{\varphi,\tau}^{\star,\mathrm{fh}}(T)$
induces isomorphisms on $H^0$ and $H^1$.
\end{theoremA}

For effective $\fkM=\fkM(T)$, put
$\overline{\fkM}=\fkM/u\fkM$, with induced
Frobenius $\overline\varphi$.
Theorem~\ref{s6 prop-mod-u} gives a natural integral quasi-isomorphism
\[
\calC_{\varphi,\tau}^{\star}(T)
\simeq
[\overline{\fkM}\xrightarrow{\overline\varphi-1}
                         \overline{\fkM}]
\]
in $D(\Z_p)$. Thus $H^i$ vanishes for $i\geq2$ in the
effective case. No such vanishing is asserted for the full
localized complex.

\noindent
\textbf{Prismatic nature of the construction.}
We now explain the prismatic origin of the coefficient rings
and descent operator appearing in Construction A.
Put $X=\Spf(\calO_K)$, and let
$M_X=\calO_X\cap\calO_X[1/p]^\times$ be its canonical log structure.
Write $X_{\Prism}$ for the absolute prismatic site of $X$ and
$(X,M_X)_{\Prism}$ for the corresponding absolute logarithmic
prismatic site. On either site, $\calO_{\Prism}$ denotes the
structure sheaf, which assigns to a prism its underlying ring
(see Section~\ref{sec:prismatic-F-crystals}).

The coefficient rings $A^{(2)}$ and $A_{\st}^{(2)}$ introduced by
the first author and Liu in
\cite[\S\S2.2--2.3]{du-liu-prismaticphiGhatmodule}
have a natural interpretation on these sites.
The ring $A^{(2)}$ is the underlying ring of the self-product of
the Breuil--Kisin prism $(\fkS,(E(u)))$ on $X_{\Prism}$.
Equipping this prism with the log structure associated with
$\N\to\fkS$, $1\mapsto u$, gives its logarithmic analogue,
whose self-product on $(X,M_X)_{\Prism}$ has underlying ring
$A_{\st}^{(2)}$
(see \cite[Proposition~4.1.3 and Lemma~5.0.12]
{du-liu-prismaticphiGhatmodule}).
Both rings embed into $A_{\inf}$, and $A_{\st}^{(2)}$ is stable
under the $G_K$-action
(see \cite[Corollary~2.4.5 and Lemma~3.3.1]
{du-liu-prismaticphiGhatmodule} and
Subsection~\ref{section Prismatic period rings}).
In this description, the semilinear operator $\tau_{\fkM}$
encodes the descent isomorphism over the corresponding
Breuil--Kisin self-product
(see Theorem~\ref{thm:exact-structures-module-categories}
and Corollary~\ref{cor:crystal-tau-fh}).

\noindent
\textbf{Herr complexes and the Kummer tower.}
The preceding construction is an analogue of the classical Herr
complex for local conditions. In classical $p$-adic Hodge theory,
Fontaine gave classifications of $p$-adic Galois representations
and important subcategories using several semilinear algebraic
objects, including notably the $(\varphi,\Gamma)$-modules
(see \cite{fontaine-p-adic-rep-i}).
Herr showed that Galois cohomology is computed by an explicit
three-term complex built from the associated $(\varphi,\Gamma)$-module
(see \cite[Theorem~2.1]{Her98}).
When $K/\Q_p$ is finite, his study of the $\psi$-operator gives
finiteness and Euler--Poincar\'e characteristic results, including
finite generation of $H^i(K,T)$ over $\Z_p$ for a $G_K$-stable
$\Z_p$-lattice $T$
(see \cite[p.~565 and \S\S3--4]{Her98}).
In particular, the integral Herr complex is perfect over $\Z_p$.
Herr also gave a new proof of local Tate duality using
$(\varphi,\Gamma)$-modules (see \cite{Her01}).

To study prismatic and syntomic cohomology and their relations with
Galois cohomology and its local conditions over both ramified and
unramified bases, we use the Kummer analogues of this theory.
These replace the cyclotomic tower
$K_{p^{\infty}}=\bigcup_{n\geq0}K(\zeta_{p^n})$ by
$K_{\infty}=\bigcup_{n\geq0}K(\pi_n)$, where the roots of unity
are compatible, $\pi_0$ is the fixed uniformizer, and
$\pi_{n+1}^p=\pi_n$.

\begin{rem}\label{rem:intro-herr-background}
\leavevmode
\begin{enumerate}
\item
In the Kummer setting, Breuil and Kisin studied semistable
$G_K$-representations using semilinear algebraic structures
constructed via the Kummer tower
(see \cite{Bre98,Bre99,kisin-crystalline}). Their work led to the theory
of \emph{Breuil--Kisin modules}, finite free $\fkS$-modules equipped
with a Frobenius structure of finite height.
For every prime $p$, the category of $p$-divisible groups over
$\calO_K$ is anti-equivalent to the category of finite free
effective Breuil--Kisin modules of height at most one
(see \cite[\S2.2]{kisin-crystalline} for $p>2$ and
\cite[Theorem~1.0.1]{LiuBT2} for $p=2$).
For a $G_K$-stable lattice in a semistable representation,
its Breuil--Kisin module describes the restriction to $G_{K_\infty}$.
Liu's $(\varphi,\widehat G)$-modules add a semilinear
$\widehat G$-action on a suitable scalar extension, subject to
Frobenius compatibility and additional invariance conditions.
They classify $G_K$-stable $\Z_p$-lattices in semistable
representations, and their realization restricts to the usual
$G_{K_\infty}$-realization of the underlying Breuil--Kisin module
(see \cite[Definition~2.2.3 and Theorem~2.3.1]{Liu10}).
In the same Kummer setting, Caruso establishes an equivalence
between the category of all finite-dimensional continuous
$\Q_p$-representations of $G_K$ and the corresponding category of
$(\varphi,\tau)$-modules, as well as an integral equivalence for
all finite free $\Z_p$-modules with continuous $G_K$-action
(see \cite[Theorem~1]{Caruso_phitau}).

\item
For the computation of Galois cohomology, Colmez and Liu extended
the classical cyclotomic theory to overconvergent and Robba-ring
coefficients (see \cite[\S2.1]{Col08} and
\cite[Theorems~0.2 and~2.8]{Liu07}).
In the Kummer setting, Tavares Ribeiro constructed a four-term
complex from $(\varphi,G_\infty)$-modules and proved that it
computes the Galois cohomology of the corresponding representation
(see \cite[Theorem~1.5]{TR11}).
The second author constructed a three-term Herr(-style) complex
from $(\varphi,\tau)$-modules computing the same cohomology,
together with an explicit quasi-isomorphism to Tavares Ribeiro's
complex (see \cite{Zha25}).
Gao and the second author developed an axiomatic cohomology theory
for the \'etale, overconvergent and rigid-overconvergent settings,
and established comparisons with $(\varphi,\Gamma)$-module
cohomology and Galois cohomology
(see \cite[Theorems~10.6 and~10.10]{GZ25}).

\item
Over an unramified base, Abhinandan constructed an explicit
complex using Wach modules whose first cohomology, after
inverting $p$, is naturally isomorphic to the local Bloch--Kato
Selmer group of the associated crystalline representation
(see \cite[Theorem~1.2]{Abhi25}).
Our prismatic construction also treats ramified bases and
semistable coefficients.
\end{enumerate}
\end{rem}

\begin{rem}\label{rem:intro-nekovar}
For a finite extension $K/\Q_p$ and a $p$-divisible group
$G/\calO_K$, Nekov\'a\v r's complex in \cite{Nek18}
is naturally isomorphic to the good truncation in degrees $0,1$
of our $E$-localized prismatic Herr complex.
Put $T=T_pG$.
Its Hodge--Tate weights lie in $[-1,0]$, so the associated
covariant prismatic module is generally non-effective.
The prismatic Herr complexes themselves are defined without
the assumption that $K/\Q_p$ is finite.
Propositions~\ref{prop nek-cochain-map} and~\ref{prop embedding}
give natural maps
\[
C_{\mathrm{Nek}}^\bullet(G)
\xrightarrow{\sim}
\tau_{\leq1}\calC_{\varphi,\tau}^{\cris,\mathrm{fh}}(T)
\lhook\joinrel\longrightarrow C_{\varphi,\tau}(T)
\simeq R\Gamma(K,T),
\]
where $\tau_{\leq1}$ denotes good truncation,
$C_{\mathrm{Nek}}^\bullet(G)$ is Nekov\'a\v r's complex in
\cite{Nek18}, and $C_{\varphi,\tau}(T)$ is the full three-term
$(\varphi,\tau)$-complex computing Galois cohomology in \cite{Zha25}.
The first arrow is an isomorphism of complexes, and the middle
arrow is termwise injective.
Nekov\'a\v r uses extensions with quotient $\fkS\{-1\}$.
Tensoring the extension by $\fkS\{1\}$ replaces this quotient
by the unit object and identifies the two extension problems.
Under this comparison, Theorem A recovers the identification of
$H^1(C_{\mathrm{Nek}}^\bullet(G))$ with the subgroup
$H_f^1(K,T)\subset H^1(K,T)$.
\end{rem}

\begin{rem}
The subscript-zero terms in Construction A reflect a structural
difference between the cyclotomic and Kummer-theoretic Herr
complexes recalled above.
The classical Herr complex is built from copies of a single
module on which both Frobenius and the cyclotomic action
are defined.
In the Kummer setting, the $\tau$-action is defined only after
extending the coefficient ring, and the three-term complex involves
a $\mathbf{Z}_p$-submodule of this larger module cut out by a
\emph{subscript-zero condition} (see \cite{Zha25}).
The appearance of this additional condition makes the analogy
with the classical Herr complex less direct, and its geometric
meaning is not apparent from the formula alone.

In their setting, Gao and the second author use cohomological
vanishing to obtain a two-term complex computing the continuous
$\widehat G$-cohomology of the enlarged coefficient module, with
degree-one term cut out by the subscript-zero condition
(see \cite[\S1.4 and Theorem~5.3]{GZ25}).
For the logarithmic prismatic coefficients considered here, the
first cohomology of the analogous two-term complex is naturally
identified with the kernel of restriction to $H_K$
(see Corollary~\ref{cor:relative-Galois-prismatic-H1}). Moreover, our prismatic construction gives the subscript-zero condition
a geometric interpretation: under the first-vertex identifications,
$N_\star(\fkM)$ is identified with the group of degree-one
\v{C}ech--Alexander cocycles for the corresponding Breuil--Kisin
cover
(see Proposition~\ref{Prop ker d^1_prism}).
In Construction A, the localized subgroup
$N_\star^{\mathrm{fh}}(\fkM)$ is defined by imposing the same
condition in the localized ambient module, rather than by
localizing $N_\star(\fkM)$ itself.
\end{rem}

\noindent
\textbf{Prismatic cohomology via $\tau$-complexes.}
With this cocycle interpretation in hand, we turn to the two-term
$\tau$-complexes underlying Construction A.
Their Frobenius mapping fibres are the three-term complexes
already introduced. More precisely,
\[
\calC_{\varphi,\tau}^{\star,\mathrm{fh}}(T)
\simeq\operatorname{fib}\!\left(
\varphi-\mathsf{can}\colon
[\fkM\xrightarrow{\tau_{\fkM}-1}N_\star(\fkM)]
\longrightarrow
[\fkM_E\xrightarrow{\tau_{\fkM}-1}N_\star^{\mathrm{fh}}(\fkM)]
\right),
\]
where $\mathsf{can}$ is induced by localization.
If $\fkM$ is effective, Frobenius preserves the integral terms,
and the mapping fibre of $\varphi-1$ on the source complex
is the effective prismatic Herr complex
$\calC_{\varphi,\tau}^{\star}(T)$ of Construction A.
The following comparison identifies the underlying two-term
complexes with prismatic cohomology and is of independent
interest, before taking the Frobenius fibre.

Prismatic cohomology, introduced in the seminal work of
Bhatt--Scholze, provides a unifying framework for integral
$p$-adic Hodge theory
(see \cite{bhatt-scholze-prismaticcohom,
bhatt-morrow-scholze-integralpadic}).
A prismatic $F$-crystal is a finite locally free
$\calO_{\Prism}$-crystal equipped, after inverting the prism ideal,
with an isomorphism from its Frobenius pullback to itself
(see Definition~\ref{def:general-prismatic-F-crystal}).
Bhatt--Scholze's classification gives the following result
(see \cite[Theorem~1.2]{bhatt-scholze-prismaticFcrystal}).

\begin{thm}
\label{Thm BS cat equivalenc of pris F-crystal and crystalline lattices_intro}
There is an equivalence between the category of prismatic $F$-crystals
in finite locally free $\calO_{\Prism}$-modules on
$X_{\Prism}$ and the category of $G_K$-stable $\Z_p$-lattices
in crystalline representations.
\end{thm}

The first author and Liu obtained the corresponding classification
of semistable lattices using the absolute logarithmic prismatic
site $(X,M_X)_{\Prism}$
(see \cite[Theorem~5.0.18]{du-liu-prismaticphiGhatmodule}).
Our crystalline and semistable prismatic
$(\varphi,\tau)$-modules
are equivalent to prismatic $F$-crystals on the two sites,
respectively (see Corollary~\ref{cor:crystal-tau-fh}), and through
covariant \'etale realization to lattices in crystalline and
semistable representations.
We describe the derived global sections of the corresponding
underlying prismatic crystals in $\calO_\Prism$-modules by two-term complexes.

We want to note that among the comparison theorems in prismatic
cohomology, the \'etale comparison is the least direct.
The Hodge--Tate, de Rham, and crystalline comparisons can be
obtained directly by evaluating prismatic cohomology on concrete
bounded prisms. By contrast, existing formulations of the \'etale
comparison pass to the perfect prismatic site, invert the
Hodge--Tate ideal sheaf, and invoke Katz's equivalence together
with $v$-descent
(see \cite[Section~3]{bhatt-scholze-prismaticFcrystal}).
This indirect route is unsatisfactory for explicit problems,
and inverting the Hodge--Tate ideal sheaf discards filtration
and torsion data. A central new input of this paper is the
construction of an explicit cochain complex over explicit period
rings that directly relates prismatic cohomology to certain
Galois cohomology groups
(see Section~\ref{sec:prismatic-cohomology-tau}, and \textit{cf.}
Theorem~\ref{prop compute prismatic coho}). 

Related explicit descriptions are available in the cyclotomic
setting. Over unramified bases, Abhinandan relates prismatic
$F$-crystals to Wach modules by identifying the Galois action
with a prismatic stratification
(see \cite[Theorems~1.1 and~1.4]{AbhinandanPrismaticWach}).
For the cyclotomic base $W(k)[\zeta_p]$, Gros--Le Stum--Quir\'os
give an explicit two-term description of prismatic cohomology
using absolute calculus, and relate their connection data to
$\Gamma$-modules
(see \cite[Theorem~7.12 and Proposition~7.14]{GLQAbsolute}).
Our contribution is a Breuil--Kisin approach via the Kummer
tower, with no restriction on the absolute ramification index
of $K$. It provides such a comparison for prismatic crystals
(in $\calO_\Prism$-modules) over $\calO_K$, uniformly in the
crystalline and semistable settings, without passing to the
perfect prismatic site. We expect this method to have further applications in integral
$p$-adic Hodge theory, which we plan to pursue in future work.

\begin{ConstructionB}
[see Theorem~\ref{prop compute prismatic coho}]
Let $T$ be a $G_K$-stable $\Z_p$-lattice in a semistable
representation, $\fkM_{\Prism}$ its associated prismatic
$F$-crystal on $(X,M_X)_{\Prism}$, and
$(\fkM,\varphi_{\fkM},\tau_{\fkM})$ its
semistable prismatic $(\varphi,\tau)$-module. Define
\[
\calC_{\Prism,\log}(T)
\coloneqq
\big[
\fkM\xrightarrow{\,\tau_{\fkM}-1\,}
\big(A_{\st}^{(2)}\otimes_{i_0,\fkS}\fkM\big)_0
\big].
\]
If $T[1/p]$ is crystalline, let $\fkM_{\Prism}$ also denote the
corresponding crystal on $X_{\Prism}$, whose logarithmic pullback
is the preceding crystal. Define
\[
\calC_{\Prism}(T)
\coloneqq
\big[
\fkM\xrightarrow{\,\tau_{\fkM}-1\,}
\big(A^{(2)}\otimes_{p_0,\fkS}\fkM\big)_0
\big].
\]
Both complexes are in degrees $0,1$, with differential
$m\mapsto\tau_{\fkM}(1\otimes m)-1\otimes m$.
The degree-one terms are $N_\st(\fkM)$ and $N_\cris(\fkM)$
from Construction A, respectively.
Under the first-vertex identifications, each is the subgroup of
degree-one \v{C}ech--Alexander cocycles for the corresponding
Breuil--Kisin cover (see Proposition~\ref{Prop ker d^1_prism}).
\end{ConstructionB}

\begin{theoremB}
[{\normalfont = Theorem~\ref{prop compute prismatic coho}}]
In the notation of Construction B, there is a natural
quasi-isomorphism in $D(\Z_p)$
\[
R\Gamma\bigl((X,M_X)_{\Prism},\fkM_{\Prism}\bigr)
\xrightarrow{\sim}\calC_{\Prism,\log}(T).
\]
If $T[1/p]$ is crystalline, there is also a natural quasi-isomorphism
\[
R\Gamma\bigl(X_{\Prism},\fkM_{\Prism}\bigr)
\xrightarrow{\sim}\calC_{\Prism}(T).
\]
In particular, these prismatic cohomology complexes are concentrated
in degrees $0,1$. No effectivity assumption is required here.
\end{theoremB}

\noindent
\textbf{Relation with $F$-gauges and syntomic cohomology.}
In the effective crystalline case, the Frobenius fibre above also
has a syntomic interpretation.
Bhatt--Lurie define perfect prismatic $F$-gauges on a $p$-adic
formal scheme $X$ as perfect complexes on its syntomification
$X^{\mathrm{Syn}}$
(see \cite[Definition~6.1.1]{bhatt-F-gauge}).
Their classification theorem gives equivalences
(see \cite[Theorem~6.6.13]{bhatt-F-gauge})
\[
\operatorname{Coh}^{\mathrm{refl}}(\calO_K^{\mathrm{Syn}})
\simeq\operatorname{Vect}^{\varphi}((\calO_K)_{\Prism})
\simeq\operatorname{Rep}^{\cris}_{\Z_p}(G_K).
\]
Over $\Z_p$, their cohomological comparison recovers
Bloch--Kato local conditions after inverting $p$
(see \cite[Proposition~6.7.3 and Remark~6.7.4]{bhatt-F-gauge}).
The question of extending this comparison to ramified bases is
discussed in Bhatt's ICM survey (see \cite[\S2.4]{Bha21}).
We obtain the following integral comparison in the effective
crystalline case.

\begin{theoremC}
[{\normalfont = Theorems~\ref{thm main on coho of F-gauge}
 and~\ref{s6 prop-mod-u}}]
Let $T$ be a $G_K$-stable $\Z_p$-lattice in a crystalline
representation with non-negative Hodge--Tate weights. Let
$\calE_T\in\operatorname{Coh}^{\mathrm{refl}}
(\calO_K^{\mathrm{Syn}})$ be the associated reflexive $F$-gauge.
Put $\fkM=\fkM(T)$ and
$\overline{\fkM}=\fkM/u\fkM$.
There are natural integral quasi-isomorphisms in $D(\Z_p)$
\[
R\Gamma(\calO_K^{\mathrm{Syn}},\calE_T)
\simeq\calC_{\varphi,\tau}^{\cris}(T)
\simeq[\overline{\fkM}\xrightarrow{\overline\varphi-1}\overline{\fkM}],
\]
where $\overline\varphi$ is induced by $\varphi_{\fkM}$
and the last complex is in degrees $0,1$.
\end{theoremC}

Combining Theorems A and C, we obtain natural isomorphisms

\[
H^0(\calO_K^{\mathrm{Syn}},\calE_T)
  \simeq H^0(K,T), \qquad
H^1(\calO_K^{\mathrm{Syn}},\calE_T)[1/p]
  \simeq H_f^1(K,T[1/p]).
\]

\noindent
\textbf{Some related works.}
Comparisons with earlier Herr, Wach-module and Nekov\'a\v{r}
constructions appear in
Remarks~\ref{rem:intro-herr-background}--\ref{rem:intro-nekovar}.
Gao--Min--Wang's cohomological-amplitude bound for finite
locally free $\mathcal{O}_\Prism$-crystals on the absolute
and logarithmic prismatic sites is an input to Theorem~B
\cite[Theorem~4.10]{gao2023hodge}.
Our contribution is the explicit presentation by
$\tau$-complexes and the identification of the
subscript-zero condition with the degree-one
\v{C}ech--Alexander cocycle condition.

For finite extensions $K/\mathbf{Q}_p$, Pentland gives
an integral Nygaard-filtered total-fibre description of
syntomic cohomology with reflexive $F$-gauge coefficients
\cite[Lemma~4.11]{PentlandSyntomification}, and identifies
its first cohomology with $H_f^1$ after inverting $p$
\cite[Proposition~4.3]{PentlandSyntomification}.
His effective-case argument also gives integral
degree-two vanishing
\cite[proof of Proposition~4.12]{PentlandSyntomification}.
Theorem~C provides an explicit prismatic Herr realization
of the effective syntomic complex and its natural
integral reduction to
$[\overline{\mathfrak M}
\xrightarrow{\bar\varphi-1}\overline{\mathfrak M}]$,
where $\overline{\mathfrak M}=\mathfrak M/u\mathfrak M$.
The vanishing statement alone is not claimed new.

Theorem~A uses \v{C}oupek--Gazaki--Marmora's exactness
results \cite[Theorem~3.12 and Remark~3.13(3)]{CGM26}.
Their Corollary~A.9 identifies syntomic $H^1$ integrally
with $H_f^1(K,T)$ for non-positive Hodge--Tate weights
in our convention.
Theorem~C instead gives an integral complex computing
syntomic cohomology in the non-negative, effective range.
In this range, our comparison with $H_f^1(K,T)$ is an
integral injection that becomes an isomorphism after inverting $p$. The surjectivity without inverting $p$ can fail, cf. Example~\ref{s6 example-missing-torsion}.

Min proves conditional congruence results for reflexive
$F$-gauges and constructs a natural injection from
syntomic $H^1$ into $H_f^1$
\cite[Theorem~1.1 and Lemma~1.2]{Min26}.
In the effective case, Corollary~\ref{s6 cor-congruences}
gives a direct reduction formula using the reduced
Frobenius module
$(\overline{\mathfrak M}/p^n\overline{\mathfrak M},
\bar\varphi)$.

\noindent
\textbf{Acknowledgements.}
This work was motivated by two questions raised by
Laurent Berger and Tong Liu, who served as referees for
the second author's doctoral thesis.
Laurent Berger asked about the geometric meaning of the
subscript-zero condition in the $(\varphi,\tau)$-complex,
while Tong Liu asked whether a variant of this complex
built from Breuil--Kisin modules could compute local
Bloch--Kato Selmer groups.
The second author is grateful to both for these stimulating questions. The project began around the start of 2023, when we observed that the subscript-zero condition admits a prismatic interpretation. We also thank Abhinandan, Yangyu Fan, Yu Min, and Yupeng Wang for helpful correspondence.

Part of this work was carried out during the second
author's postdoctoral stays at Purdue University and
Peking University. The second author thanks both
institutions for their hospitality and excellent
research environments.

The first author was partially supported by the National Key R\&D Program of China No.~2023YFA1009703, the National Natural Science Foundation of China No.~12275066, and the Beijing Natural Science Foundation under Grant No.~1254044. The second author was supported by the National Natural Science Foundation of China No. 12601022, and the Fundamental Research Funds for the Central Universities No. DUT26RC(3)014.

\noindent
\textbf{Development of the paper and use of AI.}
The main results of this paper, in their original form, had been established by late 2025 and presented by the authors in several talks. Subsequent revisions extended Theorem~\ref{prop compute prismatic coho}  to include also the semistable case, strengthened integral results in Section~\ref{sec:prismatic-Herr-complex}, and added the comparison with Nekov\'a\v r's complex in \cite{Nek18} at Yangyu Fan's suggestion. ChatGPT and Codex (OpenAI) were used to check arguments for gaps, locate and supplement references, and improve exposition and grammar. In particular, AI suggested the idea in Lemma~\ref{lem:faithful-log-Taylor}, which led to the more unified proof of Theorem~\ref{prop compute prismatic coho} presented here. The authors take full responsibility for the mathematical arguments, references, and final text.

\section{Notation}\label{sec notation}

\begin{notn}\label{notation}
Let $p>2$ be a prime. Let $K$ be a complete discretely valued field of
mixed characteristic $(0,p)$, with ring of integers $\mathcal O_K$ and
perfect residue field $k$, and put
$K_0\coloneqq W(k)[1/p]$. Fix an algebraic closure $\overline K$ of $K$. Let $C_p$ be its
$p$-adic completion, with ring of integers $\mathcal O_{C_p}$, and set
\[
\mathcal O_{C_p}^{\flat}
\coloneqq
\varprojlim_{x\mapsto x^p}\mathcal O_{C_p}/p,
\qquad
A_{\inf}\coloneqq W(\mathcal O_{C_p}^{\flat}).
\]
We denote Fontaine's map by $\theta\colon A_{\inf}\longrightarrow \mathcal O_{C_p}$. Fix a uniformizer $\pi\in\mathcal O_K$, and let
$E(u)\in W(k)[u]$ be its monic minimal polynomial over $K_0$, thus
$E(u)$ is an Eisenstein polynomial. We also fix compatible systems
\[
\widetilde\pi=(\pi_n)_{n\geq 0}\in\mathcal O_{C_p}^{\flat},
\qquad
\pi_0=\pi,
\qquad
\pi_{n+1}^p=\pi_n,
\]
and
\[
\varepsilon=(\zeta_{p^n})_{n\geq 0}\in\mathcal O_{C_p}^{\flat},
\qquad
\zeta_{p^0}=1,
\qquad
\zeta_p\neq 1,
\qquad
\zeta_{p^{n+1}}^p=\zeta_{p^n}.
\]
Write $[\widetilde\pi]$ and $[\varepsilon]$ for their Teichm\"uller lifts
in $A_{\inf}$, and abbreviate $\epsilon\coloneqq[\varepsilon]$.

Put
\[
K_{\infty}\coloneqq\bigcup_{n\geq 0}K(\pi_n),
\qquad
K_{p^\infty}\coloneqq\bigcup_{n\geq 1}K(\zeta_{p^n}),
\qquad
L\coloneqq K_{\infty}K_{p^\infty},
\]
and
\[
\begin{aligned}
G_{K_\infty}&\coloneqq\Gal(\overline K/K_\infty),
&
G_{K_{p^\infty}}&\coloneqq\Gal(\overline K/K_{p^\infty}), & G_L \coloneqq\Gal(\overline K/L).
\end{aligned}
\]

Let $\mathfrak S\coloneqq W(k)[[u]]$ equipped with the Frobenius extending the Witt vector Frobenius and satisfying
$\varphi(u)=u^p$. The assignment
\[
\mathfrak S\longrightarrow A_{\inf},
\qquad
u\longmapsto[\widetilde\pi],
\]
is $\varphi$-equivariant. Under this map, $\xi\coloneqq E([\widetilde\pi])$ generates $\Ker(\theta)$.

Let $\chi\colon G_K\to\mathbf Z_p^\times$
be the cyclotomic character and put $\Gamma\coloneqq\Gal(K_{p^\infty}/K)$. The character $\chi$ identifies $\Gamma$ with an open subgroup of
$\mathbf Z_p^\times$. Since $p>2$, the group $\Gamma$ is procyclic, moreover in this case one has $K_\infty\cap K_{p^\infty}=K$ and restriction identifies $H_K\coloneqq\Gal(L/K_\infty)\simeq\Gamma$. We can choose a topological generator
$\gamma$ such that $\chi(\gamma)$ is a positive integer,
and use the same symbol for the corresponding element of $H_K$.

We further choose a topological generator $\tau\in\Gal(L/K_{p^\infty})\simeq\mathbf Z_p$ such that
\[
\tau(\pi_n)=\zeta_{p^n}\pi_n
\quad(n\geq 0),
\]
or equivalently
\[
\tau(\widetilde\pi)=\varepsilon\widetilde\pi,
\qquad
\tau([\widetilde\pi])=\epsilon[\widetilde\pi].
\]
Then
\[
\widehat G\coloneqq\Gal(L/K)
\simeq
\overline{\langle\tau\rangle}\rtimes H_K
\text{ and } 
\gamma\tau\gamma^{-1}=\tau^{\chi(\gamma)}.
\]

The fields and Galois groups used below are summarized by
\[
\begin{tikzcd}[every arrow/.append style={dash}]
&\overline K
\arrow[ldd,bend right,"G_{K_\infty}"']
\arrow[d,"G_L"]
\arrow[rdd,bend left,"G_{K_{p^\infty}}"]
&\\
&L
\arrow[ld,"H_K"']
\arrow[dd,"\widehat G"]
\arrow[rd,"\overline{\langle\tau\rangle}"]
&\\
K_\infty\arrow[rd]&&K_{p^\infty}\arrow[ld,"\Gamma=\overline{\langle\gamma\rangle}"]\\
&K&
\end{tikzcd}
\]

\noindent
\textbf{Convention.} Throughout the paper we use $\operatorname{HT}(\mathbf Q_p(1))=-1$. Our \'etale realization is covariant.  Thus an effective prismatic
$F$-crystal of height at most $h$ corresponds to a lattice with
Hodge--Tate weights in $[0,h]$.  Notice that Du--Liu use a contravariant Hom
realization in \cite{du-liu-prismaticphiGhatmodule}, when using their functor, we will naturally pass to its dual (see Notation \ref{notn:covariant-realization}).
\end{notn}

\section{\texorpdfstring{Prismatic period rings and prismatic
$F$-crystals}{Prismatic period rings and prismatic F-crystals}}
\label{sec:prismatic-F-crystals}

We begin by recalling the definition and basic facts about the prismatic period rings and prismatic
$F$-crystals used in this paper, following
Bhatt--Scholze~\cite{bhatt-scholze-prismaticcohom,bhatt-scholze-prismaticFcrystal},
Koshikawa~\cite{koshikawa-log-prismatic}, and
Du--Liu~\cite{du-liu-prismaticphiGhatmodule}. One new input is Theorem~\ref{thm:multivariable-DL-injectivity}, which can be viewed as a multivariable analog of \cite[Proposition~2.4.3]{du-liu-prismaticphiGhatmodule}.
This will be used to compare prismatic complexes with certain complexes defined by Galois modules.

\subsection{The two prismatic sites}

Put $X=\operatorname{Spf}(\mathcal O_K)$. Let $M_X=\mathcal O_X\cap\mathcal O_X[1/p]^\times$ be the log structure associated with $\mathbf N\to\mathcal O_K$, $1\mapsto\pi$, and let $M_X^{\mathrm{triv}}=\mathcal O_X^\times$ be the trivial log structure. For $N_X\in\{M_X^{\mathrm{triv}},M_X\}$, we consider the strict absolute logarithmic prismatic site of \cite[Definition~2.3]{DLMS2}. Its opposite category consists of bounded log prisms $(B,J,M_B)$ equipped with a strict morphism
\[
\bigl(\operatorname{Spf}(B/J),M_B|_{B/J}\bigr)
\longrightarrow (X,N_X),
\]
where $M_B|_{B/J}$ denotes the pullback log structure. A covering is a strict morphism $(B,J,M_B)\to(C,JC,M_C)$ such that $B\to C$ is $(p,J)$-completely faithfully flat. By faithfully flat descent, the assignments
\[
\mathcal O_{\Prism}(B,J,M_B)=B,
\qquad
\mathcal I_{\Prism}(B,J,M_B)=J
\]
define sheaves on this site, and the $\delta$-structure induces the Frobenius endomorphism $\varphi$ of $\mathcal O_{\Prism}$.

\begin{rem}
\label{rem:trivial-log-site}
We omit the superscript $\mathrm{str}$ used in \cite{DLMS2}. When $N_X=M_X^{\mathrm{triv}}$, strictness implies that $M_B|_{B/J}$ is trivial, so by \cite[Lemma~2.1]{DLMS2} applied to the trivial chart, $M_B$ is trivial as well. Forgetting the log structure therefore identifies the site for $N_X=M_X^{\mathrm{triv}}$ with the absolute prismatic site $X_{\Prism}$. For $N_X=M_X$, we denote the site by $(X,M_X)_{\Prism}$ or $X_{\Prism_{\log}}$.
\end{rem}

\begin{eg}
The Breuil--Kisin prism and its logarithmic counterpart are
\[
(\mathfrak S,(E(u)))\in X_{\Prism},
\qquad
(\mathfrak S,(E(u)),\mathbf N)^a\in(X,M_X)_{\Prism},
\quad 1\longmapsto u,
\]
where $(-)^a$ denotes the associated log prism and the identification $\mathfrak S/(E(u))\simeq\mathcal O_K$ sends $u$ to $\pi$. The sheaves they represent cover the final objects of the respective topoi (see \cite[Example~2.6]{bhatt-scholze-prismaticFcrystal} and \cite[Lemma~4.1.8 and Proposition~5.0.16]{du-liu-prismaticphiGhatmodule})\footnote{We note that \cite[Proposition~5.0.16]{du-liu-prismaticphiGhatmodule} is due to Teruhisa Koshikawa.}. Put $\mathbf u=[\widetilde\pi]\in A_{\inf}$ and $\xi=E(\mathbf u)$, where $\theta\colon A_{\inf}\to\mathcal O_{C_p}$ has kernel $(\xi)$ and identifies $A_{\inf}/(\xi)$ with $\mathcal O_{C_p}$. Sending $u\mapsto\mathbf u$ defines morphisms from $(\mathfrak S,(E(u)))$ and $(\mathfrak S,(E(u)),\mathbf N)^a$ to $(A_{\inf},(\xi))$ and $(A_{\inf},(\xi),\mathbf N)^a$ (with prelog structure $1\mapsto\mathbf u$), respectively.
\end{eg}

\subsection{\texorpdfstring{$F$-crystals and Breuil--Kisin
modules}{F-crystals and Breuil--Kisin modules}}

We recall the notions of prismatic $F$-crystals and Breuil--Kisin
modules (see \cite[Definition~4.1 and Example~4.3]
{bhatt-scholze-prismaticFcrystal}).
For the logarithmic setting, compare
\cite[Definition~5.0.17]{du-liu-prismaticphiGhatmodule}.
We allow non-effective objects and specify our height conventions below.

\begin{defn}
\label{def:general-prismatic-F-crystal}
A \emph{crystal in finite locally free $\mathcal O_{\Prism}$-modules} is a sheaf $\mathcal E$ whose value on each $(B,J,M_B)$ is a finite projective $B$-module such that every morphism $(B,J,M_B)\to(C,JC,M_C)$ induces an isomorphism
\[
C\otimes_B\mathcal E(B,J,M_B)
\xrightarrow{\sim}\mathcal E(C,JC,M_C).
\]
A \emph{prismatic $F$-crystal} is such a crystal equipped with an isomorphism
\[
\Phi_{\mathcal E}\colon
(\varphi^*\mathcal E)[1/\mathcal I_{\Prism}]
\xrightarrow{\sim}\mathcal E[1/\mathcal I_{\Prism}],
\qquad
\text{where }\varphi^*\mathcal E=
\mathcal O_{\Prism}\otimes_{\varphi,\mathcal O_{\Prism}}\mathcal E.
\]
For integers $a\leq b$, we say that $(\mathcal E,\Phi_{\mathcal E})$ has \emph{height in $[a,b]$} if
\[
\mathcal I_{\Prism}^{b}\mathcal E
\subseteq\Phi_{\mathcal E}(\varphi^*\mathcal E)
\subseteq\mathcal I_{\Prism}^{a}\mathcal E
\]
inside $\mathcal E[1/\mathcal I_{\Prism}]$, where negative powers of the invertible ideal $\mathcal I_{\Prism}$ are taken inside $\mathcal O_{\Prism}[1/\mathcal I_{\Prism}]$. We say that $(\mathcal E,\Phi_{\mathcal E})$ is \emph{effective} if $\Phi_{\mathcal E}(\varphi^*\mathcal E)\subseteq\mathcal E$. In that case, $\Phi_{\mathcal E}\colon\varphi^*\mathcal E\to\mathcal E$ is injective since $\varphi^*\mathcal E$ is $\mathcal I_{\Prism}$-torsion-free, and $(\mathcal E,\Phi_{\mathcal E})$ has height at most $h\geq0$ if it has height in $[0,h]$, or equivalently if $\mathcal I_{\Prism}^h$ annihilates $\operatorname{Coker}(\Phi_{\mathcal E}\colon\varphi^*\mathcal E\to\mathcal E)$.

Morphisms of prismatic $F$-crystals are morphisms of crystals commuting with $\Phi$. We denote the resulting categories on $X_{\Prism}$ and $(X,M_X)_{\Prism}$ by $\operatorname{Vect}^{\varphi}(X_{\Prism})$ and $\operatorname{Vect}^{\varphi}((X,M_X)_{\Prism})$, and write $\operatorname{Vect}_{[a,b]}^{\varphi}(X_{\Prism})$ and $\operatorname{Vect}_{[a,b]}^{\varphi}((X,M_X)_{\Prism})$ for the full subcategories of objects with height in $[a,b]$.
\end{defn}

\begin{defn}
\label{def:Breuil-Kisin-module}
A \emph{finite free Breuil--Kisin module} is a finite free $\mathfrak S$-module $\mathfrak M$ equipped with an isomorphism
\[
\Phi_{\mathfrak M}\colon
(\varphi^*\mathfrak M)[1/E(u)]\xrightarrow{\sim}
\mathfrak M[1/E(u)],
\qquad \text{where }
\varphi^*\mathfrak M=
\mathfrak S\otimes_{\varphi,\mathfrak S}\mathfrak M.
\]
Its semilinear Frobenius is $\varphi_{\mathfrak M}(m)=\Phi_{\mathfrak M}(1\otimes m)$. For integers $a\leq b$, it has height in $[a,b]$ if
\[
E(u)^b\mathfrak M\subseteq
\Phi_{\mathfrak M}(\varphi^*\mathfrak M)
\subseteq E(u)^a\mathfrak M
\]
inside $\mathfrak M[1/E(u)]$. It is \emph{effective} if $\varphi_{\mathfrak M}(\mathfrak M)\subseteq\mathfrak M$, and an effective object has height at most $h\geq0$ if it has height in $[0,h]$. Morphisms are $\mathfrak S$-linear maps commuting with $\Phi$.
\end{defn}

\begin{rem}\label{rmk def of fh case}
\item (1) Every finite free Breuil--Kisin module has height in some finite interval, since the matrices of $\Phi_{\mathfrak M}$ and $\Phi_{\mathfrak M}^{-1}$ have bounded $E(u)$-denominators. We write $\operatorname{Mod}^{\varphi}_{\mathfrak S}$ for the category obtained by allowing finitely generated $\mathfrak S$-modules. Because $\varphi\colon\mathfrak S\to\mathfrak S$ is flat, $\operatorname{Mod}^{\varphi}_{\mathfrak S}$ is abelian with kernels, cokernels, and cohomology of complexes computed on underlying $\mathfrak S$-modules.

\item (2) Evaluating a prismatic $F$-crystal on $(\mathfrak S,(E(u)))$ or $(\mathfrak S,(E(u)),\mathbf N)^a$ yields a finite free Breuil--Kisin module. By faithfully flat descent along the Breuil--Kisin cover, effectivity and height in $[a,b]$ can be checked after this evaluation (see \cite[Proposition~3.26]{DLMS} and its proof). In particular, every prismatic $F$-crystal has height in some finite interval.
\end{rem}

For $\star\in\{\cris,\st\}$, let $\operatorname{Rep}^{\star}_{\mathbf Z_p}(G_K)$ denote the category of $G_K$-stable $\mathbf Z_p$-lattices in crystalline, respectively semistable, $\mathbf Q_p$-representations of $G_K$.

\begin{thm}
\label{thm:BS-DL-crystalline}
\label{thm:BS-prismatic-F-crystals}
\label{Thm BS cat equivalenc of pris F-crystal and crystalline lattices}
The covariant \'etale realization functor induces equivalences
\[
\begin{aligned}
\operatorname{Vect}^{\varphi}(X_{\Prism})
&\xrightarrow{\sim}\operatorname{Rep}^{\cris}_{\mathbf Z_p}(G_K),\\
\operatorname{Vect}^{\varphi}((X,M_X)_{\Prism})
&\xrightarrow{\sim}\operatorname{Rep}^{\st}_{\mathbf Z_p}(G_K),
\end{aligned}
\]
under which effective objects of height at most $h$ correspond to lattices with Hodge--Tate weights in $[0,h]$.
\end{thm}

\begin{proof}
The equivalences are \cite[Theorem~1.2]{bhatt-scholze-prismaticFcrystal} and \cite[Theorem~1.3(2)]{DLMS2}, and the height assertions follow from \cite[Theorems~4.1.10 and~5.0.18]{du-liu-prismaticphiGhatmodule} after taking $\mathbf Z_p$-duals in accordance with Section~\ref{sec notation}.
\end{proof}

\subsection{\texorpdfstring{The rings $A^{(2)}$ and
$A_{\st}^{(2)}$}{The rings A(2) and Ast(2)}}
\label{subsec:period-rings-A2-Ast2}
\label{section Prismatic period rings}

Recall $\epsilon=[\varepsilon]$ and put $\mu=\epsilon-1$ in $A_{\inf}$. Following \cite[\S\S2.1--2.3]{du-liu-prismaticphiGhatmodule}, set
\[
\begin{aligned}
A^{(2)}&=W(k)[[u_0,u_1]]
\left\{\frac{u_1-u_0}{E(u_0)}\right\}_{\delta}^{\wedge},
&z&=\frac{u_1-u_0}{E(u_0)},\\
A_{\st}^{(2)}&=W(k)[[u_0,y]]
\left\{\frac{y}{E(u_0)}\right\}_{\delta}^{\wedge},
&w&=\frac{y}{E(u_0)},
\end{aligned}
\]
where $\delta(u_0)=\delta(u_1)=0$ on $W(k)[[u_0,u_1]]$, while $\delta(u_0)=0$ and $\varphi(y)=(1+y)^p-1$ on $W(k)[[u_0,y]]$. Here $\{-\}_{\delta}^{\wedge}$ denotes the $(p,E(u_0))$-completed $\delta$-envelope. Put $u_1=u_0(1+y)\in A_{\st}^{(2)}$, and define morphisms of $\delta$-rings
\[
\begin{aligned}
p_0,p_1\colon\mathfrak S&\longrightarrow A^{(2)},
&p_j(u)&=u_j,\\
i_0,i_1\colon\mathfrak S&\longrightarrow A_{\st}^{(2)},
&i_j(u)&=u_j,
\end{aligned}
\]
so that $i_0(E(u))=E(u_0)$ and $i_1(E(u))=E(u_1)$. Define $\delta$-ring homomorphisms
\[
A^{(2)}\xrightarrow{\iota_{\st}}A_{\st}^{(2)}
\xrightarrow{\iota_{\inf}^{\st}}A_{\inf}
\]
by $(u_0,u_1,z)\mapsto(u_0,u_0(1+y),u_0w)$ and $(u_0,y,w)\mapsto(\mathbf u,\mu,\mu/\xi)$, and put $\iota_A=\iota_{\inf}^{\st}\circ\iota_{\st}$. In particular, $\iota_A(u_1)=\epsilon\mathbf u$ and $\iota_A(z)=\mu\mathbf u/\xi$.

\begin{prop}
\label{prop:DL-period-rings}
The prisms $(A^{(2)},(E(u_0)))$ and $(A_{\st}^{(2)},(E(u_0)))$, with the exactified log structure in the second case, are the self-products of $(\mathfrak S,(E(u)))$ and $(\mathfrak S,(E(u)),\mathbf N)^a$ in $X_{\Prism}$ and $(X,M_X)_{\Prism}$, respectively. The maps $\iota_{\st}$ and $\iota_{\inf}^{\st}$ are injective, and $A_{\st}^{(2)}\subset A_{\inf}$ is $G_K$-stable with $G_K$-action factoring through $\widehat G$.
\end{prop}

\begin{proof}
The first assertion is \cite[Proposition~4.1.3 and Lemma~5.0.12]{du-liu-prismaticphiGhatmodule}, and the remaining claims are \cite[Corollary~2.4.5 and Lemma~3.3.1]{du-liu-prismaticphiGhatmodule}.
\end{proof}

\subsection{Self-products and evaluation maps}
\label{subsec:multivariable-DL-injectivity}
\label{subsec:multivariable-evaluation-injectivity}

Set $A^{(1)}=A_{\st}^{(1)}=\mathfrak S$. For $m\geq1$, put
\[
\begin{aligned}
A^{(m+1)}&=W(k)[[u_0,\ldots,u_m]]
\left\{z_j=\frac{u_j-u_0}{E(u_0)}\right\}_{1\leq j\leq m,\delta}^{\wedge},\\
A_{\st}^{(m+1)}&=W(k)[[u_0,y_1,\ldots,y_m]]
\left\{w_j=\frac{y_j}{E(u_0)}\right\}_{1\leq j\leq m,\delta}^{\wedge},
\end{aligned}
\]
with $(p,E(u_0))$-adic completions and $\delta$-structures given by $\delta(u_j)=0$ on $W(k)[[u_0,\ldots,u_m]]$, and by $\delta(u_0)=0$ and $\varphi(y_j)=(1+y_j)^p-1$ on $W(k)[[u_0,y_1,\ldots,y_m]]$. In $A_{\st}^{(m+1)}$, put $u_j=u_0(1+y_j)$ for $1\leq j\leq m$. For $m=1$, we write $z=z_1$, $y=y_1$, and $w=w_1$ as in Subsection~\ref{section Prismatic period rings}, so for maps $A_{\st}^{(2)}\to A_{\st}^{(3)}$ we write $y$ in $A_{\st}^{(2)}$ and $y_1,y_2$ in $A_{\st}^{(3)}$. Define
\[
j_m\colon A^{(m+1)}\longrightarrow A_{\st}^{(m+1)},
\qquad u_k\longmapsto u_k,\quad z_j\longmapsto u_0w_j,
\]
so that $j_1=\iota_{\st}$.

Write $\mathcal O_K\langle T_1,\ldots,T_m\rangle_{\mathrm{PD}}^{\wedge}$ for the $p$-adic completion of $\mathcal O_K[T_j^n/n!\mid n\geq1,\,1\leq j\leq m]$ inside $K[T_1,\ldots,T_m]$. Every element of $\mathcal O_K\langle T_1,\ldots,T_m\rangle_{\mathrm{PD}}^{\wedge}$ has a unique expansion $\sum_{\nu\in\mathbf N^m}b_\nu T^\nu/\nu!$ with $b_\nu\in\mathcal O_K$ tending to zero $p$-adically, where $\nu!=\prod_{j=1}^m\nu_j!$.

\begin{prop}
\label{prop:multivariable-coproducts}
The rings $A^{(m+1)}$ and $A_{\st}^{(m+1)}$ represent the $(m+1)$-fold self-products of $(\mathfrak S,(E(u)))$ and $(\mathfrak S,(E(u)),\mathbf N)^a$ in $X_{\Prism}$ and $(X,M_X)_{\Prism}$, respectively, and are flat over $\mathfrak S$ under each projection $p_k\colon\mathfrak S\to A^{(m+1)}$ and $i_k\colon\mathfrak S\to A_{\st}^{(m+1)}$ ($0\leq k\leq m$). Moreover, the elements $E(u_k)$ generate the same prism ideal for all $0\leq k\leq m$, and
\[
\begin{aligned}
A^{(m+1)}/E(u_0)&\simeq
\mathcal O_K\langle z_1,\ldots,z_m\rangle_{\mathrm{PD}}^{\wedge},\\
A_{\st}^{(m+1)}/E(u_0)&\simeq
\mathcal O_K\langle w_1,\ldots,w_m\rangle_{\mathrm{PD}}^{\wedge}.
\end{aligned}
\]
Both rings and their quotients by $E(u_0)$ are $p$-torsion-free.
\end{prop}

\begin{proof}
The arguments of \cite[Proposition~4.1.3 and Lemma~5.0.12]{du-liu-prismaticphiGhatmodule} extend to arbitrary $m\geq1$ by taking the $(p,E(u_0))$-completed prismatic envelope of the diagonal (after exactifying the diagonal in the logarithmic case), corresponding to the relations $u_j-u_0=E(u_0)z_j$ and $y_j=E(u_0)w_j$. Since $(u_1-u_0,\ldots,u_m-u_0)$ and $(y_1,\ldots,y_m)$ are regular sequences modulo $(p,E(u_0))$ in $W(k)[[u_0,\ldots,u_m]]$ and $W(k)[[u_0,y_1,\ldots,y_m]]$, \cite[Proposition~3.13]{bhatt-scholze-prismaticcohom} shows that $A^{(m+1)}$ and $A_{\st}^{(m+1)}$ are bounded prisms and $(p,E(u_0))$-completely flat over $\mathfrak S$ via $p_0$ and $i_0$ (see also \cite[Lemmas~2.9--2.10 and Corollary~2.12]{DLMS2} for iterated coproducts in the logarithmic case). Because $\mathfrak S$ is Noetherian, complete flatness implies flatness as in \cite[Proposition~2.2.8(5)]{du-liu-prismaticphiGhatmodule}. Symmetry of the factors gives flatness for every $p_k$ or $i_k$. Since each $p_k$ or $i_k$ is a morphism of prisms, all $E(u_k)$ generate the same prism ideal.

Computing modulo $E(u_0)$ in each variable $z_j$ as in \cite[Proposition~2.2.8(1)--(2)]{du-liu-prismaticphiGhatmodule} identifies the uncompleted quotient with $\mathcal O_K[z_j^n/n!\mid n\geq1,\,1\leq j\leq m]$, and $p$-adic completion gives the isomorphism for $A^{(m+1)}/E(u_0)$. The isomorphism for $A_{\st}^{(m+1)}/E(u_0)$ follows similarly from \cite[Proposition~2.14]{DLMS2}. The quotients modulo $E(u_0)$ are $p$-torsion-free by their divided-power expansions, and $A^{(m+1)}$ and $A_{\st}^{(m+1)}$ are $p$-torsion-free by flatness over $\mathfrak S$.
\end{proof}

\begin{prop}
\label{Prop BK prism covers the final object}
Evaluation on $(\mathfrak S,(E(u)))$ identifies $\operatorname{Vect}^{\varphi}(X_{\Prism})$ with the category of pairs $(\mathfrak M,\varepsilon_{\mathfrak M}^{\cris})$ consisting of a finite free Breuil--Kisin module $\mathfrak M$ and a Frobenius-compatible descent isomorphism
\[
\varepsilon_{\mathfrak M}^{\cris}\colon
A^{(2)}\otimes_{p_1,\mathfrak S}\mathfrak M
\xrightarrow{\sim}A^{(2)}\otimes_{p_0,\mathfrak S}\mathfrak M
\]
satisfying the cocycle condition over $A^{(3)}$. Likewise, evaluation on $(\mathfrak S,(E(u)),\mathbf N)^a$ identifies $\operatorname{Vect}^{\varphi}((X,M_X)_{\Prism})$ with the category of pairs $(\mathfrak M,\varepsilon_{\mathfrak M}^{\st})$ where
\[
\varepsilon_{\mathfrak M}^{\st}\colon
A_{\st}^{(2)}\otimes_{i_1,\mathfrak S}\mathfrak M
\xrightarrow{\sim}A_{\st}^{(2)}\otimes_{i_0,\mathfrak S}\mathfrak M
\]
satisfies the cocycle condition over $A_{\st}^{(3)}$. Both equivalences preserve effectivity and height intervals.
\end{prop}

\begin{proof}
This follows from faithfully flat descent along the Breuil--Kisin cover and its self-products (see \cite[Proposition~4.1.6 and Corollary~4.1.9]{du-liu-prismaticphiGhatmodule} and \cite[Lemma~3.5]{DLMS2}), together with Remark~\ref{rmk def of fh case}(2) for effectivity and height intervals.
\end{proof}

\begin{rem}
\label{rem:parallel-log-nonlog-conventions}
For an object of $\operatorname{Vect}^{\varphi}(X_{\Prism})$, pullback to $(X,M_X)_{\Prism}$ corresponds to base changing $\varepsilon_{\mathfrak M}^{\cris}$ along $j_1\colon A^{(2)}\to A_{\st}^{(2)}$. Throughout, we write $p_k\colon\mathfrak S\to A^{(m+1)}$ and $i_k\colon\mathfrak S\to A_{\st}^{(m+1)}$ ($0\leq k\leq m$) for the $k$-th projection.
\end{rem}

For $\mathbf a=(a_1,\ldots,a_m)\in\mathbf Z^m$, the evaluation maps
\[
\begin{aligned}
\iota_{\mathbf a}\colon A^{(m+1)}&\longrightarrow A_{\inf},
&u_0&\longmapsto\mathbf u,
&u_j&\longmapsto\epsilon^{a_j}\mathbf u,\\
\iota_{\mathbf a}^{\st}\colon A_{\st}^{(m+1)}&\longrightarrow A_{\inf},
&u_0&\longmapsto\mathbf u,
&1+y_j&\longmapsto\epsilon^{a_j}
\end{aligned}
\]
define morphisms of prisms, sending $z_j$ to $(\epsilon^{a_j}-1)\mathbf u/\xi$ and $w_j$ to $(\epsilon^{a_j}-1)/\xi$.

\begin{thm}
\label{thm:multivariable-DL-injectivity}
For every $m\geq1$, the maps
\[
(\iota_{\mathbf{a}})_{\mathbf{a}\in\mathbf{Z}^m}
\colon A^{(m+1)}
\longrightarrow
\prod_{\mathbf{a}\in\mathbf{Z}^m} A_{\inf}, \qquad
(\iota_{\mathbf{a}}^{\st})_{\mathbf{a}\in\mathbf{Z}^m}
\colon A_{\st}^{(m+1)}
\longrightarrow
\prod_{\mathbf{a}\in\mathbf{Z}^m} A_{\inf}
\]
are injective.
\end{thm}

\begin{proof}
Put $R=A_{\st}^{(m+1)}$, and let $f\in R$ satisfy $\iota_{\mathbf a}^{\st}(f)=0$ for all $\mathbf a\in\mathbf Z^m$. By Proposition~\ref{prop:multivariable-coproducts}, the image $\overline f$ of $f$ in $R/E(u_0)\simeq\mathcal O_K\langle w_1,\ldots,w_m\rangle_{\mathrm{PD}}^{\wedge}$ has a unique expansion
\[
\overline f=\sum_{\nu\in\mathbf N^m}b_\nu\frac{w^\nu}{\nu!},
\qquad b_\nu\in\mathcal O_K,\quad b_\nu\longrightarrow0,
\]
where $w^\nu=\prod_{j=1}^m w_j^{\nu_j}$ and $\nu!=\prod_{j=1}^m\nu_j!$. Put $c=\theta((\epsilon-1)/\xi)\in\mathcal O_{C_p}$. Since $\epsilon-1\in\operatorname{Ker}(\theta)=(\xi)$ and $\epsilon^a-1\equiv a(\epsilon-1)\pmod{(\epsilon-1)^2}$ in $A_{\inf}$ for every $a\in\mathbf Z$, we have $\theta\bigl((\epsilon^a-1)/\xi\bigr)=ac$.

For each $\mathbf a=(a_1,\ldots,a_m)\in\mathbf Z^m$, the substitution $w_j\mapsto a_j w$ defines a $p$-adically continuous $\mathcal O_K$-algebra homomorphism $\psi_{\mathbf a}\colon\mathcal O_K\langle w_1,\ldots,w_m\rangle_{\mathrm{PD}}^{\wedge}\to\mathcal O_K\langle w\rangle_{\mathrm{PD}}^{\wedge}$ sending $\overline f$ to
\[
\psi_{\mathbf a}(\overline f)
=\sum_{d\geq0}\Biggl(\sum_{|\nu|=d}\frac{d!}{\nu!}\,b_\nu a^\nu\Biggr)\frac{w^d}{d!},
\]
where $|\nu|=\nu_1+\cdots+\nu_m$ and $\frac{d!}{\nu!}\in\mathbf Z$. Since $\iota_{\mathbf a}^{\st}(E(u_0))=\xi\in\operatorname{Ker}(\theta)$, the composite $\theta\circ\iota_{\mathbf a}^{\st}\colon R\to\mathcal O_{C_p}$ factors through $R/E(u_0)$ as $\overline{\iota}_{\inf}^{\st}\circ\psi_{\mathbf a}$, where $\overline{\iota}_{\inf}^{\st}\colon\mathcal O_K\langle w\rangle_{\mathrm{PD}}^{\wedge}\to\mathcal O_{C_p}$ is induced by $\theta\circ\iota_{\inf}^{\st}$ and sends $w$ to $c$. By the proof of \cite[Proposition~2.4.3]{du-liu-prismaticphiGhatmodule}, the Weierstrass preparation theorem over $\mathcal O_K$ together with the transcendence of $c$ over $K$ (\cite[Lemma~2.4.4]{du-liu-prismaticphiGhatmodule}) implies that $\overline{\iota}_{\inf}^{\st}$ is injective. From $\theta\bigl(\iota_{\mathbf a}^{\st}(f)\bigr)=0$ we therefore obtain $\psi_{\mathbf a}(\overline f)=0$ for all $\mathbf a\in\mathbf Z^m$. Hence for every $d\geq0$, the homogeneous polynomial $\sum_{|\nu|=d}\frac{d!}{\nu!}b_\nu X^\nu\in K[X_1,\ldots,X_m]$ vanishes on $\mathbf Z^m$, so all $b_\nu=0$ and $f=E(u_0)f_1$ for some $f_1\in R$.

Since $\iota_{\mathbf a}^{\st}(E(u_0))=\xi$ is a non-zero-divisor in $A_{\inf}$, we also have $\iota_{\mathbf a}^{\st}(f_1)=0$ for all $\mathbf a\in\mathbf Z^m$. Iterating gives
\[
f\in\bigcap_{n\geq1}E(u_0)^nR
\subseteq\bigcap_{n\geq1}(p,E(u_0))^nR=0.
\]
For $R=A^{(m+1)}$, the expansion in $R/E(u_0)\simeq\mathcal O_K\langle z_1,\ldots,z_m\rangle_{\mathrm{PD}}^{\wedge}$ and $\theta\bigl(\iota_{\mathbf a}(z_j)\bigr)=\pi a_j c$ factor $\theta\circ\iota_{\mathbf a}$ through $z_j\mapsto\pi a_j w$ followed by $\overline{\iota}_{\inf}^{\st}$, and since $\pi\neq0$ in $K$ the same argument gives $f=0$.
\end{proof}

\section{\texorpdfstring{Crystalline and semistable prismatic
$(\varphi,\tau)$-modules}{Crystalline and semistable prismatic
(phi,tau)-modules}}
\label{sec:prismatic-module-categories}
We retain the notation of Section~\ref{sec:prismatic-F-crystals}.
Building on Du--Liu's prismatic $(\varphi,\widehat G)$-modules
\cite{du-liu-prismaticphiGhatmodule}, we introduce crystalline
and semistable prismatic $(\varphi,\tau)$-modules, and we show the data of the $\widehat G$-structure can be recovered by a single $\tau$-semilinear structure with additional conditions.

\subsection{Tate twists}

Put $c_0=E(0)/p\in W(k)^\times$. For $n\in\mathbf Z$, let $\mathfrak S\{n\}=\mathfrak S e_n$ be the Breuil--Kisin module with
\[
\Phi_{\mathfrak S\{n\}}(1\otimes e_n)
=\Bigl(\frac{c_0}{E(u)}\Bigr)^n e_n.
\]
Thus $\mathfrak S\{-1\}$ is effective of height one, and $\mathfrak S\{n\}^{\vee}=\mathfrak S\{-n\}$ (see \cite[Example~2.38]{CGM26}). Put $\mathfrak M\{n\}=\mathfrak M\otimes_{\mathfrak S}\mathfrak S\{n\}$. Then $\mathfrak M$ has height in $[a,b]$ if and only if $\mathfrak M\{a\}$ is effective of height at most $b-a$, so any Breuil--Kisin module becomes effective after a suitable twist.

\subsection{\texorpdfstring{Du--Liu's prismatic
$(\varphi,\widehat G)$-modules}{Du--Liu's prismatic
(phi,Ghat)-modules}}
\label{subsec:prismatic-phi-Ghat-modules}

We keep the embeddings of $A^{(2)}$ and $A_{\st}^{(2)}$ into $A_{\inf}$ from Section~\ref{sec:prismatic-F-crystals}, so that $i_1=\tau\circ i_0$, and in $A_{\st}^{(2)}$ we have
\[
u_1=u_0(1+y),\qquad y=E(u_0)w.
\]
Set
\[
I_+=A_{\st}^{(2)}\cap W(\mathfrak m_{C_p}^{\flat}).
\]
Then $A_{\st}^{(2)}$ is $\widehat G$-stable with $A_{\st}^{(2)}/I_+=W(k)$ (see \cite[Lemma~3.3.1]{du-liu-prismaticphiGhatmodule}). We equip $A^{(2)}$, $A_{\st}^{(2)}$, and finite free modules over them with the $(p,E(u_0))$-adic topology.

\begin{defn}
\label{def:DL-prismatic-phi-Ghat}
\label{def Du-Liu's prismatic phi, G^hat modules}
A \emph{prismatic $(\varphi,\widehat G)$-module of height at most $h$} is an effective finite free Breuil--Kisin module $\mathfrak M$ of height at most $h$, together with a continuous semilinear $\widehat G$-action on
\[
\widehat{\mathfrak M}
:=A_{\st}^{(2)}\otimes_{i_0,\mathfrak S}\mathfrak M
\]
which commutes with Frobenius $\varphi$ and satisfies
\[
1\otimes\mathfrak M\subset\widehat{\mathfrak M}^{H_K},
\qquad
\widehat G\text{ acts trivially on }
\widehat{\mathfrak M}/I_+\widehat{\mathfrak M}.
\]
A morphism is a morphism of underlying Breuil--Kisin modules whose base change to $A_{\st}^{(2)}$ is $\widehat G$-equivariant. We denote this category by $\operatorname{Mod}^{\st,[0,h]}_{\varphi,\widehat G}(\mathfrak S)$.
\end{defn}

\begin{thm}
\label{thm:DL-phi-Ghat-semistable}
\label{Thm cat equivalence to reps}
The functor
\[
\widehat T_{\mathrm{DL}}(\mathfrak M)
=\operatorname{Hom}_{A_{\st}^{(2)},\varphi}
(\widehat{\mathfrak M},A_{\inf})
\]
of \cite{du-liu-prismaticphiGhatmodule}, equipped with the $G_K$-action $(g f)(v)=g\bigl(f(g^{-1}v)\bigr)$ via $G_K\to\widehat G$ on $\widehat{\mathfrak M}$, is an anti-equivalence from $\operatorname{Mod}^{\st,[0,h]}_{\varphi,\widehat G}(\mathfrak S)$ to the category of $G_K$-stable $\mathbf Z_p$-lattices in semistable representations with Hodge--Tate weights in $[-h,0]$. Its dual gives a covariant equivalence with the category of semistable lattices with Hodge--Tate weights in $[0,h]$.
\end{thm}

\begin{proof}
This is \cite[Theorem~3.3.3]{du-liu-prismaticphiGhatmodule}, with Hodge--Tate weights expressed in our sign convention.
\end{proof}

\begin{notn}
\label{notn:covariant-realization}
Let $T_{\mathrm{cov}}$ be the functor $T_{\mathrm{cov}}(\mathfrak M)=\widehat T_{\mathrm{DL}}(\mathfrak M)^\vee$, which identifies $\operatorname{Mod}^{\st,[0,h]}_{\varphi,\widehat G}(\mathfrak S)$ with the category of semistable lattices with Hodge--Tate weights in $[0,h]$, and denote by $\mathfrak M(T)$ the module attached to such a lattice $T$. This agrees with the covariant prismatic realization (see \cite[Remark~4.3.5]{du-liu-prismaticphiGhatmodule}).
\end{notn}

\begin{prop}
\label{prop:DL-crystalline-criterion}
An object $\mathfrak M\in\operatorname{Mod}^{\st,[0,h]}_{\varphi,\widehat G}(\mathfrak S)$ corresponds under $T_{\mathrm{cov}}$ to a crystalline lattice if and only if the action $\tau_{\mathfrak M}$ of $\tau$ on $\widehat{\mathfrak M}$ satisfies
\[
\tau_{\mathfrak M}(1\otimes\mathfrak M)
\subset A^{(2)}\otimes_{p_0,\mathfrak S}\mathfrak M
\subset\widehat{\mathfrak M}.
\]
We denote the resulting full subcategory by $\operatorname{Mod}^{\cris,[0,h]}_{\varphi,\widehat G}(\mathfrak S)$.
\end{prop}

\begin{proof}
This is \cite[Corollary~3.3.4]{du-liu-prismaticphiGhatmodule}.
\end{proof}

\subsection{\texorpdfstring{The $(\varphi,\tau)$-presentation}
{The (phi,tau)-presentation}}

For a finite free $\mathfrak S$-module $\mathfrak M$ and $g\in H_K$, write $G_g$ for the $g$-semilinear automorphism $G_g(a\otimes m)=g(a)\otimes m$ of $\widehat{\mathfrak M}=A_{\st}^{(2)}\otimes_{i_0,\mathfrak S}\mathfrak M$.

\begin{defn}
\label{def:prismatic-phi-tau}
\label{def effective phi Ghat module for crystalline repns}
An \emph{effective semistable prismatic $(\varphi,\tau)$-module of height at most $h$} is a triple $(\mathfrak M,\varphi_{\mathfrak M},\tau_{\mathfrak M})$ satisfying the following conditions:
\begin{enumerate}
\item $(\mathfrak M,\varphi_{\mathfrak M})$ is an effective finite free Breuil--Kisin module of height at most $h$.
\item $\tau_{\mathfrak M}$ is a $\tau$-semilinear automorphism of $\widehat{\mathfrak M}$ commuting with Frobenius.
\item For all $m\in\mathfrak M$ and $g\in H_K$ with $\chi(g)\in\mathbf Z_{>0}$,
\[
G_g\tau_{\mathfrak M}(1\otimes m)
=\tau_{\mathfrak M}^{\chi(g)}(1\otimes m).
\]
\item $(\tau_{\mathfrak M}-1)(1\otimes\mathfrak M)
\subset I_+\widehat{\mathfrak M}$.
\end{enumerate}
We say that $(\mathfrak M,\varphi_{\mathfrak M},\tau_{\mathfrak M})$ is \emph{crystalline} if $\tau_{\mathfrak M}(1\otimes\mathfrak M) \subset A^{(2)}\otimes_{p_0,\mathfrak S}\mathfrak M$. A morphism is a morphism of underlying Breuil--Kisin modules whose base change to $A_{\st}^{(2)}$ commutes with $\tau_{\mathfrak M}$. For $\star\in\{\cris,\st\}$, we denote the resulting category by $\operatorname{Mod}^{\star,[0,h]}_{\varphi,\tau}(\mathfrak S)$ and put
\[
\operatorname{Mod}^{\star,\mathrm{eff}}_{\varphi,\tau}(\mathfrak S)
=\bigcup_{h\geq0}
\operatorname{Mod}^{\star,[0,h]}_{\varphi,\tau}(\mathfrak S).
\]
\end{defn}

\begin{rem}
In Caruso's definition of $(\varphi,\tau)$-modules (see \cite[Definition~1.11]{Caruso_phitau}), one should also require the $\tau$-semilinear map $\tau_M$ to be an automorphism, as in Definition~\ref{def:prismatic-phi-tau}(2). Otherwise the zero map is $\tau$-semilinear, commutes with Frobenius, and satisfies the remaining conditions, yet cannot come from a group action.
\end{rem}

\begin{lem}
\label{lem:uniform-unipotence-coefficients}
The $\widehat G$-action on $A_{\st}^{(2)}$ is continuous, and for $n\geq0$,
\[
(\tau^{p^n}-1)A_{\st}^{(2)}
\subset(p,E(u_0))^nA_{\st}^{(2)}.
\]
\end{lem}

\begin{proof}
Recall from \cite[\S2.3]{du-liu-prismaticphiGhatmodule} the ring $A_{\st,\max}^{(2)}$ over $\mathcal O_{\max}=\bigl(\mathfrak S[E(u)/p]\bigr)^{\wedge}_p$, where $\mathcal O_{\max}$ maps to $A_{\st,\max}^{(2)}$ via $u\mapsto u_0$. Every element of $A_{\st,\max}^{(2)}$ has a unique expansion $\sum_{j\geq0}b_j\gamma_j(w)$ with $\gamma_j(w)=w^j/j!$ and $b_j\in\mathcal O_{\max}$ tending to zero $p$-adically. By the proof of \cite[Proposition~2.3.1(5)]{du-liu-prismaticphiGhatmodule}, which uses \cite[Lemma~2.2.10]{du-liu-prismaticphiGhatmodule},
\begin{equation}
A_{\st}^{(2)}\cap p^nA_{\st,\max}^{(2)}
=(p,E(u_0))^nA_{\st}^{(2)}.
\label{eq:max-induced-topology}
\end{equation}
Indeed, the equality on the uncompleted prismatic envelope passes to the $(p,E(u_0))$-adic completion by approximating modulo $(p,E(u_0))^N$ for $N\geq n$ and using that $(p,E(u_0))^nA_{\st}^{(2)}$ is closed.

For $n\geq1$, set $\alpha_n=(1+y)^{p^n}$ and $b_n=(\alpha_n-1)/y$. Since $y=E(u_0)w\in pA_{\st,\max}^{(2)}$, the relation $b_{n+1}=b_n\sum_{j=0}^{p-1}\alpha_n^j$ gives $b_n\in p^nA_{\st,\max}^{(2)}$. Writing $E(U)-E(V)=(U-V)Q(U,V)$, we obtain
\[
E(\alpha_nu_0)=E(u_0)(1+q_n),\qquad
q_n=u_0wb_nQ(\alpha_nu_0,u_0)
\in p^nA_{\st,\max}^{(2)}.
\]
Thus $\tau^{p^n}$ sends $u_0$ to $\alpha_nu_0$, fixes $y$, and sends $w$ to $(1+q_n)^{-1}w$, so it acts trivially modulo $p^n$ on the topological generators $u_0,y,E(u_0)/p,\gamma_j(w)$ of $A_{\st,\max}^{(2)}$. By \eqref{eq:max-induced-topology}, this gives the stated inclusion for $n\geq1$, while the case $n=0$ is immediate.

For $g\in H_K$, the action fixes $u_0$ and sends $y$ to $(1+y)^{\chi(g)}-1$, defined by its convergent binomial expansion in $A_{\st,\max}^{(2)}$. As $\chi(g)\to 1$, the ratio $((1+y)^{\chi(g)}-1)/y$ tends to $1$ in $A_{\st,\max}^{(2)}$, so $g$ tends to the identity on each $\gamma_j(w)$ and hence on $A_{\st,\max}^{(2)}$. By \eqref{eq:max-induced-topology}, the $H_K$-action on $A_{\st}^{(2)}$ is continuous. Combined with the estimate for $\tau^{p^n}$, this shows that $\widehat G$ acts continuously on $A_{\st}^{(2)}$.
\end{proof}

\begin{prop}
\label{prop:automatic-admissibility}
\label{lem extend-Zp-action continuous}
Let $(\mathfrak M,\varphi_{\mathfrak M},\tau_{\mathfrak M})\in\operatorname{Mod}^{\st,[0,h]}_{\varphi,\tau}(\mathfrak S)$. As $n\to\infty$, $\tau_{\mathfrak M}^{p^n}$ converges uniformly to the identity on $\widehat{\mathfrak M}$, and there is a unique continuous semilinear $\widehat G$-action on $\widehat{\mathfrak M}$ through which $\tau$ acts as $\tau_{\mathfrak M}$ and each $g\in H_K$ acts as $G_g$. This equips $\mathfrak M$ with the structure of an object of $\operatorname{Mod}^{\st,[0,h]}_{\varphi,\widehat G}(\mathfrak S)$.
\end{prop}

\begin{proof}
By Proposition~\ref{prop:multivariable-coproducts}, $\tau(E(u_0))=E(u_1)$ generates $(E(u_0))$ in $A_{\st}^{(2)}$, while $H_K$ fixes $u_0$. Semilinearity therefore implies that $\tau_{\mathfrak M}^{\pm1}$ and all $G_g$ preserve $(p,E(u_0))^r\widehat{\mathfrak M}$ for every $r\geq0$, so $\tau_{\mathfrak M}^{\pm1}$ are continuous, and the $H_K$-action on $\widehat{\mathfrak M}$ is continuous by Lemma~\ref{lem:uniform-unipotence-coefficients}. For $g\in H_K$ with $\chi(g)\in\mathbf Z_{>0}$, the relation $g\tau g^{-1}=\tau^{\chi(g)}$ in $\widehat G$ and Definition~\ref{def:prismatic-phi-tau}(3) give
\[
\begin{aligned}
(G_g\tau_{\mathfrak M}G_g^{-1})(b\otimes m)
&=(g\tau g^{-1})(b)\,G_g\tau_{\mathfrak M}(1\otimes m)\\
&=\tau^{\chi(g)}(b)\,\tau_{\mathfrak M}^{\chi(g)}(1\otimes m)
=\tau_{\mathfrak M}^{\chi(g)}(b\otimes m)
\end{aligned}
\]
for all $b\in A_{\st}^{(2)}$ and $m\in\mathfrak M$.

Since $\chi\colon H_K\to\mathbf Z_p^\times$ has open image, we may choose $g_n\in H_K$ with $g_n\to1$ and $\chi(g_n)=1+p^n$ for all sufficiently large $n$. Then $\tau_{\mathfrak M}^{p^n}=G_{g_n}\tau_{\mathfrak M}G_{g_n}^{-1}\tau_{\mathfrak M}^{-1}$, so for every $x\in\widehat{\mathfrak M}$,
\[
(\tau_{\mathfrak M}^{p^n}-1)x
=G_{g_n}\tau_{\mathfrak M}(G_{g_n}^{-1}-1)
 \tau_{\mathfrak M}^{-1}x+(G_{g_n}-1)x.
\]
Because $(G_{g_n}^{-1}-1)\tau_{\mathfrak M}^{-1}x\to0$ and $(G_{g_n}-1)x\to0$, while $\tau_{\mathfrak M}$ and all $G_{g_n}$ preserve $(p,E(u_0))^r\widehat{\mathfrak M}$, we obtain $(\tau_{\mathfrak M}^{p^n}-1)x\to0$. Choosing an $A_{\st}^{(2)}$-basis $e_1,\ldots,e_d$ of $\widehat{\mathfrak M}$ and writing $x=\sum_i b_i e_i$ gives
\[
(\tau_{\mathfrak M}^{p^n}-1)x
=\sum_i\tau^{p^n}(b_i)(\tau_{\mathfrak M}^{p^n}-1)e_i
 +\sum_i(\tau^{p^n}(b_i)-b_i)e_i.
\]
For each $r$, the first sum lies in $(p,E(u_0))^r\widehat{\mathfrak M}$ for all large $n$ independently of the $b_i$, and so does the second for $n\geq r$ by Lemma~\ref{lem:uniform-unipotence-coefficients}. Hence $\tau_{\mathfrak M}^{p^n}\to 1$ uniformly on $\widehat{\mathfrak M}$.

Writing $\tau_{\mathfrak M}^k-\tau_{\mathfrak M}^\ell=\tau_{\mathfrak M}^\ell(\tau_{\mathfrak M}^{k-\ell}-1)$ and factoring $\tau_{\mathfrak M}^{k-\ell}-1$ through $\tau_{\mathfrak M}^{p^n}-1$ when $p^n\mid(k-\ell)$ shows that, for each $r$, there exists $n$ such that
\[
k\equiv\ell\pmod{p^n}
\quad\Longrightarrow\quad
(\tau_{\mathfrak M}^{k}-\tau_{\mathfrak M}^{\ell})
\widehat{\mathfrak M}
\subseteq(p,E(u_0))^r\widehat{\mathfrak M}
\qquad(k,\ell\in\mathbf Z).
\]
By completeness, $\tau_{\mathfrak M}^a x=\lim_j\tau_{\mathfrak M}^{k_j}x$ is well defined for $a\in\mathbf Z_p$ and any sequence of integers $k_j\to a$, and $(a,x)\mapsto\tau_{\mathfrak M}^a x$ is a continuous semilinear $\mathbf Z_p$-action on $\widehat{\mathfrak M}$. Continuity in the exponent and density of $\{g\in H_K:\chi(g)\in\mathbf Z_{>0}\}$ in $H_K$ extend the conjugation formula to $G_g\tau_{\mathfrak M}^a G_g^{-1}=\tau_{\mathfrak M}^{\chi(g)a}$ for all $g\in H_K$ and $a\in\mathbf Z_p$. Since $(\tau_{\mathfrak M}^a G_g)(\tau_{\mathfrak M}^b G_h)=\tau_{\mathfrak M}^{a+\chi(g)b}G_{gh}$, setting $\rho(\tau^a g)=\tau_{\mathfrak M}^a G_g$ defines a continuous semilinear action of $\widehat G=\overline{\langle\tau\rangle}\rtimes H_K$.

By construction, $H_K$ fixes $1\otimes\mathfrak M$. Since $\mathfrak M$ is effective and $\varphi((p,E(u_0)))\subseteq(p,E(u_0))$, Frobenius on $\widehat{\mathfrak M}$ is continuous, so its commutation with $\tau_{\mathfrak M}$ and all $G_g$ extends to $\widehat G$. By Definition~\ref{def:prismatic-phi-tau}(4) and the trivial $\widehat G$-action on $A_{\st}^{(2)}/I_+=W(k)$, both $\tau_{\mathfrak M}$ and all $G_g$ act trivially modulo $I_+$. Moreover, $I_+\widehat{\mathfrak M}$ is the kernel of $\widehat{\mathfrak M}\to W(k)\otimes_{\mathfrak S}\mathfrak M$, which is continuous because $E(u_0)$ maps to $E(0)\in pW(k)$ and $W(k)\otimes_{\mathfrak S}\mathfrak M$ is $p$-adically separated. Thus $I_+\widehat{\mathfrak M}$ is closed, so $\widehat G$ acts trivially on $\widehat{\mathfrak M}/I_+\widehat{\mathfrak M}$. Uniqueness follows because $\tau$ and $H_K$ topologically generate $\widehat G$.
\end{proof}

\begin{cor}
\label{cor:restriction-tau-equivalence}
\label{Cor Caruso Prop 1.3 on gluing galois action}
\label{Prop Caruso Prop 1.3 on gluing galois action}
Restriction to $\tau$ gives equivalences
\[
\operatorname{Mod}^{\star,[0,h]}_{\varphi,\widehat G}(\mathfrak S)
\xrightarrow{\sim}
\operatorname{Mod}^{\star,[0,h]}_{\varphi,\tau}(\mathfrak S),
\qquad \star\in\{\cris,\st\}.
\]
\end{cor}

\begin{proof}
The condition $1\otimes\mathfrak M\subset\widehat{\mathfrak M}^{H_K}$ forces every $g\in H_K$ to act on $\widehat{\mathfrak M}$ as $G_g$. Proposition~\ref{prop:automatic-admissibility} therefore provides a quasi-inverse to restriction to $\tau$, and the crystalline condition matches Proposition~\ref{prop:DL-crystalline-criterion}.
\end{proof}

\begin{thm}
\label{thm:exact-structures-module-categories}
Evaluation on $(\mathfrak S,(E(u)))$, followed by evaluation of the descent datum at $\tau$, gives equivalences
\[
\begin{aligned}
\operatorname{Vect}^{\varphi}_{[0,h]}((X,M_X)_{\Prism})
&\xrightarrow{\sim}
\operatorname{Mod}^{\st,[0,h]}_{\varphi,\tau}(\mathfrak S),\\
\operatorname{Vect}^{\varphi}_{[0,h]}(X_{\Prism})
&\xrightarrow{\sim}
\operatorname{Mod}^{\cris,[0,h]}_{\varphi,\tau}(\mathfrak S).
\end{aligned}
\]
Under the first equivalence, the descent isomorphism over $A_{\st}^{(2)}$ is
\[
\begin{aligned}
\varepsilon_{\mathfrak M}^{\st}\colon
A_{\st}^{(2)}\otimes_{i_1,\mathfrak S}\mathfrak M
&\xrightarrow{\sim}
A_{\st}^{(2)}\otimes_{i_0,\mathfrak S}\mathfrak M,\\
a\otimes_{i_1} m&\longmapsto a\tau_{\mathfrak M}(1\otimes_{i_0} m),
\end{aligned}
\]
and for a crystalline object it is the base change of a descent isomorphism
\[
\varepsilon_{\mathfrak M}^{\cris}\colon
A^{(2)}\otimes_{p_1,\mathfrak S}\mathfrak M
\xrightarrow{\sim}
A^{(2)}\otimes_{p_0,\mathfrak S}\mathfrak M.
\]
Both reduce to the identity along the diagonal and satisfy the cocycle condition over $A_{\st}^{(3)}$ and $A^{(3)}$, respectively.
\end{thm}

\begin{proof}
By \cite[Corollary~4.1.9 and Proposition~5.0.16]{du-liu-prismaticphiGhatmodule}, effective prismatic $F$-crystals on $(X,M_X)_{\Prism}$ and $X_{\Prism}$ are equivalent to Breuil--Kisin modules with descent data over $A_{\st}^{(\bullet)}$ and $A^{(\bullet)}$, and \cite[\S4.3 and Theorem~5.0.18]{du-liu-prismaticphiGhatmodule} identify these descent data with $(\varphi,\widehat G)$-modules via evaluation at $\tau$. Combining this with Corollary~\ref{cor:restriction-tau-equivalence} yields both equivalences, while compatibility along $A^{(2)}\hookrightarrow A_{\st}^{(2)}$ relates $\varepsilon_{\mathfrak M}^{\cris}$ and $\varepsilon_{\mathfrak M}^{\st}$.
\end{proof}

\subsection{Finite-height objects and exact structures}\label{subsection finite height}

Let $\mathcal O_{\Prism}\{-1\}\in\operatorname{Vect}^{\varphi}_{[0,1]}(X_{\Prism})$ be the object whose value on $(\mathfrak S,(E(u)))$ is $\mathfrak S\{-1\}$. For $n\in\mathbf Z$, put $\mathcal O_{\Prism}\{n\}:=(\mathcal O_{\Prism}\{-1\})^{\otimes(-n)}$, where the zeroth tensor power is $\mathcal O_{\Prism}$ and negative tensor powers are formed via duals in $\operatorname{Vect}^{\varphi}(X_{\Prism})$. We use the same notation for its pullback to $(X,M_X)_{\Prism}$ and for the associated $(\varphi,\widehat G)$- and $(\varphi,\tau)$-modules, whose realizations under $T_{\mathrm{cov}}$ are $\mathbf Z_p(n)$ (see \cite[Example~2.38 and Remark~3.13(3)]{CGM26}). For a crystal $\mathcal E$, put $\mathcal E\{n\}=\mathcal E\otimes\mathcal O_{\Prism}\{n\}$.

For a finite free Breuil--Kisin module $\mathfrak M$, put $\mathfrak M_E:=\mathfrak M[1/E(u)]$. By Proposition~\ref{prop:multivariable-coproducts}, $E(u_0)$ and $E(u_1)$ generate the same ideal in both $A^{(2)}$ and $A_{\st}^{(2)}$, so $p_0,p_1$ and $i_0,i_1$ localize to
\[
\begin{aligned}
  p_0,p_1\colon
  \mathfrak S[1/E(u)]
  &\longrightarrow A^{(2)}[1/E(u_0)],\\
  i_0,i_1\colon
  \mathfrak S[1/E(u)]
  &\longrightarrow A_{\st}^{(2)}[1/E(u_0)].
\end{aligned}
\]
Since $i_0(E(u))=E(u_0)$, localization gives canonical identifications
\[
  \widehat{\mathfrak M}[1/E(u_0)]
  \simeq
  A_{\st}^{(2)}[1/E(u_0)]
  \otimes_{i_0,\mathfrak S}\mathfrak M
  \simeq
  A_{\st}^{(2)}[1/E(u_0)]
  \otimes_{i_0,\mathfrak S[1/E(u)]}\mathfrak M_E,
\]
and similarly with $A^{(2)}$ and $p_0$ in place of $A_{\st}^{(2)}$ and $i_0$. We use these identifications throughout. Moreover, $i_1=\tau\circ i_0$ gives $\tau(E(u_0))=E(u_1)$, so $\tau$ preserves $(E(u_0))$ and extends to $A_{\st}^{(2)}[1/E(u_0)]$. Consequently, $\tau_{\mathfrak M}$ extends uniquely to $\widehat{\mathfrak M}[1/E(u_0)]$. We keep the same notation for these extensions, as well as for the localized descent isomorphisms $\varepsilon_{\mathfrak M}^{\cris}$ and $\varepsilon_{\mathfrak M}^{\st}$.

A finite-height prismatic $(\varphi,\tau)$-module is defined as in Definition~\ref{def:prismatic-phi-tau}, except that in Condition~(1) we allow $\mathfrak M$ to be a finite free Breuil--Kisin module of arbitrary finite height, and in Condition~(2) Frobenius compatibility means
\[
  \varphi\circ\tau_{\mathfrak M}
  =
  \tau_{\mathfrak M}\circ\varphi
  \colon
  \widehat{\mathfrak M}
  \longrightarrow
  \widehat{\mathfrak M}[1/E(u_0)],
\]
where $\varphi$ is induced by $\varphi_{\mathfrak M}$. The remaining conditions, including the crystalline condition and the notion of morphism, are unchanged, and finite-height $(\varphi,\widehat G)$-modules are defined analogously from Definition~\ref{def:DL-prismatic-phi-Ghat}.

For $\star\in\{\cris,\st\}$ and $a\leq b$, put
\[
\operatorname{Mod}^{\star,[a,b]}_{\varphi,\tau}(\mathfrak S)
=\{\mathfrak M:\mathfrak M\{a\}\in
\operatorname{Mod}^{\star,[0,b-a]}_{\varphi,\tau}(\mathfrak S)\},
\qquad
\operatorname{Mod}^{\star,\mathrm{fh}}_{\varphi,\tau}(\mathfrak S)
=\bigcup_{a\leq b}
\operatorname{Mod}^{\star,[a,b]}_{\varphi,\tau}(\mathfrak S),
\]
and similarly for $(\varphi,\widehat G)$-modules and prismatic $F$-crystals. Since the Frobenius of any prismatic $F$-crystal and its inverse have bounded $E(u)$-denominators on $(\mathfrak S,(E(u)))$, every object of $\operatorname{Vect}^{\varphi}(Y_{\Prism})$ has finite height, so
\[
\operatorname{Vect}^{\varphi,\mathrm{fh}}(Y_{\Prism})
=\operatorname{Vect}^{\varphi}(Y_{\Prism})
\qquad\text{for } Y_{\Prism}\in\{X_{\Prism},(X,M_X)_{\Prism}\}.
\]

\begin{cor}
\label{cor:restriction-tau-fh}
Restriction to $\tau$ gives height-preserving equivalences
\[
\operatorname{Mod}^{\star,\mathrm{fh}}_{\varphi,\widehat G}(\mathfrak S)
\xrightarrow{\sim}
\operatorname{Mod}^{\star,\mathrm{fh}}_{\varphi,\tau}(\mathfrak S),
\qquad \star\in\{\cris,\st\},
\]
and every object on the right carries the unique continuous $\widehat G$-action of Proposition~\ref{prop:automatic-admissibility}.
\end{cor}

\begin{proof}
Twisting by $\mathfrak S\{n\}$ for suitable $n$ reduces to the effective case in Corollary~\ref{cor:restriction-tau-equivalence}.
\end{proof}

\begin{cor}
\label{cor:crystal-tau-fh}
The equivalences of Theorem~\ref{thm:exact-structures-module-categories} extend to height-preserving equivalences
\[
\begin{aligned}
\operatorname{Vect}^{\varphi}((X,M_X)_{\Prism})
&\xrightarrow{\sim}
\operatorname{Mod}^{\st,\mathrm{fh}}_{\varphi,\tau}(\mathfrak S),\\
\operatorname{Vect}^{\varphi}(X_{\Prism})
&\xrightarrow{\sim}
\operatorname{Mod}^{\cris,\mathrm{fh}}_{\varphi,\tau}(\mathfrak S),
\end{aligned}
\]
with descent isomorphisms given by the same formulas.
\end{cor}

\begin{proof}
Twisting reduces to Theorem~\ref{thm:exact-structures-module-categories}, and the result is independent of the choice of twist since evaluation and descent commute with tensor products and duals.
\end{proof}

We equip each module category above with the exact structure in which a sequence is exact if and only if its underlying sequence of finite free $\mathfrak S$-modules is exact (and hence split over $\mathfrak S$). Pushouts of admissible monomorphisms and pullbacks of admissible epimorphisms are formed on the underlying $\mathfrak S$-modules with the induced Frobenius and $\widehat G$- or $\tau$-actions. After twisting to $[0,h]$, the height bound amounts to the existence of $V=E(u)^h\Phi_{\mathfrak M}^{-1}$ satisfying
\[
\Phi_{\mathfrak M}V=E(u)^h\operatorname{id}_{\mathfrak M},
\qquad
V\Phi_{\mathfrak M}=E(u)^h\operatorname{id}_{\varphi^*\mathfrak M},
\]
so pushouts and pullbacks preserve every height interval, verifying the axioms of an exact category. On $\operatorname{Vect}^{\varphi}(X_{\Prism})$ and $\operatorname{Vect}^{\varphi}((X,M_X)_{\Prism})$, a sequence is exact if it is exact on every prism, which by faithfully flat descent along the Breuil--Kisin cover is equivalent to exactness on $(\mathfrak S,(E(u)))$. Thus the equivalences of Theorem~\ref{thm:exact-structures-module-categories} and Corollary~\ref{cor:crystal-tau-fh} are exact.

Finally, $T_{\mathrm{cov}}$ extends to finite-height objects via
\[
T_{\mathrm{cov}}(\mathfrak M)
=T_{\mathrm{cov}}(\mathfrak M\{n\})(-n)
\]
for any $n$ such that $\mathfrak M\{n\}$ is effective. This is independent of $n$ and induces equivalences of categories
\[
\operatorname{Mod}^{\star,\mathrm{fh}}_{\varphi,\tau}(\mathfrak S)
\xrightarrow{\sim}\operatorname{Rep}^{\star}_{\mathbf Z_p}(G_K),
\qquad \star\in\{\cris,\st\},
\]
with the category of $G_K$-stable $\mathbf Z_p$-lattices in crystalline, respectively semistable, representations. We still write $\mathfrak M(T)$ for the module attached to $T$, which has height in $[a,b]$ whenever $T$ has Hodge--Tate weights in $[a,b]$. Note, however, that $\mathfrak M(-)$ need not preserve arbitrary short exact sequences of $G_K$-lattices (see \cite[\S3.1]{CGM26}).

\subsection{Formal coordinates on the logarithmic diagonal}\label{subsection Formal coordinates on the logarithmic diagonal}

To analyze degree-one \v{C}ech--Alexander cocycles in Proposition~\ref{Prop ker d^1_prism}, we embed the localizations of $A_{\st}^{(2)}$ and $A_{\st}^{(3)}$ into formal power series rings. Let
$\partial_{01},\partial_{02},\partial_{12}
\colon A_{\st}^{(2)}
\longrightarrow A_{\st}^{(3)}$
be the cosimplicial maps induced by the pairs of indices $(0,1)$, $(0,2)$, and $(1,2)$. Retain $\mathcal O_{\max}=\bigl(\mathfrak S[E(u)/p]\bigr)^{\wedge}_p$ from the proof of Lemma~\ref{lem:uniform-unipotence-coefficients}, and set $B=\mathcal O_{\max}[1/p,1/E(u)]$.

\begin{lem}
\label{lem:faithful-log-Taylor}
The substitutions
\[
\begin{aligned}
\eta_1\colon
A_{\st}^{(2)}[1/p,1/E(u_0)]&\longrightarrow B[[t]], & u_0&\mapsto u,\quad 1+y\mapsto e^t,\\
\eta_2\colon
A_{\st}^{(3)}[1/p,1/E(u_0)]&\longrightarrow B[[s,t]], & u_0&\mapsto u,\quad 1+y_1\mapsto e^s,\quad 1+y_2\mapsto e^{s+t}
\end{aligned}
\]
define injective ring homomorphisms, and $u\mapsto e^su$ defines a ring homomorphism $\rho_s\colon B\to B[[s]]$. For $f\in A_{\st}^{(2)}[1/p,1/E(u_0)]$ with $F_f(u,t)=\eta_1(f)$, we have
\[
\begin{aligned}
\eta_2(\partial_{01}f)&=F_f(u,s),\\
\eta_2(\partial_{02}f)&=F_f(u,s+t),\\
\eta_2(\partial_{12}f)&=F_f(e^su,t)
=\rho_s(F_f(u,t)),
\end{aligned}
\]
where $\rho_s$ acts on coefficients in $B$ and fixes $t$.
\end{lem}

\begin{proof}
Let $D=u\,d/du$ on $\mathfrak S[E(u)/p]$. Since $D(\mathfrak S)\subset\mathfrak S$ and $D(E(u)/p)=uE'(u)/p$, we have $D(\mathfrak S[E(u)/p])\subset p^{-1}\mathfrak S[E(u)/p]$, so $D$ extends continuously to $\mathcal O_{\max}\to p^{-1}\mathcal O_{\max}$ and hence to $B\to B$. By the Leibniz rule,
\[
\rho_s(b)=\sum_{j\geq0}\frac{D^j(b)}{j!}s^j
\]
defines a ring homomorphism $\rho_s\colon B\to B[[s]]$ sending $u$ to $e^su$. In particular, $\rho_s(E(u))=E(e^su)$ has invertible constant term $E(u)\in B^\times$, so $\rho_s(E(u)^{-1})=E(e^su)^{-1}$.

Put $\mathscr P_2=(\mathcal O_{\max}\langle w_1,w_2\rangle_{\mathrm{PD}})^{\wedge}_p$. By the proof of \cite[Proposition~2.3.1(2)]{du-liu-prismaticphiGhatmodule}, there is an embedding $A_{\st}^{(2)}\hookrightarrow A_{\st,\max}^{(2)}=(\mathcal O_{\max}\langle w\rangle_{\mathrm{PD}})^{\wedge}_p$ sending $u_0$ to $u$ and $y/E(u_0)$ to $w$. Since $i_0\colon\mathfrak S\to A_{\st}^{(2)}$ is flat and strict on log structures, \cite[Remark~2.4]{DLMS2} identifies $A_{\st}^{(3)}$ with the $(p,E(u_0))$-completed tensor product $A_{\st}^{(2)}\widehat{\otimes}_{i_0,\mathfrak S,i_0}A_{\st}^{(2)}$. Because $E(u)\in p\mathcal O_{\max}$, the two copies of $A_{\st}^{(2)}\to(\mathcal O_{\max}\langle w_j\rangle_{\mathrm{PD}})^{\wedge}_p$ ($j=1,2$) induce a continuous ring homomorphism
\[
A_{\st}^{(3)}\longrightarrow\mathscr P_2,
\qquad u_0\mapsto u,\quad y_j/E(u_0)\mapsto w_j.
\]

To see that $A_{\st}^{(3)}\to\mathscr P_2$ is injective, set $A_{\max}=\widehat{A_{\inf}[\xi/p]}_p$, and for $a\in\mathbf Z$ put
\[
q_a=\frac{\epsilon^a-1}{\xi}
=\frac{\epsilon^a-1}{\epsilon-1}\frac{\epsilon-1}{\xi}.
\]
By \cite[Proposition~2.4.3]{du-liu-prismaticphiGhatmodule}, all divided powers of $(\epsilon-1)/\xi$ lie in $A_{\max}$, and hence so do those of $q_a$ since $(\epsilon^a-1)/(\epsilon-1)\in A_{\inf}$. For $(a_1,a_2)\in\mathbf Z^2$, sending $u\mapsto\mathbf u$ and $(w_1,w_2)\mapsto(q_{a_1},q_{a_2})$ defines a continuous homomorphism $\mathscr P_2\to A_{\max}$. We check that its composite with $A_{\st}^{(3)}\to\mathscr P_2$ agrees with the corresponding evaluation $A_{\st}^{(3)}\to A_{\inf}\hookrightarrow A_{\max}$ of Section~\ref{sec:prismatic-F-crystals}. Put $c=\varphi(E(u))/p\in\mathcal O_{\max}^\times$. The Frobenius on $A_{\st,\max}^{(2)}$ is compatible with $\varphi$ on $\mathcal O_{\max}$ and sends $w$ to $\bigl((1+E(u)w)^p-1\bigr)/(pc)$, which under $u\mapsto\mathbf u$ and $w\mapsto q_a$ maps to
\[
\frac{(1+\xi q_a)^p-1}{\varphi(\xi)}
=\frac{\epsilon^{ap}-1}{\varphi(\xi)}=\varphi(q_a).
\]
Since $A_{\max}$ is $p$-torsion-free and $p\delta(x)=\varphi(x)-x^p$ in $A_{\st}^{(3)}$, agreement on $u_0$ and $w_1,w_2$ extends to all $\delta$-iterates and hence to $A_{\st}^{(3)}$. Any element in the kernel of $A_{\st}^{(3)}\to\mathscr P_2$ therefore vanishes under all such evaluations and is zero by Theorem~\ref{thm:multivariable-DL-injectivity}.

Every element of $\mathscr P_2$ has a unique expansion
\[
\sum_{a,b\geq0}b_{a,b}\gamma_a(w_1)\gamma_b(w_2),
\qquad b_{a,b}\in\mathcal O_{\max},\quad b_{a,b}\to0
\quad(a+b\to\infty).
\]
Substituting $w_1=(e^s-1)/E(u)$ and $w_2=(e^{s+t}-1)/E(u)$ contributes only finitely many terms to each degree in $(s,t)$ and hence defines a ring homomorphism $\mathscr P_2\to B[[s,t]]$. This map is injective because the lowest nonzero homogeneous part of a nonzero series is transformed by the invertible linear substitution $(w_1,w_2)\mapsto\bigl(s/E(u),(s+t)/E(u)\bigr)$ and cannot vanish. Inverting $p$ and $E(u_0)$ yields the injective map $\eta_2$, and the same argument with $A_{\st}^{(2)}\hookrightarrow A_{\st,\max}^{(2)}$ and $w\mapsto(e^t-1)/E(u)$ gives $\eta_1$.

Finally, on $A_{\st}^{(2)}$ the maps $\partial_{01},\partial_{02},\partial_{12}$ act by
\[
\begin{array}{c|cc}
&u_0&1+y\\ \hline
\partial_{01}&u_0&1+y_1\\
\partial_{02}&u_0&1+y_2\\
\partial_{12}&u_0(1+y_1)&(1+y_2)/(1+y_1).
\end{array}
\]
Writing $E(U)-E(V)=(U-V)Q(U,V)$ gives $E(u_1)=E(u_0)\bigl(1+u_0w_1Q(u_1,u_0)\bigr)$, where the second factor is a unit in $A_{\st}^{(3)}$ since $u_0$ is topologically nilpotent for the $(p,E(u_0))$-adic topology. Thus $\partial_{01},\partial_{02},\partial_{12}$ localize to $A_{\st}^{(2)}[1/p,1/E(u_0)]\to A_{\st}^{(3)}[1/p,1/E(u_0)]$. Since $\eta_2$ sends $u_0$, $u_0(1+y_1)$, and $u_0(1+y_2)$ to $u$, $e^su$, and $e^{s+t}u$, the formulas for $\eta_2(\partial_{ij}f)$ follow.
\end{proof}

\section{Prismatic cohomology via \texorpdfstring{$\tau$}{tau}-complexes}
\label{sec:prismatic-cohomology-tau}

In this section, we construct explicit two-term
$\tau$-complexes computing the cohomology of prismatic crystals on the absolute (logarithmic) prismatic
sites (see Theorem~\ref{prop compute prismatic coho}).
The key observation is to identify the subscript-zero condition with the degree-one \v{C}ech--Alexander cocycle condition
for the corresponding Breuil--Kisin cover
(see Proposition~\ref{Prop ker d^1_prism}). These two-term complexes is a key input for the prismatic Herr complexes of the next section.

\subsection{The \v{C}ech--Alexander complexes}
\label{subsection cech nerve}

Let $\mathfrak M_{\Prism}$ be a prismatic $F$-crystal on $X_{\Prism}$ or $(X,M_X)_{\Prism}$, and put $\mathfrak M=\mathfrak M_{\Prism}(\mathfrak S,(E(u)))$, where on $(X,M_X)_{\Prism}$ the Breuil--Kisin prism is endowed with its canonical log structure. Unless stated otherwise, tensor products of $\mathfrak M$ with $A^{(n)}$ or $A_{\st}^{(n)}$ are taken via $p_0$ or $i_0$, so the prism ideal is generated by $E(u_0)$ (with $u_0=u$ in degree $0$). Consider the \v{C}ech--Alexander complexes
\begin{align*}
C_{\Prism}^{\bullet}&\colon \mathfrak M\xrightarrow{d_{\Prism}^{0}}
A^{(2)}\otimes_{p_0,\mathfrak S}\mathfrak M
\xrightarrow{d_{\Prism}^{1}}
A^{(3)}\otimes_{p_0,\mathfrak S}\mathfrak M\longrightarrow\cdots,\\
C_{\Prism,\log}^{\bullet}&\colon
\mathfrak M\xrightarrow{d_{\Prism,\log}^{0}}
A_{\st}^{(2)}\otimes_{i_0,\mathfrak S}\mathfrak M
\xrightarrow{d_{\Prism,\log}^{1}}
A_{\st}^{(3)}\otimes_{i_0,\mathfrak S}\mathfrak M\longrightarrow\cdots
\end{align*}
of $\mathfrak M_{\Prism}$ on $X_{\Prism}$ and $(X,M_X)_{\Prism}$, respectively. Write
\[
\begin{aligned}
\varepsilon_{\mathfrak M}^{\cris}&\colon
A^{(2)}\otimes_{p_1,\mathfrak S}\mathfrak M
\xrightarrow{\;\sim\;}A^{(2)}\otimes_{p_0,\mathfrak S}\mathfrak M,\\
\varepsilon_{\mathfrak M}^{\st}&\colon
A_{\st}^{(2)}\otimes_{i_1,\mathfrak S}\mathfrak M
\xrightarrow{\;\sim\;}A_{\st}^{(2)}\otimes_{i_0,\mathfrak S}\mathfrak M
\end{aligned}
\]
for the descent isomorphisms. When $\mathfrak M_{\Prism}$ is defined on $X_{\Prism}$, $\varepsilon_{\mathfrak M}^{\st}$ is the base change of $\varepsilon_{\mathfrak M}^{\cris}$. In both cases,
\[
d^0(m)=\varepsilon_{\mathfrak M}(1\otimes m)-1\otimes m,
\qquad
d^1(z)=\varepsilon_{01}(\partial_{12}z)-\partial_{02}z+\partial_{01}z,
\]
where $\partial_{ab}$ is induced by the $(a,b)$-projection and $\varepsilon_{01}=\partial_{01}^*\varepsilon_{\mathfrak M}$.

Both complexes compute derived global sections. Indeed, the Breuil--Kisin prisms cover the final objects of $X_{\Prism}$ and $(X,M_X)_{\Prism}$ by \cite[Example~2.6]{bhatt-scholze-prismaticFcrystal} and \cite[Lemma~2.8]{DLMS2}, and on each term of the \v{C}ech nerve a finite locally free crystal has vanishing higher cohomology by complete faithfully flat descent: the argument of \cite[Corollary~3.12]{bhatt-scholze-prismaticcohom} applies to both sites since coverings are strict on log structures and completely faithfully flat on underlying rings. The \v{C}ech-to-derived spectral sequence therefore yields
\[
R\Gamma(X_{\Prism},\mathfrak M_{\Prism})\simeq C_{\Prism}^{\bullet}, \qquad
R\Gamma((X,M_X)_{\Prism},\mathfrak M_{\Prism})\simeq C_{\Prism,\log}^{\bullet}.
\]

\subsection{Continuous cochains and the Kedlaya--Liu comparison}
\label{section The prismatic complex}

Recall from Corollary~\ref{cor:crystal-tau-fh} the continuous $\widehat G$-action on $\widehat{\mathfrak M}$, under which $H_K$ acts through $A_{\st}^{(2)}$ and fixes $1\otimes\mathfrak M$. We equip $A_{\st}^{(2)}$ and $\widehat{\mathfrak M}$ with the $(p,E(u_0))$-adic topology.

\begin{lem}
\label{lem:HK-invariants-prismatic}
We have
\[
\bigl(A_{\st}^{(2)}\bigr)^{H_K}=i_0(\mathfrak S),
\qquad
\widehat{\mathfrak M}^{H_K}=1\otimes\mathfrak M.
\]
\end{lem}

\begin{proof}
Let $\operatorname{aug}\colon A_{\st}^{(2)}\to\mathfrak S$ be the augmentation map sending $u_0\mapsto u$ and $y\mapsto0$, which kills $w=y/E(u_0)$ and all its $\delta$-iterates. Given $f\in(A_{\st}^{(2)})^{H_K}$, put
\[
F_f(u,t)=\eta_1(f)\in B[[t]],
\qquad B=\mathcal O_{\max}[1/p,1/E(u)],
\]
as in Lemma~\ref{lem:faithful-log-Taylor}. For $q=\chi(\gamma)>1$, the element $\gamma$ fixes $u_0$ and sends $1+y$ to $(1+y)^q$, so $F_f(u,qt)=F_f(u,t)$. Since $q^n-1\in B^\times$ for every $n\geq1$, the series $F_f(u,t)$ is constant in $t$. Let $\kappa\colon\mathfrak S\to B$ be the natural map. Specializing at $t=0$ agrees with $\operatorname{aug}$, whereas $\eta_1(i_0(a))$ is the constant series $\kappa(a)$ for $a\in\mathfrak S$. Viewing the middle terms below as constant series in $B[[t]]$, we obtain
\[
\eta_1(f)
=F_f(u,0)
=\kappa(\operatorname{aug}(f))
=\eta_1\bigl(i_0(\operatorname{aug}(f))\bigr).
\]
By injectivity of $\eta_1$ on $A_{\st}^{(2)}$, it follows that $f=i_0(\operatorname{aug}(f))\in i_0(\mathfrak S)$. Since the reverse inclusion is clear, this proves the first equality, and choosing an $\mathfrak S$-basis of $\mathfrak M$ gives the second.
\end{proof}

For a topological module $V$ with continuous $\widehat G$-action, let $C^r_{\mathrm{cont}}(\widehat G,V)$ denote the group of continuous maps $\widehat G^r\to V$, equipped with the differential
\[
\begin{aligned}
(df)(g_1,\ldots,g_{r+1})
={}&g_1f(g_2,\ldots,g_{r+1})+\sum_{j=1}^{r}(-1)^j
f(g_1,\ldots,g_jg_{j+1},\ldots,g_{r+1})\\
&+(-1)^{r+1}f(g_1,\ldots,g_r),
\end{aligned}
\]
and write $H^r_{\mathrm{cont}}(\widehat G,V)$ for its cohomology (and similarly for $H_K$). For $r\geq0$, set
\[
K^r(V)=C^r_{\mathrm{cont}}(\widehat G,V)^{H_K^{r+1}},
\]
where $(h_0,\ldots,h_r)\in H_K^{r+1}$ acts via
\[
\bigl((h_0,\ldots,h_r)f\bigr)(g_1,\ldots,g_r)
=h_0f(h_0^{-1}g_1h_1,\ldots,h_{r-1}^{-1}g_rh_r).
\]
A direct check shows that $d$ preserves $K^\bullet(V)$, making $K^\bullet(V)$ a subcomplex of $C^\bullet_{\mathrm{cont}}(\widehat G,V)$. This is the relative continuous cochain complex of \cite[Examples~5.5.6--5.5.7]{KL16} (in \cite[Remark~4.2.10]{du-liu-prismaticphiGhatmodule}, the two $H_K$-factors are ordered oppositely).

When $V=\widehat{\mathfrak M}$, Lemma~\ref{lem:HK-invariants-prismatic} identifies $K^0(\widehat{\mathfrak M})$ with $\mathfrak M$. Writing $d_{\mathrm{KL}}$ for the differential on $K^\bullet(\widehat{\mathfrak M})$, we have in low degrees
\[
(d_{\mathrm{KL}}^0m)(g)=g(1\otimes m)-1\otimes m,\qquad
(d_{\mathrm{KL}}^1f)(g_1,g_2)=g_1f(g_2)-f(g_1g_2)+f(g_1).
\]

We next compare prismatic cochains with $K^\bullet(\widehat{\mathfrak M})$ over $A_{\st}^{(2)}$. Every $g\in\widehat G$ admits a unique decomposition $g=\tau^{a(g)}h$ with $a(g)\in\mathbf Z_p$ and $h\in H_K$, and $g\mapsto a(g)$ is continuous. Given $\mathbf g=(g_1,\ldots,g_r)\in\widehat G^r$, put $b_j=a(g_1\cdots g_j)$ and define a ring homomorphism
\[
e_{\mathbf g}^{\st}\colon A_{\st}^{(r+1)}\longrightarrow A_{\st}^{(2)},
\qquad
u_0\longmapsto u_0,\quad 1+y_j\longmapsto(1+y)^{b_j}.
\]
Indeed, $w_j=y_j/E(u_0)$ maps to the $(p,E(u_0))$-adically convergent series
\[
w_j\longmapsto w\sum_{n\geq1}\binom{b_j}{n}y^{n-1},
\qquad w=\frac{y}{E(u_0)},
\]
and since these assignments commute with Frobenius on the base ring, they extend uniquely to the $\delta$-envelope and its $(p,E(u_0))$-adic completion. Moreover, $e_{\mathbf g}^{\st}(x)$ depends continuously on $\mathbf g$: the binomial series above converge uniformly in $\mathbf g$, the $\delta$-operation is continuous with $\delta(I^N)\subset I^{N-1}$ for $I=(p,E(u_0))$ and $N\geq1$, and each $e_{\mathbf g}^{\st}$ preserves $I^N$, so continuity passes to the completion by $I$-adic approximation. Setting $e_{\mathbf g}^{\cris}=e_{\mathbf g}^{\st}\circ j_r$ in the non-logarithmic case, we obtain maps
\[
\begin{aligned}
\iota_r^{\st}\colon
A_{\st}^{(r+1)}&\longrightarrow K^r(A_{\st}^{(2)}),
& x&\longmapsto\bigl[\mathbf g\longmapsto e_{\mathbf g}^{\st}(x)\bigr],\\
\iota_r^{\cris}\colon
A^{(r+1)}&\longrightarrow K^r(A_{\st}^{(2)}),
& x&\longmapsto\bigl[\mathbf g\longmapsto e_{\mathbf g}^{\cris}(x)\bigr].
\end{aligned}
\]
To check that $\iota_r^{\st}(x)$ is $H_K^{r+1}$-invariant, note that replacing $\mathbf g$ by $\mathbf g'=(h_0^{-1}g_1h_1,\ldots,h_{r-1}^{-1}g_rh_r)$ sends $g_1\cdots g_j$ to $h_0^{-1}(g_1\cdots g_j)h_j$ and thus replaces $b_j$ by $\chi(h_0)^{-1}b_j$. Since $h_0$ fixes $u_0$ and sends $1+y$ to $(1+y)^{\chi(h_0)}$, we have $h_0(e_{\mathbf g'}^{\st}(x))=e_{\mathbf g}^{\st}(x)$.

In degree $1$, the maps $q_0,q_1\colon\mathfrak S\to K^1(A_{\st}^{(2)})$ given by $q_0(s)(g)=i_0(s)$ and $q_1(s)(g)=g(i_0(s))$ fit into the diagram
\[
\begin{tikzcd}[column sep=large,row sep=large]
\mathfrak S
 \arrow[r,shift left=.7ex,"p_0"]
 \arrow[r,shift right=.7ex,"p_1"']
 \arrow[d,shift left=.7ex,"q_0"]
 \arrow[d,shift right=.7ex,"q_1"']
&A^{(2)}\arrow[d,"\iota_1^{\cris}"]
 \arrow[r,hook,"j_1"]
&A_{\st}^{(2)}\arrow[dl,"\iota_1^{\st}"']
 \arrow[d,equal]\\
K^1(A_{\st}^{(2)})\arrow[r,equal]
&K^1(A_{\st}^{(2)})\arrow[r,"\operatorname{ev}_{\tau}"']
&A_{\st}^{(2)},
\end{tikzcd}
\]
in which the left square commutes for each pair of parallel arrows ($p_0$ with $q_0$, and $p_1$ with $q_1$). In particular,
\[
\operatorname{ev}_{\tau}\circ\iota_1^{\st}=\operatorname{id}_{A_{\st}^{(2)}},
\qquad
\operatorname{ev}_{\tau}\circ\iota_1^{\cris}=j_1.
\]

\begin{prop}
\label{lemma key injectivity}
\label{prop:integral-cochain-comparison}
The maps $\iota_r^{\st}$ and $\iota_r^{\cris}$ are injective for all $r\geq0$. They induce termwise injective morphisms of complexes
\[
C_{\Prism,\log}^{\bullet}\longrightarrow K^{\bullet}(\widehat{\mathfrak M}),
\qquad
C_{\Prism}^{\bullet}\longrightarrow K^{\bullet}(\widehat{\mathfrak M})
\]
in the logarithmic and crystalline cases, respectively, both equal to the identity in degree $0$.
\end{prop}

\begin{proof}
Composing $\iota_r^{\st}$ or $\iota_r^{\cris}$ with evaluation at $(\tau^{a_1},\tau^{a_2-a_1},\ldots,\tau^{a_r-a_{r-1}})$ for $(a_1,\ldots,a_r)\in\mathbf Z^r$ and $A_{\st}^{(2)}\hookrightarrow A_{\inf}$ gives the map of Theorem~\ref{thm:multivariable-DL-injectivity}, which is injective. The statement in degree $0$ is Lemma~\ref{lem:HK-invariants-prismatic}.

On modules, the map sends $a\otimes m$ to $\mathbf g\mapsto e_{\mathbf g}(a)(1\otimes m)$ (where $e_{\mathbf g}=e_{\mathbf g}^{\st}$ or $e_{\mathbf g}^{\cris}$), which is injective since $\mathfrak M$ is finite free over $\mathfrak S$. Compatibility with all coface maps except $\partial_0$ follows directly from the definition of $e_{\mathbf g}$. For $\partial_0$, $\varepsilon_{\mathfrak M}^{\st}$ evaluates at $g_1$ to the action of $g_1$ on $\widehat{\mathfrak M}$: by Corollary~\ref{cor:crystal-tau-fh}, the $\widehat G$-action on $\widehat{\mathfrak M}$ is the unique continuous action fixing $1\otimes\mathfrak M$ on $H_K$ and acting as $\tau_{\mathfrak M}$ at $\tau$. Thus $\partial_0$ matches the term $g_1f(g_2,\ldots,g_{r+1})$ of $d_{\mathrm{KL}}$.
\end{proof}

\subsection{Relative cocycles and the subscript-zero condition}

\begin{defn}
\label{def of sub-zero}
Define the $\mathbf Z_p$-submodule $N_{\st}(\mathfrak M)\subset\widehat{\mathfrak M}$ by
\[
N_{\st}(\mathfrak M)
=\bigl(A_{\st}^{(2)}\otimes_{i_0,\mathfrak S}\mathfrak M\bigr)_0
:=\big\{z\in\widehat{\mathfrak M}\;\big|\;
(g\otimes1)z=\sum_{j=0}^{\chi(g)-1}\tau_{\mathfrak M}^{j}(z)\text{ for all }g\in H_K\text{ with }\chi(g)\in\mathbf Z_{>0}\big\}.
\]
When $\mathfrak M_{\Prism}$ is defined on $X_{\Prism}$, set
\[
N_{\cris}(\mathfrak M)
=\bigl(A^{(2)}\otimes_{p_0,\mathfrak S}\mathfrak M\bigr)_0
:=N_{\st}(\mathfrak M)\cap
\bigl(A^{(2)}\otimes_{p_0,\mathfrak S}\mathfrak M\bigr),
\]
where the intersection is formed inside $\widehat{\mathfrak M}$.
\end{defn}

\begin{prop}
\label{Lem Ker^1_KL}
\label{prop:tau-evaluation-KL}
Evaluation at $\tau$ induces a natural isomorphism
\[
\operatorname{ev}_{\tau}\colon
\Ker d_{\mathrm{KL}}^1\xrightarrow{\;\sim\;}N_{\st}(\mathfrak M).
\]
\end{prop}

\begin{proof}
An element $f\in K^1(\widehat{\mathfrak M})$ lies in $\Ker d_{\mathrm{KL}}^1$ if and only if $f\colon\widehat G\to\widehat{\mathfrak M}$ is a continuous $1$-cocycle vanishing on $H_K$: any $1$-cocycle satisfies $f(1)=0$, so $H_K^2$-invariance implies $f|_{H_K}=0$, and conversely $f|_{H_K}=0$ together with $f(xy)=f(x)+xf(y)$ gives $f(hgk)=hf(g)$ for all $h,k\in H_K$. Iterating $f(xy)=f(x)+xf(y)$ gives
\[
f(\tau^n)=\sum_{j=0}^{n-1}\tau_{\mathfrak M}^j\bigl(f(\tau)\bigr)
\qquad(n\geq1).
\]
For $g\in H_K$ with $q=\chi(g)\in\mathbf Z_{>0}$, the relation $g\tau=\tau^q g$ and $f|_{H_K}=0$ give
\[
(g\otimes1)f(\tau)
=f(g\tau)
=f(\tau^q g)
=f(\tau^q)
=\sum_{j=0}^{q-1}\tau_{\mathfrak M}^j\bigl(f(\tau)\bigr),
\]
so $f(\tau)\in N_{\st}(\mathfrak M)$. Since $\{\tau^nh:n\in\mathbf Z,\ h\in H_K\}$ is dense in $\widehat G$, $f$ is uniquely determined by $f(\tau)$.

Conversely, for $z\in N_{\st}(\mathfrak M)$ and $n\in\mathbf Z$, put
\[
S_n(z)=
\begin{cases}
\sum_{j=0}^{n-1}\tau_{\mathfrak M}^j(z),&n>0,\\
0,&n=0,\\
-\sum_{j=n}^{-1}\tau_{\mathfrak M}^j(z),&n<0,
\end{cases}
\]
using the continuous $\widehat G$-action from Proposition~\ref{prop:automatic-admissibility} and Corollary~\ref{cor:restriction-tau-fh}. For $r\gg0$, pick $g_r\in H_K$ with $\chi(g_r)=1+p^r$, so that $g_r\to1$ in $H_K$ as $r\to\infty$. Since $z\in N_{\st}(\mathfrak M)$, we have
\[
S_{p^r}(z)=\tau_{\mathfrak M}^{-1}\bigl((g_r\otimes1)z-z\bigr)\longrightarrow0
\qquad(r\to\infty).
\]
Because the submodules $(p,E(u_0))^m\widehat{\mathfrak M}$ form a $\tau_{\mathfrak M}$-stable neighborhood basis of $0$ and
\[
S_{a+b}(z)=S_a(z)+\tau_{\mathfrak M}^{a}S_b(z)
\qquad(a,b\in\mathbf Z),
\]
the map $n\mapsto S_n(z)$ is uniformly continuous for the $p$-adic topology on $\mathbf Z$. By completeness of $\widehat{\mathfrak M}$, it extends uniquely to a continuous map $a\mapsto S_a(z)$ on $\mathbf Z_p$ satisfying $S_{a+b}(z)=S_a(z)+\tau_{\mathfrak M}^{a}S_b(z)$ for all $a,b\in\mathbf Z_p$.

For $h\in H_K$ with $q=\chi(h)\in\mathbf Z_{>0}$ and any $n>0$, conjugation and $z\in N_{\st}(\mathfrak M)$ give
\[
(h\otimes1)S_n(z)
=\sum_{j=0}^{n-1}\tau_{\mathfrak M}^{qj}\bigl((h\otimes1)z\bigr)
=\sum_{j=0}^{n-1}\tau_{\mathfrak M}^{qj}S_q(z)
=S_{qn}(z).
\]
Since $\{h\in H_K:\chi(h)\in\mathbf Z_{>0}\}$ contains the nonnegative powers of $\gamma$, it is dense in $H_K$. By continuity of the $\widehat G$-action, of $\chi$, and of $a\mapsto S_a(z)$, we obtain $(h\otimes1)S_a(z)=S_{\chi(h)a}(z)$ for all $(h,a)\in H_K\times\mathbf Z_p$. Setting $f(g):=S_{a(g)}(z)$ for $g=\tau^{a(g)}h$ defines a continuous map $f\colon\widehat G\to\widehat{\mathfrak M}$ vanishing on $H_K$ with $f(\tau)=z$. Using $(\tau^ah)(\tau^bk)=\tau^{a+\chi(h)b}hk$ and the two identities for $S_a(z)$, we find
\[
\begin{aligned}
f\bigl((\tau^a h)(\tau^b k)\bigr)
&=S_{a+\chi(h)b}(z)
=S_a(z)+\tau_{\mathfrak M}^{a}S_{\chi(h)b}(z)
=S_a(z)+\tau_{\mathfrak M}^{a}\bigl((h\otimes1)S_b(z)\bigr)\\
&=f(\tau^a h)+(\tau^a h)f(\tau^b k)
\end{aligned}
\]
for all $a,b\in\mathbf Z_p$ and $h,k\in H_K$, so $f\in\Ker d_{\mathrm{KL}}^1$.
\end{proof}

\subsection{From continuous cocycles to prismatic cocycles}

Propositions~\ref{prop:integral-cochain-comparison} and~\ref{prop:tau-evaluation-KL} immediately give $\Ker d^1_{\Prism,\log}\subset N_{\st}(\mathfrak M)$. For the reverse inclusion in Proposition~\ref{Prop ker d^1_prism}, we use Subsection~\ref{subsection Formal coordinates on the logarithmic diagonal} together with Lemma~\ref{lem:distribution-rigidity}.

\begin{lem}
\label{lem:distribution-rigidity}
Let $B$ be a $\mathbf Q_p$-algebra with a ring homomorphism $\rho_t\colon B\to B[[t]]$ satisfying $\rho_0=\operatorname{id}_B$ and $\rho_{s+t}=\rho_s\circ\rho_t$, and let $D$ be a finite $B$-module with a $\rho_t$-semilinear map $T_t\colon D\to D[[t]]$ satisfying $T_0=\operatorname{id}_D$ and $T_{s+t}=T_s\circ T_t$, where $s,t$ are independent variables.\footnote{Here, our convention is that when considering $\rho_s\circ\rho_t$ and $T_s\circ T_t$, we view $\rho_s$ and $T_s$ as their unique extensions $\rho_s\colon B[[t]]\to (B[[t]])[[s]]$ and $T_s\colon D[[t]]\to (D[[t]])[[s]]$ that keep the powers of $t$ unchanged.}

For an integer $q>1$ and $F(t)\in D[[t]]$, one has $F(qt)=\sum_{j=0}^{q-1}T_{jt}(F(t))$ if and only if $F(s+t)=F(s)+T_s(F(t))$, where $T_{jt}(F(t)):=\left.T_s(F(t))\right|_{s=jt}$.
\end{lem}

\begin{proof}
Setting $s=jt$ in $F(s+t)=F(s)+T_s(F(t))$ and summing over $0\leq j<q$ gives $F(qt)=\sum_{j=0}^{q-1}T_{jt}(F(t))$. Conversely, suppose $F(qt)=\sum_{j=0}^{q-1}T_{jt}(F(t))$. Evaluating at $t=0$ gives $(q-1)F(0)=0$, so $F(0)=0$. Put $v=F'(0)$, write $T_t(v)=\sum_{n\geq0}v_nt^n$, and set
\[
C(t)=\sum_{n\geq0}\frac{v_n}{n+1}t^{n+1}\in D[[t]].
\]
Since $C'(t)=T_t(v)$ and $T_s$ acts coefficientwise while fixing $t$, we have
\[
\partial_t\bigl(C(s)+T_s(C(t))\bigr)
=T_s(C'(t))=T_s(T_t(v))=T_{s+t}(v)=\partial_t C(s+t).
\]
Both sides equal $C(s)$ at $t=0$, so $C(s+t)=C(s)+T_s(C(t))$ and hence $C(qt)=\sum_{j=0}^{q-1}T_{jt}(C(t))$. By construction, $C(0)=0=F(0)$ and $C'(0)=T_0(v)=v=F'(0)$, so $F-C\in t^2D[[t]]$. If $F\neq C$, write $F-C=h_nt^n+O(t^{n+1})$ with $n\geq2$ minimal and $h_n\neq0$. Since $(F-C)(qt)=\sum_{j=0}^{q-1}T_{jt}((F-C)(t))$ and $T_{jt}(h_n)=h_n+O(t)$, comparing coefficients of $t^n$ gives $(q^n-q)h_n=0$. As $q>1$, $n\geq2$, and $B$ is a $\mathbf Q_p$-algebra, $q^n-q\in B^\times$, contradicting $h_n\neq0$. Thus $F=C$.
\end{proof}

\begin{prop}
\label{Prop ker d^1_prism}
We have
\[
\Ker d_{\Prism,\log}^1=N_{\st}(\mathfrak M),
\qquad
\Ker d_{\Prism}^1=N_{\cris}(\mathfrak M)
\]
in the logarithmic and crystalline cases, respectively, and in both cases $d^0(m)=\tau_{\mathfrak M}(1\otimes m)-1\otimes m$.
\end{prop}

\begin{proof}
When $\mathfrak M_{\Prism}$ is defined on $X_{\Prism}$, Proposition~\ref{prop:integral-cochain-comparison} gives the commutative diagram
\[
\begin{tikzcd}[column sep=large,row sep=large]
C_{\Prism}^1\arrow[r,hook,"j_1\otimes1"]
 \arrow[d,"d_{\Prism}^1"']
&C_{\Prism,\log}^1\arrow[r,hook]
 \arrow[d,"d_{\Prism,\log}^1"]
&K^1(\widehat{\mathfrak M})\arrow[d,"d_{\mathrm{KL}}^1"]\\
C_{\Prism}^2\arrow[r,hook,"j_2\otimes1"']
&C_{\Prism,\log}^2\arrow[r,hook]
&K^2(\widehat{\mathfrak M}).
\end{tikzcd}
\]
Since $\iota_{\mathbf a}^{\st}\circ j_m=\iota_{\mathbf a}$ for $\mathbf a\in\mathbf Z^m$ ($m=1,2$), Theorem~\ref{thm:multivariable-DL-injectivity} implies that $j_m$ (and hence $j_m\otimes1$, as $\mathfrak M$ is finite free over $\mathfrak S$) is injective. When $\mathfrak M_{\Prism}$ is defined only on $(X,M_X)_{\Prism}$, we use the right-hand square. Injectivity of the bottom horizontal arrows gives
\[
\Ker d_{\Prism,\log}^1=C_{\Prism,\log}^1\cap\Ker d_{\mathrm{KL}}^1,\qquad
\Ker d_{\Prism}^1=C_{\Prism}^1\cap\Ker d_{\mathrm{KL}}^1
\]
inside $K^1(\widehat{\mathfrak M})$. Together with Proposition~\ref{prop:tau-evaluation-KL} and the identities $\operatorname{ev}_{\tau}\circ\iota_1^{\st}=\operatorname{id}$ and $\operatorname{ev}_{\tau}\circ\iota_1^{\cris}=j_1$, this yields $\Ker d_{\Prism,\log}^1\subset N_{\st}(\mathfrak M)$ and $\Ker d_{\Prism}^1\subset N_{\cris}(\mathfrak M)$.

Conversely, let $z\in N_{\st}(\mathfrak M)$ and put
\[
B=\mathcal O_{\max}[1/p,1/E(u)],\qquad
D=B\otimes_{\mathfrak S}\mathfrak M,\qquad
F(u,t)=(\eta_1\otimes1)(z)\in D[[t]]
\]
as in Lemma~\ref{lem:faithful-log-Taylor}, where $\eta_1(u_0)=u$. Fix an $\mathfrak S$-basis $e_1,\ldots,e_d$ of $\mathfrak M$ to identify $D[[t]]\simeq B[[t]]^d$. By Corollary \ref{cor:crystal-tau-fh}, $\varepsilon_{\mathfrak M}^{\st}$ is determined by $\tau_{\mathfrak M}$. Let $A(u,t)\in\mathrm{GL}_d(B[[t]])$ be the matrix whose $j$-th column is $(\eta_1\otimes1)(\tau_{\mathfrak M}(1\otimes e_j))$. Since $\tau(u_0)=u_0(1+y)$ and $\tau(y)=y$, the action of $\tau$ induces under $\eta_1$ the map $\rho_t\colon B\to B[[t]]$ sending $u\mapsto e^tu$, and $\tau$-semilinearity of $\tau_{\mathfrak M}$ gives
\[
(\eta_1\otimes1)\bigl(\tau_{\mathfrak M}(v)\bigr)
=A(u,t)\rho_t\bigl((\eta_1\otimes1)(v)\bigr), \text{ for any }v\in\widehat{\mathfrak M}[1/p,1/E(u_0)].
\]
Put $T_s(w)=A(u,s)\rho_s(w)$ on $D$. We have $A(u,0)=I_d$, so $T_0=\operatorname{id}_D$, and the cocycle condition for $\varepsilon_{\mathfrak M}^{\st}$ gives $A(u,s+t)=A(u,s)\rho_s(A(u,t))$ and hence $T_{s+t}=T_s\circ T_t$. In particular, $T_t^j=T_{jt}$ for all $j\geq0$, so $(\eta_1\otimes1)(\tau_{\mathfrak M}^jz)=T_{jt}(F(u,t))$. For $q=\chi(\gamma)>1$, the action of $\gamma$ fixes $u$ and sends $t$ to $qt$, giving $(\eta_1\otimes1)((\gamma\otimes1)z)=F(u,qt)$. Thus $z\in N_{\st}(\mathfrak M)$ implies
\[
F(u,qt)=\sum_{j=0}^{q-1}T_{jt}\bigl(F(u,t)\bigr).
\]
By Lemmas~\ref{lem:faithful-log-Taylor} and~\ref{lem:distribution-rigidity},
\[
(\eta_2\otimes1)(d_{\Prism,\log}^1z)
=T_s\bigl(F(u,t)\bigr)-F(u,s+t)+F(u,s)=0.
\]
Since $\eta_2$ is injective after inverting $pE(u_0)$, and $pE(u_0)$ is a non-zero-divisor on $A_{\st}^{(3)}$ (which embeds into a product of copies of $A_{\inf}$ by Theorem~\ref{thm:multivariable-DL-injectivity}) and hence on $C_{\Prism,\log}^2$, we obtain $d_{\Prism,\log}^1z=0$.

If $z\in N_{\cris}(\mathfrak M)$, then $d_{\Prism,\log}^1(j_1\otimes1)(z)=0$, so injectivity of $j_2\otimes1$ in the left square gives $d_{\Prism}^1z=0$. Finally, evaluating $\varepsilon_{\mathfrak M}^{\st}$ at $\tau$ gives $\tau_{\mathfrak M}$, whence the formula for $d^0$.
\end{proof}

\subsection{The cohomological comparisons}

Passing to good truncations $\tau_{\leq1}$, we obtain a commutative diagram of two-term complexes:
\[
\begin{tikzcd}[column sep=large,row sep=large]
\tau_{\leq1}C_{\Prism}^{\bullet}\arrow[r,hook]\arrow[d,"\sim"']
&\tau_{\leq1}C_{\Prism,\log}^{\bullet}\arrow[r,"\sim"]\arrow[d,"\sim"']
&\tau_{\leq1}K^{\bullet}(\widehat{\mathfrak M})
 \arrow[d,"\operatorname{ev}_{\tau}","\sim"']\\
{[\mathfrak M\to N_{\cris}(\mathfrak M)]}\arrow[r,hook]
&{[\mathfrak M\to N_{\st}(\mathfrak M)]}\arrow[r,equal]
&{[\mathfrak M\to N_{\st}(\mathfrak M)]},
\end{tikzcd}
\]
where the bottom differentials are $\tau_{\mathfrak M}-1$, the vertical maps are the identity in degree $0$, and every arrow marked $\sim$ is an isomorphism of complexes (the top-right arrow by Propositions~\ref{prop:tau-evaluation-KL} and~\ref{Prop ker d^1_prism}). When $\mathfrak M_{\Prism}$ is defined only on $(X,M_X)_{\Prism}$, only the last two columns apply.

\begin{thm}
\label{prop compute prismatic coho}
Let $\mathfrak M_{\Prism}$ be a prismatic $F$-crystal on $(X,M_X)_{\Prism}$. There are natural quasi-isomorphisms in $D(\mathbf Z_p)$:
\[
R\Gamma((X,M_X)_{\Prism},\mathfrak M_{\Prism})
\simeq\bigl[\mathfrak M\xrightarrow{\tau_{\mathfrak M}-1}N_{\st}(\mathfrak M)\bigr]
\simeq\tau_{\leq1}K^{\bullet}(\widehat{\mathfrak M}),
\]
with $[\mathfrak M\xrightarrow{\tau_{\mathfrak M}-1}N_{\st}(\mathfrak M)]$ placed in degrees $0$ and $1$. If $\mathfrak M_{\Prism}$ is the pullback of a prismatic $F$-crystal on $X_{\Prism}$, then
\[
R\Gamma(X_{\Prism},\mathfrak M_{\Prism})
\simeq\bigl[\mathfrak M\xrightarrow{\tau_{\mathfrak M}-1}N_{\cris}(\mathfrak M)\bigr].
\]
\end{thm}

\begin{proof}
By \cite[Theorem~4.10]{gao2023hodge}, the cohomology of any finite locally free $\mathcal O_{\Prism}$-crystal on $X_{\Prism}$ (even without a Frobenius structure) is concentrated in degrees $0$ and $1$. The same argument applies on $(X,M_X)_{\Prism}$ by \cite[Lemma~2.2 and Notation~2.8]{gao2023hodge} (where $y_j$ differs by a sign from our convention). Thus $\tau_{\leq1}C^{\bullet}\hookrightarrow C^{\bullet}$ is a quasi-isomorphism for both $C_{\Prism}^{\bullet}$ and $C_{\Prism,\log}^{\bullet}$, and the theorem follows from Proposition~\ref{Prop ker d^1_prism} and the diagram above.
\end{proof}

\begin{cor}
\label{cor:relative-Galois-prismatic-H1}
There are natural isomorphisms
\[
\begin{aligned}
H^0((X,M_X)_{\Prism},\mathfrak M_{\Prism})
&\simeq\widehat{\mathfrak M}^{\widehat G},\\
H^1((X,M_X)_{\Prism},\mathfrak M_{\Prism})
&\simeq\Ker\bigl(H^1_{\mathrm{cont}}(\widehat G,\widehat{\mathfrak M})\longrightarrow H^1_{\mathrm{cont}}(H_K,\widehat{\mathfrak M})\bigr).
\end{aligned}
\]
When $\mathfrak M_{\Prism}$ comes from $X_{\Prism}$, the pullback map $H^i(X_{\Prism},\mathfrak M_{\Prism})\to H^i((X,M_X)_{\Prism},\mathfrak M_{\Prism})$ is an isomorphism for $i=0$ and an injection for $i=1$, whose image consists of classes represented by continuous $1$-cocycles $f\colon\widehat G\to\widehat{\mathfrak M}$ satisfying
\[
f|_{H_K}=0,\qquad
f(\tau)\in A^{(2)}\otimes_{p_0,\mathfrak S}\mathfrak M.
\]
\end{cor}

\begin{proof}
By Theorem~\ref{prop compute prismatic coho}, $H^i((X,M_X)_{\Prism},\mathfrak M_{\Prism})\simeq H^i(K^\bullet(\widehat{\mathfrak M}))$ for $i=0,1$. Since $K^0(\widehat{\mathfrak M})=\widehat{\mathfrak M}^{H_K}=1\otimes\mathfrak M$ by Lemma~\ref{lem:HK-invariants-prismatic}, its degree-$0$ cocycles are $\widehat{\mathfrak M}^{\widehat G}$. In degree $1$, $\Ker d_{\mathrm{KL}}^1$ consists of continuous $1$-cocycles $\widehat G\to\widehat{\mathfrak M}$ vanishing on $H_K$, so $K^\bullet(\widehat{\mathfrak M})\hookrightarrow C^\bullet_{\mathrm{cont}}(\widehat G,\widehat{\mathfrak M})$ induces a map $H^1(K^\bullet(\widehat{\mathfrak M}))\to\Ker\bigl(H^1_{\mathrm{cont}}(\widehat G,\widehat{\mathfrak M})\to H^1_{\mathrm{cont}}(H_K,\widehat{\mathfrak M})\bigr)$.

If $f$ is a continuous $1$-cocycle whose restriction to $H_K$ is $h\mapsto hv-v$ for some $v\in\widehat{\mathfrak M}$, then $f_0(g):=f(g)-(gv-v)$ is a cohomologous $1$-cocycle vanishing on $H_K$, hence in $\Ker d_{\mathrm{KL}}^1$, proving surjectivity. If two cocycles in $\Ker d_{\mathrm{KL}}^1$ differ by $g\mapsto gv-v$ ($v\in\widehat{\mathfrak M}$), restricting to $H_K$ gives $v\in\widehat{\mathfrak M}^{H_K}=1\otimes\mathfrak M$, so they differ by a coboundary in $K^\bullet(\widehat{\mathfrak M})$, proving injectivity.

When $\mathfrak M_{\Prism}$ comes from $X_{\Prism}$, the complexes $[\mathfrak M\to N_{\cris}(\mathfrak M)]$ and $[\mathfrak M\to N_{\st}(\mathfrak M)]$ share the same degree-$0$ term and differential while $N_{\cris}(\mathfrak M)\subset N_{\st}(\mathfrak M)$, giving the isomorphism on $H^0$ and injection on $H^1$. Proposition~\ref{prop:tau-evaluation-KL} describes the image on $H^1$, and this condition is independent of the choice of $f$ since adding $d_{\mathrm{KL}}^0(m)$ ($m\in\mathfrak M$) changes $f(\tau)$ by $(\tau_{\mathfrak M}-1)m\in A^{(2)}\otimes_{p_0,\mathfrak S}\mathfrak M$.
\end{proof}

\begin{rem}
\label{Cor total complex gives C_cris}
Suppose $\mathfrak M_{\Prism}$ is effective, and let $\star=\st$, or $\star=\cris$ when $\mathfrak M_{\Prism}$ comes from $X_{\Prism}$. Then Frobenius preserves $\mathfrak M$ and $N_\star(\mathfrak M)$ and commutes with $\tau_{\mathfrak M}-1$, so the fiber of $\varphi-1$ on $[\mathfrak M\xrightarrow{\tau_{\mathfrak M}-1}N_\star(\mathfrak M)]$ is the total complex of
\[
\begin{tikzcd}[column sep=large,row sep=large]
\mathfrak M\arrow[r,"\varphi_{\mathfrak M}-1"]
 \arrow[d,"\tau_{\mathfrak M}-1"']
&\mathfrak M\arrow[d,"\tau_{\mathfrak M}-1"]\\
N_\star(\mathfrak M)\arrow[r,"\varphi-1"']
&N_\star(\mathfrak M),
\end{tikzcd}
\]
namely the complex
\[
\big[\mathfrak M\xrightarrow{d^0}
\mathfrak M\oplus N_\star(\mathfrak M)\xrightarrow{d^1}
N_\star(\mathfrak M)\big]
\]
in degrees $0,1,2$, with $d^0(m)=\bigl((\varphi_{\mathfrak M}-1)m,(\tau_{\mathfrak M}-1)m\bigr)$ and $d^1(a,b)=(\tau_{\mathfrak M}-1)a-(\varphi-1)b$. This is the effective prismatic Herr complex studied in the next section. While the concentration of prismatic cohomology in degrees $0$ and $1$ implies that this complex has no cohomology above degree $2$, the vanishing of $H^2$ requires a separate argument. When $\mathfrak M_{\Prism}$ is not effective, Frobenius need not preserve $\mathfrak M$ or $N_\star(\mathfrak M)$, which leads to the $E$-localized complex of the next section.
\end{rem}

\section{Bloch--Kato local conditions and prismatic Herr complexes}
\label{sec:prismatic-Herr-complex}

In this section, we attach to crystalline or semistable prismatic $(\varphi,\tau)$-modules two explicit three-term complexes similar to the role of Herr complex for $(\varphi,\Gamma)$-modules. We will relate their first cohomology groups to the local Bloch--Kato groups $H_f^1$ and $H_g^1$. In subsections \ref{Section Crystalline Herr complex} and \ref{subsectino BK groups}, we begin with an $E$-localized complex for finite-height modules whose $H^1$ injects into the integral Bloch--Kato groups, and induces an isomorphism after inverting $p$. In Subsection \ref{s6:nekovar}, for crystalline representations attached to $p$-divisible groups over $\mathcal O_K$ we also show the the good truncation of the complex is canonically isomorphic to Nekov\'a\v r's complex in \cite{Nek18}. In Subsection \ref{s6:effective}, for the effective case, we study another three-term complex, which is a subcomplex of the one mentioned above, and establish its vanishing in degree $2$. Finally, in Subsection \ref{s6:classical-comparison} we compare both constructions with the \'etale $(\varphi,\tau)$-Herr complex.

\subsection{The $E$-localized prismatic Herr complex}
\label{Section Crystalline Herr complex}
\label{s6 finite-height}

Let $T$ be a lattice in a semistable $G_K$-representation $V=T[1/p]$, and let $\mathfrak M=\mathfrak M(T)$ be its finite-height prismatic $(\varphi,\tau)$-module as in Subsection~\ref{subsection finite height}. We say that $\mathfrak M$ is \emph{of type $\st$} when $T$ is semistable, and \emph{of type $\cris$} when $T$ is crystalline. In the latter case, $\mathfrak M$ satisfies the crystalline condition of Definition~\ref{def:prismatic-phi-tau}.

In what follows, fix $\star\in\{\cris,\st\}$, so that $T\in\operatorname{Rep}^{\star}_{\mathbf Z_p}(G_K)$ and $\mathfrak M\in\operatorname{Mod}^{\star,\mathrm{fh}}_{\varphi,\tau}(\mathfrak S)$ as in Subsection~\ref{subsection finite height}. Write $\mathfrak M_{\Prism}(T)$ for the associated prismatic $F$-crystal, which is defined on $X_{\Prism}$ when $\star=\cris$ and on $(X,M_X)_{\Prism}$ when $\star=\st$. We will use $\mathfrak S$ also to denote the prismatic $(\varphi,\tau)$-module corresponding to the trivial representation $\mathbf Z_p$.

\begin{defn}\label{defn: N_st^fh}
For $T$ and $\mathfrak M=\mathfrak M(T)$ as above, set
\[
N_{\st}^{\mathrm{fh}}(\mathfrak M)
=\big\{z\in A_{\st}^{(2)}[1/E(u_0)]\otimes_{i_0,\mathfrak S}\mathfrak M\big|
(g\otimes1)z=\sum_{j=0}^{\chi(g)-1}\tau_{\mathfrak M}^{j}(z)\text{ for all }g\in H_K\text{ with }\chi(g)\in\mathbf Z_{>0}\big\}.
\]
When $T$ is crystalline, we set
\[
N_{\cris}^{\mathrm{fh}}(\mathfrak M)
=N_{\st}^{\mathrm{fh}}(\mathfrak M)\cap
\left(A^{(2)}[1/E(u_0)]\otimes_{p_0,\mathfrak S}\mathfrak M\right),
\]
where the intersection is taken inside $A_{\st}^{(2)}[1/E(u_0)]\otimes_{i_0,\mathfrak S}\mathfrak M$. 

\smallskip

Clearly for the submodule $N_\star(\mathfrak M)$ of Definition~\ref{def of sub-zero}, we have $N_\star(\mathfrak M)\subset N_\star^{\mathrm{fh}}(\mathfrak M)$.
\end{defn}

\begin{defn}
\label{def ambient-finite-height-Herr}
The \emph{$E$-localized prismatic Herr complex} of $\mathfrak M\in\operatorname{Mod}^{\star,\mathrm{fh}}_{\varphi,\tau}(\mathfrak S)$ is the complex
\begin{equation}\label{s6:eq-fh}
\mathcal C_{\varphi,\tau}^{\star,\mathrm{fh}}(\mathfrak M)
=\big[\mathfrak M\xrightarrow{d^0_{\mathrm{fh}}}
\mathfrak M_E\oplus N_\star(\mathfrak M)
\xrightarrow{d^1_{\mathrm{fh}}}N_\star^{\mathrm{fh}}(\mathfrak M)\big]
\end{equation}
in degrees $0,1,2$, with differentials
\[
d^0_{\mathrm{fh}}(m)
=\bigl((\varphi_{\mathfrak M}-1)m,(\tau_{\mathfrak M}-1)m\bigr),\qquad
d^1_{\mathrm{fh}}(a,b)=(\tau_{\mathfrak M}-1)a-(\varphi-1)b.
\]
For $T\in\operatorname{Rep}^{\star}_{\mathbf Z_p}(G_K)$ and $V=T[1/p]$, we write $\mathcal C_{\varphi,\tau}^{\star,\mathrm{fh}}(T)=\mathcal C_{\varphi,\tau}^{\star,\mathrm{fh}}(\mathfrak M(T))$ and put
\[
\mathcal C_{\varphi,\tau}^{\star,\mathrm{fh}}(V)
:=\mathcal C_{\varphi,\tau}^{\star,\mathrm{fh}}(T)[1/p],
\]
which depends only on $V$ since any two $G_K$-stable lattices in $V$ are commensurable.
\end{defn}

Let us check that \eqref{s6:eq-fh} is well-defined. Given $a\in\mathfrak M_E$ and $g\in H_K$ with $\chi(g)\in\mathbf Z_{>0}$, the conjugation relation between $g$ and $\tau_{\mathfrak M}$ imply that $(\tau_{\mathfrak M}-1)a\in N_{\st}^{\mathrm{fh}}(\mathfrak M)$ by direct computations:
since $H_K$ fixes $1\otimes\mathfrak M_E$, the conjugation relation gives
\[
(g\otimes1)\bigl((\tau_{\mathfrak M}-1)a\bigr)
=\tau_{\mathfrak M}^{\chi(g)}(g\otimes1)a-(g\otimes1)a\
=(\tau_{\mathfrak M}^{\chi(g)}-1)a\
=\sum_{j=0}^{\chi(g)-1}\tau_{\mathfrak M}^{j}
\bigl((\tau_{\mathfrak M}-1)a\bigr).
\] If $\star=\cris$, the crystalline descent isomorphism after inverting $E(u_0)$ further places $(\tau_{\mathfrak M}-1)a$ in $A^{(2)}[1/E(u_0)]\otimes_{p_0,\mathfrak S}\mathfrak M$, and hence in $N_{\cris}^{\mathrm{fh}}(\mathfrak M)$. Similarly, $\varphi$ maps each coefficient module into its $E(u_0)$-localization and commutes with both $H_K$ and $\tau_{\mathfrak M}$, so it preserves the defining equations of $N_\star$ and carries $N_\star(\mathfrak M)$ into $N_\star^{\mathrm{fh}}(\mathfrak M)$. Because $\varphi$ and $\tau_{\mathfrak M}$ commute, we have $d^1_{\mathrm{fh}}\circ d^0_{\mathrm{fh}}=0$.

One can equivalently view $\mathcal C_{\varphi,\tau}^{\star,\mathrm{fh}}(\mathfrak M)$ as the mapping fiber
\[
\operatorname{fib}\left(
\varphi-\mathsf{can}:
\bigl[\mathfrak M\xrightarrow{\tau_{\mathfrak M}-1}N_\star(\mathfrak M)\bigr]
\longrightarrow
\bigl[\mathfrak M_E\xrightarrow{\tau_{\mathfrak M}-1}N_\star^{\mathrm{fh}}(\mathfrak M)\bigr]
\right),
\]
where $\mathsf{can}$ are the natural inclusions $\mathfrak M\hookrightarrow\mathfrak M_E$ and $N_\star(\mathfrak M)\hookrightarrow N_\star^{\mathrm{fh}}(\mathfrak M)$ in degrees $0$ and $1$. Its good truncation is
\begin{equation}\label{s6 eq-truncation}
\tau_{\leq1}\mathcal C_{\varphi,\tau}^{\star,\mathrm{fh}}(T)
=\bigl[\mathfrak M(T)\xrightarrow{d^0_{\mathrm{fh}}}\Ker d^1_{\mathrm{fh}}\bigr].
\end{equation}

\begin{rem}
The condition defining $N_{\st}^{\mathrm{fh}}(\mathfrak M)$ in Definition~\ref{defn: N_st^fh} (and hence also $N_{\cris}^{\mathrm{fh}}(\mathfrak M)$) can be checked on a single topological generator $\gamma\in H_K$: in fact, it suffices to verify it at any single element $g_0\in H_K$ with $q_0=\chi(g_0)\in\mathbf Z_{>1}$.

Indeed, retain the notation $D$ and $T_s$ from the proof of Proposition~\ref{Prop ker d^1_prism}, and put $F(u,t)=(\eta_1\otimes1)(z)\in D[[t]]$ for $z\in\widehat{\mathfrak M}[1/E(u_0)]$. As in the proof of Proposition~\ref{Prop ker d^1_prism}, the defining relation at $g_0$ translates into
\[
F(u,q_0t)=\sum_{j=0}^{q_0-1}T_{jt}\bigl(F(u,t)\bigr).
\]
By Lemma~\ref{lem:distribution-rigidity}, this implies $F(u,s+t)=F(u,s)+T_s(F(u,t))$, which iterates to $F(u,dt)=\sum_{j=0}^{d-1}T_{jt}(F(u,t))$ for every integer $d\geq1$. Use the above equivalence again shows that the required relation holds for every $h\in H_K$ with $\chi(h)=d\in\N$.
\end{rem}

\subsection{Bloch--Kato local conditions}\label{subsectino BK groups}

Recall that for $V\in\operatorname{Rep}^{\star}_{\mathbf Q_p}(G_K)$ with $\star\in\{\cris,\st\}$,
\[
H_\star^1(K,V)\coloneqq \Ker\!\left(H^1(K,V)\longrightarrow
H^1(K,\mathbf B_\star\otimes_{\mathbf Q_p}V)\right),
\]
where $\mathbf B_{\cris}$ and $\mathbf B_{\st}$ are Fontaine's period rings, note that $H_{\cris}^1=H_f^1$ is the crystalline Bloch--Kato subgroup of \cite{Bloch-Kato90L-functions-Tamagawa-number}. For a $G_K$-stable $\mathbf Z_p$-lattice $T\subset V$, set
\begin{equation}\label{s6:eq-integral-local-condition}
H_\star^1(K,T)=H^1(K,T)\times_{H^1(K,V)}H_\star^1(K,V).
\end{equation}
When $V$ is semistable, Hyodo's theorem, extended to arbitrary perfect residue fields by Berger \cite[Th\'eor\`eme~6.2]{Ber02}, identifies $H_{\st}^1(K,V)$ with $H_g^1(K,V)$ (defined using $\mathbf B_{\dR}$ in place of $\mathbf B_{\st}$), and hence $H_{\st}^1(K,T)=H_g^1(K,T)$ under the above conventions. Concretely, $H_\star^1(K,T)$ classifies extensions $0\to T\to U\to\mathbf Z_p\to0$ of $G_K$-lattices with $U\in\operatorname{Rep}^{\star}_{\mathbf Z_p}(G_K)$:
\[
H_\star^1(K,T)
=\operatorname{Ext}^1_{\operatorname{Rep}^{\star}_{\mathbf Z_p}(G_K)}(\mathbf Z_p,T).
\]

\begin{thm}
\label{s6 thm-fh-main}
Let $T\in\operatorname{Rep}^{\star}_{\mathbf Z_p}(G_K)$. There are natural isomorphisms
\[
H^0\!\left(\mathcal C_{\varphi,\tau}^{\star,\mathrm{fh}}(T)\right)\simeq H^0(K,T), \qquad
H^1\!\left(\mathcal C_{\varphi,\tau}^{\star,\mathrm{fh}}(T)\right)
\simeq\operatorname{Ext}^1_{\operatorname{Mod}^{\star,\mathrm{fh}}_{\varphi,\tau}(\mathfrak S)}
(\mathfrak S,\mathfrak M(T)),
\]
and the latter extension group is naturally isomorphic to $\operatorname{Ext}^1_{\operatorname{Vect}^{\varphi}(X_{\Prism})}(\mathcal O_{\Prism},\mathfrak M_{\Prism}(T))$ when $\star=\cris$, and to $\operatorname{Ext}^1_{\operatorname{Vect}^{\varphi}((X,M_X)_{\Prism})}(\mathcal O_{\Prism},\mathfrak M_{\Prism}(T))$ when $\star=\st$. Furthermore, the functor $T_{\mathrm{cov}}$ induces an injection
\[
H^1\!\left(\mathcal C_{\varphi,\tau}^{\star,\mathrm{fh}}(T)\right)
\lhook\joinrel\longrightarrow H_\star^1(K,T),
\]
which becomes an isomorphism after inverting $p$. Moreover, the above inclusion is an isomorphism when $\operatorname{HT}(T)\subset[-r,0]$ for some $r\geq0$.
\end{thm}

We will prepare for the proof of  Theorem~\ref{s6 thm-fh-main} with the next three propositions.

\begin{prop}
\label{lem remove intersection with W}
Let $T\in\operatorname{Rep}^{\st}_{\mathbf Z_p}(G_K)$ and put $\mathfrak M=\mathfrak M(T)$. Then
\[
N_{\st}(\mathfrak M)
=\bigl((A_{\st}^{(2)}\cap W(\mathfrak m_{C_p}^{\flat}))
\otimes_{i_0,\mathfrak S}\mathfrak M\bigr)_0,
\]
where the subscript $0$ denotes the submodule cut out by the equations of Definition~\ref{def of sub-zero}. If $T$ is crystalline, the same holds with $N_{\cris}$, $A^{(2)}$ in place of $N_{\st}$, $A_{\st}^{(2)}$.
\end{prop}

\begin{proof}
For $z\in N_{\st}(\mathfrak M)$, let $\bar z$ be its image modulo $I_+$ in $\widehat{\mathfrak M}/I_+\widehat{\mathfrak M}\simeq W(k)\otimes_{\mathfrak S}\mathfrak M$. Because both $H_K$ and $\tau_{\mathfrak M}$ act trivially on this quotient, reducing the defining equation of $N_{\st}(\mathfrak M)$ modulo $I_+$ gives $\bar z=\chi(g)\bar z$ whenever $g\in H_K$ has $\chi(g)\in\mathbf Z_{>0}$. Picking $g\in H_K$ with $\chi(g)=1+p^n$ for $n\gg0$ gives $p^n\bar z=0$, and hence $\bar z=0$ as $W(k)\otimes_{\mathfrak S}\mathfrak M$ is $p$-torsion-free. Thus $z\in I_+\widehat{\mathfrak M}$, which proves one inclusion. The reverse inclusion is immediate, and intersecting with $A^{(2)}\otimes_{p_0,\mathfrak S}\mathfrak M$ gives the crystalline case.
\end{proof}

\begin{prop}
\label{prop ambient-Ext-Herr}
For $\mathfrak M\in\operatorname{Mod}^{\star,\mathrm{fh}}_{\varphi,\tau}(\mathfrak S)$, there is a natural isomorphism
\[
H^1\!\left(\mathcal C^{\star,\mathrm{fh}}_{\varphi,\tau}(\mathfrak M)\right)
\simeq\operatorname{Ext}^1_{\operatorname{Mod}^{\star,\mathrm{fh}}_{\varphi,\tau}(\mathfrak S)}
(\mathfrak S,\mathfrak M),
\]
and the same holds in the isogeny category \textup{(with exactness taken on the underlying $\mathfrak S[1/p]$-modules)}.
\end{prop}

\begin{proof}
Given an admissible extension $0\to\mathfrak M\to\mathfrak E\xrightarrow{\textsf{pr}}\mathfrak S\to0$, choose an $\mathfrak S$-linear splitting $s\mapsto se$ with $e\in\mathfrak E$ lifting $1\in\mathfrak S$, so that $\mathfrak E=\mathfrak M\oplus\mathfrak S e$. Write
\begin{equation}\label{eq:extension-lift}
\varphi_{\mathfrak E}(e)=e+\lambda,\qquad
\tau_{\mathfrak E}(1\otimes e)=1\otimes e+\mu.
\end{equation}
Clearly $\lambda\in\Ker(\textsf{pr}[1/E])=\mathfrak M_E$. For any $g\in H_K$ with $\chi(g)\in\mathbf Z_{>0}$, evaluating $G_g\tau_{\mathfrak E}=\tau_{\mathfrak E}^{\chi(g)}$ on the $G_g$-invariant element $1\otimes e$ gives
\[
1\otimes e+(g\otimes1)\mu
=\tau_{\mathfrak E}^{\chi(g)}(1\otimes e)
=1\otimes e+\sum_{j=0}^{\chi(g)-1}\tau_{\mathfrak M}^{j}(\mu),
\]
hence $\mu\in N_{\st}(\mathfrak M)$. When $\star=\cris$, we have $\tau_{\mathfrak E}(1\otimes e)\in A^{(2)}\otimes_{p_0,\mathfrak S}\mathfrak E$. Since $\mu\in\Ker(1\otimes\textsf{pr})=\widehat{\mathfrak M}$, it follows that $\mu\in A^{(2)}\otimes_{p_0,\mathfrak S}\mathfrak M$ and hence $\mu\in N_{\cris}(\mathfrak M)$. Moreover, commutativity of $\varphi$ and $\tau_{\mathfrak E}$,
applied to $1\otimes e$, gives $1\otimes\lambda+\varphi(\mu)
=
\mu+\tau_{\mathfrak M}(1\otimes\lambda)$
in $\widehat{\mathfrak M}[1/E(u_0)]$. Hence $(\tau_{\mathfrak M}-1)\lambda-(\varphi-1)\mu=0$, so $(\lambda,\mu)\in\Ker(d^1_{\mathrm{fh}})$.

Conversely, given $(\lambda,\mu)\in\Ker(d^1_{\mathrm{fh}})\subset\mathfrak M_E\oplus N_\star(\mathfrak M)$, define $\varphi_{\mathfrak E}$ and $\tau_{\mathfrak E}$ on $\mathfrak E=\mathfrak M\oplus\mathfrak S e$ by extending the operators on $\mathfrak M$ via \eqref{eq:extension-lift}. Then $d^1_{\mathrm{fh}}(\lambda,\mu)=0$ ensures that $\varphi_{\mathfrak E}$ and $\tau_{\mathfrak E}$ commute on $1\otimes e$, while $\mu\in N_\star(\mathfrak M)$ gives Condition~(3) of Definition~\ref{def:prismatic-phi-tau} on $1\otimes e$. Since $\lambda\in\mathfrak M_E$, both the linearized Frobenius of $\mathfrak E$ and its inverse have bounded $E(u)$-denominators, so $\mathfrak E$ has finite height. Relative to $\widehat{\mathfrak E}=\widehat{\mathfrak M}\oplus A_{\st}^{(2)}(1\otimes e)$, the matrix of $\tau_{\mathfrak E}$ is upper triangular with invertible diagonal blocks and off-diagonal entry $\mu\in I_+\widehat{\mathfrak M}$ by Proposition~\ref{lem remove intersection with W}. Thus $\tau_{\mathfrak E}$ is a $\tau$-semilinear automorphism satisfying the finite-height conditions of Definition~\ref{def:prismatic-phi-tau}, as well as the crystalline condition when $\mu\in N_{\cris}(\mathfrak M)$.

Changing the lift $e$ to $e+m$ ($m\in\mathfrak M$) shifts $(\lambda,\mu)$ by $d^0_{\mathrm{fh}}(m)$, and every isomorphism of extensions carries $e$ to a lift of this form. 

Over $\mathfrak S[1/p]$, if $(\lambda,\mu)$ is a $1$-cocycle of $\mathcal C_{\varphi,\tau}^{\star,\mathrm{fh}}(\mathfrak M)[1/p]$, choose $n\geq0$ with $p^n\lambda\in\mathfrak M_E$ and $p^n\mu\in N_\star(\mathfrak M)$. Then
\[
\varphi(p^ne)=p^ne+p^n\lambda,\qquad
\tau(1\otimes p^ne)=1\otimes p^ne+p^n\mu
\]
define an admissible integral extension on $\mathfrak M\oplus\mathfrak S(p^ne)$ with projection $p^ne\mapsto1$, whose pullback along $p^{-n}:\mathfrak S[1/p]\to\mathfrak S[1/p]$ after inverting $p$ realizes the cocycle $(\lambda,\mu)$.
\end{proof}

\begin{rem}
The isomorphism of Proposition~\ref{prop ambient-Ext-Herr} respects the group structures, matching $[(\lambda_1,\mu_1)]+[(\lambda_2,\mu_2)]=[(\lambda_1+\lambda_2,\mu_1+\mu_2)]$ on $H^1$ with the Baer sum on $\operatorname{Ext}^1$. Indeed, for extensions $\mathfrak E_1$ and $\mathfrak E_2$ of $\mathfrak S$ by $\mathfrak M$ with lifts $e_i\in\mathfrak E_i$ of $1$ and associated cocycles $(\lambda_i,\mu_i)$ as in \eqref{eq:extension-lift}, the Baer sum is constructed from $\mathfrak E_1\oplus\mathfrak E_2$ by pullback along the diagonal $\mathfrak S\to\mathfrak S\oplus\mathfrak S$ and pushout along addition $\mathfrak M\oplus\mathfrak M\to\mathfrak M$, and the image of $(e_1,e_2)$ satisfies \eqref{eq:extension-lift} with $(\lambda_1+\lambda_2,\mu_1+\mu_2)$.
\end{rem}

\begin{prop}
\label{prop exactness-trivial-quotient}
The covariant Breuil--Kisin functor $T\mapsto\mathfrak M(T)$ on crystalline and semistable representations is exact after inverting $p$. If $T\in\operatorname{Rep}^{\star}_{\mathbf Z_p}(G_K)$ has Hodge--Tate weights in $[-r,0]$ for some integer $r\geq0$ and we set $\mathfrak M^{(r)}(T)=\mathfrak M(T(-r))$, there are natural integral isomorphisms
\[
H_\star^1(K,T)
\simeq\operatorname{Ext}^1_{\operatorname{Mod}^{\star,[0,r]}_{\varphi,\tau}(\mathfrak S)}
\bigl(\mathfrak S\{-r\},\mathfrak M^{(r)}(T)\bigr)
\simeq\operatorname{Ext}^1_{\operatorname{Mod}^{\star,\mathrm{fh}}_{\varphi,\tau}(\mathfrak S)}
\bigl(\mathfrak S,\mathfrak M(T)\bigr).
\]
\end{prop}

\begin{proof}
For rational exactness, we adapt \cite[Proposition~3.4]{CGM26}. Given an exact sequence $0\to V_1\to V_2\to V_3\to0$ in $\operatorname{Rep}^{\star}_{\mathbf Q_p}(G_K)$, choose $G_K$-stable lattices $T_i\subset V_i$ such that $0\to T_1\to T_2\to T_3\to0$ is exact, and view $0\to\mathfrak M(T_1)\to\mathfrak M(T_2)\to\mathfrak M(T_3)\to0$ as a complex in the abelian category $\operatorname{Mod}^{\varphi}_{\mathfrak S}$ recalled in Remark~\ref{rmk def of fh case}. Consider the functor of \cite[Definition~2.6]{CGM26},
\[
T_{\mathfrak S}(\mathfrak N)
=\bigl(W(C_p^{\flat})\otimes_{\mathfrak S}\mathfrak N\bigr)^{\varphi=1}
\qquad(\mathfrak N\in\operatorname{Mod}^{\varphi}_{\mathfrak S}),
\]
which takes values in finitely generated $\mathbf Z_p$-modules with continuous $G_{K_\infty}$-action via $W(C_p^\flat)$ and is exact by \cite[Proposition~2.8]{CGM26}. Since forgetting the $\tau$-structure corresponds to restricting the associated Galois representation to $G_{K_\infty}$, we have functorial identifications $T_{\mathfrak S}(\mathfrak M(T_i))\simeq T_{\mathrm{cov}}(\mathfrak M(T_i))|_{G_{K_\infty}}\simeq T_i|_{G_{K_\infty}}$. Thus $T_{\mathfrak S}$ turns $0\to\mathfrak M(T_1)\to\mathfrak M(T_2)\to\mathfrak M(T_3)\to0$ into the exact sequence $0\to T_1|_{G_{K_\infty}}\to T_2|_{G_{K_\infty}}\to T_3|_{G_{K_\infty}}\to0$, so $T_{\mathfrak S}(H)=0$ for each homology module $H$. By \cite[Proposition~2.9 and Lemma~2.10]{CGM26}, each $H$ is killed by a power of $(p,u)$ and hence by a power of $p$, proving that $0\to\mathfrak M(T_1)[1/p]\to\mathfrak M(T_2)[1/p]\to\mathfrak M(T_3)[1/p]\to0$ is exact.

The first isomorphism follows from \cite[Theorem~3.12]{CGM26}, see also Remark~3.13(3) of \emph{loc.\ cit.}, notice their interval $[0,r]$ becomes $[-r,0]$
in our convention. We have
\[
T_{\mathrm{cov}}(\mathfrak M^{(r)}(T))(r)
\simeq T(-r)(r)\simeq T,\qquad
T_{\mathrm{cov}}(\mathfrak S\{-r\})(r)
\simeq\mathbf Z_p(-r)(r)\simeq\mathbf Z_p.
\] 

Because $-\otimes_{\mathfrak S}\mathfrak S\{r\}$ is exact and sends $\mathfrak S\{-r\}$ and $\mathfrak M^{(r)}(T)$ to $\mathfrak S$ and $\mathfrak M(T)$, compatibility of $T_{\mathrm{cov}}$ with Tate twists gives the commutative diagram
\[
\begin{tikzcd}[row sep=large]
\operatorname{Ext}^1_{\operatorname{Mod}^{\star,[0,r]}_{\varphi,\tau}(\mathfrak S)}
(\mathfrak S\{-r\},\mathfrak M^{(r)}(T))
\arrow[d,"-\otimes_{\mathfrak S}\mathfrak S\{r\}"']
\arrow[dd,bend left=75,"T_{\mathrm{cov}}(-)(r)","\simeq"']
\\
\operatorname{Ext}^1_{\operatorname{Mod}^{\star,\mathrm{fh}}_{\varphi,\tau}(\mathfrak S)}
(\mathfrak S,\mathfrak M(T))
\arrow[d,"T_{\mathrm{cov}}"']
\\
H_\star^1(K,T).
\end{tikzcd}
\]
Here the lower vertical arrow is well-defined and injective, since $T_{\mathrm{cov}}$ is fully faithful and preserves admissible short exact sequences (by the exactness of $T_{\mathfrak S}$ upon restriction to $G_{K_\infty}$). The above discussions give the second isomorphism.
\end{proof}

\begin{proof}[Proof of Theorem~\ref{s6 thm-fh-main}]
Elements of $H^0(\mathcal C_{\varphi,\tau}^{\star,\mathrm{fh}}(T))$ are elements of $\mathfrak M(T)$ fixed by $\varphi_{\mathfrak M}$ and $\tau_{\mathfrak M}$. These correspond to morphisms $\mathfrak S\to\mathfrak M(T)$ of prismatic $(\varphi,\tau)$-modules, and hence via $T_{\mathrm{cov}}$ to $G_K$-equivariant maps $\mathbf Z_p\to T$, or $T^{G_K}$.

In degree $1$, Proposition~\ref{prop ambient-Ext-Herr} identifies $H^1(\mathcal C_{\varphi,\tau}^{\star,\mathrm{fh}}(T))$ with $\operatorname{Ext}^1_{\operatorname{Mod}^{\star,\mathrm{fh}}_{\varphi,\tau}(\mathfrak S)}(\mathfrak S,\mathfrak M(T))$. Evaluation on the Breuil--Kisin prism together with restriction of the descent datum to $\tau$ yields exact equivalences between prismatic $F$-crystals and prismatic $(\varphi,\tau)$-modules in both the crystalline and semistable cases (see Corollary~\ref{cor:crystal-tau-fh} and the discussion of exact structures following it), inducing the two identifications with crystal extension groups. Further, $T_{\mathrm{cov}}$ takes any admissible extension $0\to\mathfrak M(T)\to\mathfrak E\to\mathfrak S\to0$ to an exact sequence $0\to T\to T_{\mathrm{cov}}(\mathfrak E)\to\mathbf Z_p\to0$ in $\operatorname{Rep}^{\star}_{\mathbf Z_p}(G_K)$. If the latter splits, full faithfulness of $T_{\mathrm{cov}}$ lifts a $G_K$-equivariant section $\mathbf Z_p\to T_{\mathrm{cov}}(\mathfrak E)$ to a splitting $\mathfrak S\to\mathfrak E$ of $\mathfrak E\to\mathfrak S$, so the induced map into $H_\star^1(K,T)$ is injective. Both the isomorphism after inverting $p$ and the integral isomorphism for weights in $[-r,0]$ follow from Proposition~\ref{prop exactness-trivial-quotient}.
\end{proof}

\begin{rem}
\label{s6 cor-saturation}
Let $L_\star(T)$ denote the image of $H^1(\mathcal C_{\varphi,\tau}^{\star,\mathrm{fh}}(T))\hookrightarrow H_\star^1(K,T)$ in Theorem~\ref{s6 thm-fh-main}, and put $Q_\star(T)=H_\star^1(K,T)/L_\star(T)$. Then $H_\star^1(K,T)$ is the $p$-saturation of $L_\star(T)$ in $H^1(K,T)$, \textit{i.e.}
$H_\star^1(K,T)=\{x\in H^1(K,T):p^nx\in L_\star(T)\text{ for some }n\geq0\}$.

Indeed, $Q_\star(T)[1/p]=0$ by Theorem~\ref{s6 thm-fh-main}, so every $x\in H_\star^1(K,T)$ has $p^nx\in L_\star(T)$ for some $n\geq0$. Conversely, if $p^nx\in L_\star(T)$, then $x$ maps into the $\mathbf Q_p$-subspace $H_\star^1(K,V)\subset H^1(K,V)$, so $x\in H_\star^1(K,T)$ by \eqref{s6:eq-integral-local-condition}. When $K/\mathbf Q_p$ is finite, $H_\star^1(K,T)$ is finitely generated over $\mathbf Z_p$ and hence $Q_\star(T)$ is finite. When $\operatorname{HT}(T)\subset[-r,0]$ ($r\geq0$), $Q_\star(T)=0$ by Theorem~\ref{s6 thm-fh-main}.
\end{rem}

\subsection{Comparison with Nekov\'a\v r's complex}
\label{s6:nekovar}

Assume throughout this subsection that $K/\mathbf Q_p$ is finite. Let $G$ be a $p$-divisible group over $\mathcal O_K$ with $T=T_pG$, and write $\mathfrak M_{\mathrm{Kis}}$ for Kisin's contravariant functor from $p$-divisible groups over $\mathcal O_K$ to Breuil--Kisin modules of height at most $1$ \cite{kisin-crystalline}. Following \cite[\S2.4]{Nek18}, set $N=\mathfrak M_{\mathrm{Nek}}(G)=\mathfrak M_{\mathrm{Kis}}(G^\vee)$ and $N^0=\operatorname{im}(\Phi_N:\varphi^*N\to N)$, so that $N\simeq\mathfrak M(T(-1))$ as Breuil--Kisin modules. Put
\[
\mathcal L=\mathfrak S\{-1\}=\mathfrak S e_{-1},\qquad
c(u)=\frac{E(u)}{c_0}=\frac{pE(u)}{E(0)},
\]
so $\varphi_{\mathcal L}(e_{-1})=c(u)e_{-1}$. Nekov\'a\v r's complex
\[
C_{\mathrm{Nek}}^\bullet(G)
=\bigl[N\xrightarrow{\varphi_N-c(u)}N^0\bigr]
\]
(in degrees $0$ and $1$) satisfies $H^1(C_{\mathrm{Nek}}^\bullet(G))\simeq H_f^1(K,T)$ by \cite[Proposition~2.5 and Corollary~2.6]{Nek18}. We shall identify $C_{\mathrm{Nek}}^\bullet(G)$ canonically with $\tau_{\leq1}\mathcal C_{\varphi,\tau}^{\cris,\mathrm{fh}}(T)$.

Let $\operatorname{Mod}^{\mathrm{Nek},[0,1]}_\varphi$ be the category of finite free Breuil--Kisin modules of height at most $1$, with exact structure inherited from $\mathfrak S$-modules. Since $E(u)N\subset N^0$, the differential $\varphi_N-c(u)$ lands in $N^0$. More generally, an $\mathfrak S$-split extension $N\oplus\mathfrak S\widetilde e$ of $\mathcal L$ by $N$ with $\varphi(\widetilde e)=c(u)\widetilde e+x$ ($x\in N$) has height at most $1$ if and only if $x\in N^0$: as $N$ has height at most $1$ and $c_0\in\mathfrak S^\times$, the vector $E(u)\widetilde e=c_0\varphi(\widetilde e)-c_0x$ lies in the image of the linearized Frobenius precisely when $x\in N^0$.

\begin{prop}
\label{prop nek-normalization}
Follow the above notations. Endow each Breuil--Kisin module of height at most $1$ with its unique crystalline prismatic $(\varphi,\tau)$-structure  compatible with its given Frobenius. Then
\[
D:\operatorname{Mod}^{\mathrm{Nek},[0,1]}_\varphi
\xrightarrow{\sim}\operatorname{Mod}^{\cris,[-1,0]}_{\varphi,\tau}(\mathfrak S),
\qquad M\longmapsto M\otimes_{\mathfrak S}\mathcal L^\vee,
\]
is an exact equivalence sending $\mathcal L$ to $\mathfrak S$ and $N$ to $\mathfrak M(T)$. Moreover, every admissible extension of $\mathfrak S$ by $\mathfrak M(T)$ in $\operatorname{Mod}^{\cris,\mathrm{fh}}_{\varphi,\tau}(\mathfrak S)$ has height in $[-1,0]$.
\end{prop}

\begin{proof}
By Kisin's classification of $p$-divisible groups and its
compatibility with Tate modules \cite[\S2.2]{kisin-crystalline},
the Breuil--Kisin functor gives an equivalence between
crystalline $G_K$-lattices with Hodge--Tate weights in $[0,1]$
and finite free Breuil--Kisin modules of height at most one. Combining this with
\cite[Theorem~4.1.10 and Remark~4.3.5]{du-liu-prismaticphiGhatmodule}
shows that each such module admits a unique crystalline
prismatic $(\varphi,\tau)$-structure.
Moreover, every morphism of the underlying Breuil--Kisin
modules commutes with the resulting $\tau$-operators after
base change to $A_{\st}^{(2)}$ along $i_0$. Tensoring with $\mathcal L^\vee=\mathfrak S\{1\}$ shifts the height interval from $[0,1]$ to $[-1,0]$, with inverse given by tensoring with $\mathcal L$ and forgetting the $\tau$-action, and both functors preserve $\mathfrak S$-exactness. With this crystalline structure on $N$, we have $T_{\mathrm{cov}}(D(N))\simeq T_{\mathrm{cov}}(N)(1)\simeq T(-1)(1)\simeq T$, and $D(N)\simeq\mathfrak M(T)$.

If $0\to\mathfrak M(T)\to\mathfrak E\to\mathfrak S\to0$ is an admissible extension in $\operatorname{Mod}^{\cris,\mathrm{fh}}_{\varphi,\tau}(\mathfrak S)$, then $0\to T\to T_{\mathrm{cov}}(\mathfrak E)\to\mathbf Z_p\to0$ is exact in $\operatorname{Rep}^{\cris}_{\mathbf Z_p}(G_K)$, so the Hodge--Tate weights of $T_{\mathrm{cov}}(\mathfrak E)$ lie in $[-1,0]$. Since $\mathfrak E\simeq\mathfrak M(T_{\mathrm{cov}}(\mathfrak E))$, the height bound of Subsection~\ref{subsection finite height} shows that $\mathfrak E$ has height in $[-1,0]$.
\end{proof}

Keep these crystalline prismatic $(\varphi,\tau)$-structures on $N$, and recall $\mathfrak M=D(N)$ and $\widehat{\mathfrak M}=A_{\st}^{(2)}\otimes_{i_0,\mathfrak S}\mathfrak M$. Let $e_1\in\mathcal L^\vee$ be dual to $e_{-1}$, and define $\alpha\in(A^{(2)})^\times$ by
\[
\tau_{\mathcal L}(1\otimes e_{-1})=\alpha(1\otimes e_{-1}).
\]
For $x\in N^0$, let $E_x=N\oplus\mathfrak S\widetilde e\in\operatorname{Mod}^{\mathrm{Nek},[0,1]}_\varphi$ have Frobenius and crystalline $\tau$-action given by
\begin{equation}\label{eq:nek-extension-lift}
\varphi_{E_x}(\widetilde e)=c(u)\widetilde e+x,\qquad
\tau_{E_x}(1\otimes\widetilde e)=\alpha(1\otimes\widetilde e)+b_x,
\end{equation}
with $b_x\in A^{(2)}\otimes_{p_0,\mathfrak S}N$. In $D(E_x)=E_x\otimes_{\mathfrak S}\mathcal L^\vee$, the element $\widetilde e\otimes e_1$ lifts $1\in\mathfrak S$. Since $\varphi_{\mathcal L^\vee}(e_1)=c(u)^{-1}e_1$ and $\tau_{\mathcal L^\vee}(1\otimes e_1)=\alpha^{-1}(1\otimes e_1)$, the cocycle $(\lambda_x,\mu_x)$ attached to $\widetilde e\otimes e_1$ by \eqref{eq:extension-lift} is
\begin{equation}\label{eq:nek-twisted-extension}
\lambda_x=c(u)^{-1}x\otimes e_1,\qquad
\mu_x=\alpha^{-1}b_x\otimes e_1.
\end{equation}
By Proposition~\ref{prop ambient-Ext-Herr}, $\mu_x\in N_{\cris}(\mathfrak M)$ and $(\tau_{\mathfrak M}-1)\lambda_x=(\varphi-1)\mu_x$.

\begin{rem}
\label{rem tate-descent-product}
\textup{(1)} The identity $(\varphi\circ\tau_{\mathcal L})(1\otimes e_{-1})=(\tau_{\mathcal L}\circ\varphi)(1\otimes e_{-1})$ gives $p_0(c(u))\varphi(\alpha)=\tau(p_0(c(u)))\alpha$. Since $c(u)=E(u)/c_0$ and $\tau$ fixes $c_0$, this becomes
\[
\varphi(\alpha)E(u_0)=E(u_1)\alpha.
\]
Moreover, applying the augmentation $A^{(3)}\to\mathfrak S$ to the cocycle condition $\partial_{02}(\alpha)=\partial_{01}(\alpha)\partial_{12}(\alpha)$ in $A^{(3)}$ gives $\operatorname{aug}(\alpha)=\operatorname{aug}(\alpha)^2$, so $\operatorname{aug}(\alpha)=1$ as $\alpha$ is a unit. Put $J=\Ker(\operatorname{aug}:A^{(2)}\to\mathfrak S)$. By the proof of Lemma~\ref{lem:small-ring-Frobenius-contraction}, $\varphi^n(J)\subset(p,E(u_0))^nJ$, hence $\varphi^n(\alpha)-1=\varphi^n(\alpha-1)\to0$ as $n\to\infty$. Rewriting the Frobenius relation as $E(u_0)/E(u_1)=\alpha/\varphi(\alpha)$ and taking the product of its $\varphi^j$-iterates for $0\leq j<n$ gives
\[
\prod_{j=0}^{n-1}\varphi^j\!\left(\frac{E(u_0)}{E(u_1)}\right)
=\alpha\varphi^n(\alpha)^{-1}
\xrightarrow{n\to\infty}
\alpha, \text{ hence } \alpha=\prod_{j\geq0}\varphi^j\!\left(\frac{E(u_0)}{E(u_1)}\right),
\]
where we consider the $(p,E(u_0))$-adic topology on $A^{(2)}$, and each factor is a unit because $E(u_0)$ and $E(u_1)$ generate the same prism ideal. 

\textup{(2)} For $x\in N^0$, the term $b_x$ in \eqref{eq:nek-extension-lift} is uniquely determined by the crystalline structure on $E_x$ from Proposition~\ref{prop nek-normalization}. Explicitly, compatibility with Frobenius and with the diagonal augmentation gives
\[
\bigl(\varphi-p_1(c(u))\bigr)b_x
=\tau_N(1\otimes x)-\varphi(\alpha)(1\otimes x),\qquad
(\operatorname{aug}\otimes1)(b_x)=0,
\]
where $\varphi$ denotes the semilinear Frobenius on $A^{(2)}\otimes_{p_0,\mathfrak S}N$ induced by $\varphi_N$.
\end{rem}

\begin{prop}
\label{prop nek-cochain-map}
The formulas
\[
\Psi_G^0(n)=n\otimes e_1,\qquad
\Psi_G^1(x)=\bigl(c(u)^{-1}x\otimes e_1,\alpha^{-1}b_x\otimes e_1\bigr)
\]
define a natural isomorphism of complexes
\[
\Psi_G^\bullet:C_{\mathrm{Nek}}^\bullet(G)
\xrightarrow{\sim}\tau_{\leq1}\mathcal C_{\varphi,\tau}^{\cris,\mathrm{fh}}(T),
\]
whose inverse in degree $1$ sends a cocycle $(\lambda,\mu)$ with $\lambda=a\otimes e_1$ \textup{($a\in N[1/E(u)]$)} to $c(u)a\in N^0$.
\end{prop}

\begin{proof}
By \eqref{eq:nek-twisted-extension}, $\Psi_G^1(x)=(\lambda_x,\mu_x)$ lies in $\Ker(d^1_{\mathrm{fh}})$. In the Baer sum of $E_x$ and $E_y$, whose underlying Breuil--Kisin extension is $E_{x+y}$, the induced lift satisfies \eqref{eq:nek-extension-lift} with $(x+y,b_x+b_y)$ in place of $(x,b_x)$. The uniqueness of the crystalline prismatic $(\varphi,\tau)$-structure therefore gives $b_{x+y}=b_x+b_y$, so $\Psi_G^1$ is additive. For any morphism $f:N\to N'$ in $\operatorname{Mod}^{\mathrm{Nek},[0,1]}_\varphi$, the induced map $E_x\to E_{f(x)}$ preserving the chosen lifts commutes with the $\tau$-operators, giving $(1\otimes f)(b_x)=b_{f(x)}$ and hence naturality of $\Psi_G^\bullet$. If $x=(\varphi_N-c(u))n$, identifying $\widetilde e$ with $e_{-1}+n$ in the split extension $N\oplus\mathcal L$ gives $b_x=\tau_N(1\otimes n)-\alpha(1\otimes n)$, hence $\Psi_G^1((\varphi_N-c(u))n)=d^0_{\mathrm{fh}}(n\otimes e_1)$. Thus $\Psi_G^\bullet$ is a morphism of complexes.

Clearly $\Psi_G^0$ is an isomorphism. For $\Psi_G^1$, take $(\lambda,\mu)\in\Ker(d^1_{\mathrm{fh}})$ and write $\lambda=a\otimes e_1$. By Propositions~\ref{prop ambient-Ext-Herr} and~\ref{prop nek-normalization}, $(\lambda,\mu)$ defines an admissible extension in $\operatorname{Mod}^{\cris,[-1,0]}_{\varphi,\tau}(\mathfrak S)$ whose tensor product with $\mathcal L$ has height at most $1$. For the lift obtained by tensoring with $e_{-1}$, \eqref{eq:nek-twisted-extension} shows that the Frobenius takes the form \eqref{eq:nek-extension-lift} with $x=c(u)a$, and the height-$1$ bound gives $x\in N^0$. Since its underlying Breuil--Kisin module is $E_x$, uniqueness of the crystalline $(\varphi,\tau)$-structure identifies its $\tau$-operator with that in \eqref{eq:nek-extension-lift}. Now, tensoring back with $\mathcal L^\vee$ recovers $(\lambda,\mu)=\Psi_G^1(x)$. Finally, $x\mapsto\Psi_G^1(x)$ is injective since $c(u)$ is invertible in $\mathfrak S[1/E(u)]$.
\end{proof}

\begin{cor}
\label{cor nek-hf-recovered}
For $T=T_pG$, there are natural isomorphisms
\[
H^1(C_{\mathrm{Nek}}^\bullet(G))
\simeq H^1\!\left(\mathcal C_{\varphi,\tau}^{\cris,\mathrm{fh}}(T)\right)
\simeq H_f^1(K,T).
\]
\end{cor}

\begin{proof}
Since $T_p G$ has Hodge--Tate weights in $[-1,0]$, the integral assertions follow from Proposition~\ref{prop nek-cochain-map} and Theorem~\ref{s6 thm-fh-main}.
\end{proof}

\subsection{The effective prismatic Herr complex}
\label{s6:effective}

\begin{defn}\label{def effective prismatic Herr complex}
Let $T\in\operatorname{Rep}^{\star}_{\mathbf Z_p}(G_K)$ have Hodge--Tate weights in $[0,h]$ for some $h\geq0$, so that $\mathfrak M=\mathfrak M(T)$ is effective. Define the \emph{effective prismatic Herr complex}
\begin{equation}\label{s6:eq-effective}
\mathcal C_{\varphi,\tau}^{\star}(\mathfrak M)
=\big[\mathfrak M\xrightarrow{d^0}
\mathfrak M\oplus N_\star(\mathfrak M)
\xrightarrow{d^1}N_\star(\mathfrak M)\big]
\end{equation}
in degrees $0,1,2$, with differentials given by the formulas of Definition~\ref{def ambient-finite-height-Herr}, and write $\mathcal C_{\varphi,\tau}^{\star}(T)=\mathcal C_{\varphi,\tau}^{\star}(\mathfrak M(T))$. 
\end{defn}

Since $\mathfrak M$ is effective, $\varphi_{\mathfrak M}(\mathfrak M)\subset\mathfrak M$ and $\varphi(N_\star(\mathfrak M))\subset N_\star(\mathfrak M)$, so \eqref{s6:eq-effective} is the total complex of the following
\[
\begin{tikzcd}[column sep=large,row sep=large]
\mathfrak M\arrow[r,"\varphi_{\mathfrak M}-1"]\arrow[d,"\tau_{\mathfrak M}-1"']
&\mathfrak M\arrow[d,"\tau_{\mathfrak M}-1"]\\
N_\star(\mathfrak M)\arrow[r,"\varphi-1"']&N_\star(\mathfrak M).
\end{tikzcd}
\]
When $T$ is crystalline, the comparison quasi-isomorphism of Theorem~\ref{prop compute prismatic coho} commutes with Frobenius, and taking fibers of $\varphi-1$ yields
\[
\mathcal C_{\varphi,\tau}^{\cris}(\mathfrak M)
\simeq\operatorname{fib}\!\left(
\varphi-1:R\Gamma(X_{\Prism},\mathfrak M_{\Prism})
\longrightarrow R\Gamma(X_{\Prism},\mathfrak M_{\Prism})\right).
\]
For $\star=\st$, the same formula holds with $(X,M_X)_{\Prism}$ in place of $X_{\Prism}$. Termwise inclusion defines a natural map of complexes
\begin{equation}\label{eq:effective-to-ambient-Herr}
\jmath_{\mathfrak M}:\mathcal C_{\varphi,\tau}^{\star}(\mathfrak M)
\longrightarrow\mathcal C_{\varphi,\tau}^{\star,\mathrm{fh}}(\mathfrak M).
\end{equation}

\begin{prop}
\label{prop Extension}\label{s6 prop-effective-extension}
Let $h\geq0$ and let $\mathfrak M\in\operatorname{Mod}^{\star,[0,h]}_{\varphi,\tau}(\mathfrak S)$. Every admissible extension of $\mathfrak S$ by $\mathfrak M$ in $\operatorname{Mod}^{\star,\mathrm{fh}}_{\varphi,\tau}(\mathfrak S)$ has an effective middle term of height at most $h$, and the inclusions of exact categories induce isomorphisms
\[
\operatorname{Ext}^1_{\operatorname{Mod}^{\star,[0,h]}_{\varphi,\tau}(\mathfrak S)}(\mathfrak S,\mathfrak M)
\simeq\operatorname{Ext}^1_{\operatorname{Mod}^{\star,\mathrm{eff}}_{\varphi,\tau}(\mathfrak S)}(\mathfrak S,\mathfrak M)
\simeq\operatorname{Ext}^1_{\operatorname{Mod}^{\star,\mathrm{fh}}_{\varphi,\tau}(\mathfrak S)}(\mathfrak S,\mathfrak M).
\]
Inside $\mathfrak M_E\oplus N_\star(\mathfrak M)$, we have $\Ker d^1=\Ker d^1_{\mathrm{fh}}$, and hence
\[
\tau_{\leq1}\mathcal C_{\varphi,\tau}^{\star}(\mathfrak M)
=\tau_{\leq1}\mathcal C_{\varphi,\tau}^{\star,\mathrm{fh}}(\mathfrak M).
\]
\end{prop}

\begin{proof}
Put $T=T_{\mathrm{cov}}(\mathfrak M)$. If $0\to\mathfrak M\to\mathfrak E\to\mathfrak S\to0$ is an admissible extension in $\operatorname{Mod}^{\star,\mathrm{fh}}_{\varphi,\tau}(\mathfrak S)$, applying $T_{\mathrm{cov}}$ gives a short exact sequence $0\to T\to U\to\mathbf Z_p\to0$ in $\operatorname{Rep}^{\star}_{\mathbf Z_p}(G_K)$ with $U=T_{\mathrm{cov}}(\mathfrak E)$. Since the Hodge--Tate weights of $T$ lie in $[0,h]$ and that of $\mathbf Z_p$ is $0$, the weights of $U$ lie in $[0,h]$ as well, so $\mathfrak E\simeq\mathfrak M(U)$ is effective of height at most $h$. The natural inclusions of exact categories induce the displayed
maps on $\operatorname{Ext}^1$. These are injective by full
faithfulness and surjective by the height bound just proved,
since admissibility is defined by the same underlying
$\mathfrak S$-module exactness. Here $T_{\mathrm{cov}}$ is used
only to establish the height bound, no exactness of its inverse
functor is needed.

Given $(\lambda,\mu)\in\Ker(d^1_{\mathrm{fh}})$, form the corresponding extension $\mathfrak E=\mathfrak M\oplus\mathfrak S e$ as in Proposition~\ref{prop ambient-Ext-Herr}. Since $\mathfrak E$ is effective, $\varphi_{\mathfrak E}(e)=e+\lambda$ lies in $\mathfrak E$, so $\lambda\in\mathfrak M$. The cocycle identity then already holds in $N_\star(\mathfrak M)\subset N_\star^{\mathrm{fh}}(\mathfrak M)$, giving $\Ker d^1_{\mathrm{fh}}\subset\Ker d^1$. The reverse inclusion and the equality of degree-$0$ terms and differentials are immediate from the definitions, so the good truncations coincide.
\end{proof}

\begin{thm}
\label{thm H0 H1}
Let $T\in\operatorname{Rep}^{\star}_{\mathbf Z_p}(G_K)$ have Hodge--Tate weights in $[0,h]$ for some $h\geq0$, and put $V=T[1/p]$. Then $H^0(\mathcal C_{\varphi,\tau}^{\star}(T))\simeq H^0(K,T)$, and in degree $1$ there are natural isomorphisms and an injection
\[
H^1\!\left(\mathcal C_{\varphi,\tau}^{\star}(T)\right)
\simeq\operatorname{Ext}^1_{\operatorname{Mod}^{\star,\mathrm{eff}}_{\varphi,\tau}(\mathfrak S)}
(\mathfrak S,\mathfrak M(T))
\simeq H^1\!\left(\mathcal C_{\varphi,\tau}^{\star,\mathrm{fh}}(T)\right)
\lhook\joinrel\longrightarrow H_\star^1(K,T),
\]
which becomes an isomorphism with $H_\star^1(K,V)$ after inverting $p$.
\end{thm}

\begin{proof}
Combine Proposition~\ref{s6 prop-effective-extension} with Theorem~\ref{s6 thm-fh-main}. 
\end{proof}

As we will see in Example~\ref{s6 example-missing-torsion}, the surjectivity without inverting $p$ can fail. See also \cite[Example~3.1 and Remark~3.2]{CGM26}.

\begin{thm}
\label{s6 prop-deformation-retract}\label{s6 prop-mod-u}
\label{s6 cor-representatives}\label{s6:cor-cohomology-formulas}
Let $\mathfrak M\in\operatorname{Mod}^{\star,\mathrm{eff}}_{\varphi,\tau}(\mathfrak S)$, and put $T=T_{\mathrm{cov}}(\mathfrak M)$ and $\overline{\mathfrak M}=\mathfrak M/u\mathfrak M$. There are natural quasi-isomorphisms of complexes of $\mathbf Z_p$-modules
\begin{equation}\label{s6:eq-mod-u}
\mathcal C_{\varphi,\tau}^{\star}(\mathfrak M)
\xrightarrow{P}\bigl[\mathfrak M\xrightarrow{\varphi_{\mathfrak M}-1}\mathfrak M\bigr]
\xrightarrow{Q}\bigl[\overline{\mathfrak M}\xrightarrow{\overline\varphi-1}\overline{\mathfrak M}\bigr],
\end{equation}
where both two-term complexes sit in degrees $0$ and $1$, $P$ is given by $P^0=\textup{id}$, $P^1(a,b)=a$, and $P^2=0$, and $Q$ is reduction modulo $u$ \textup{(with $\overline\varphi$ semilinear over the Frobenius on $W(k)$)}. In particular, $H^2(\mathcal C_{\varphi,\tau}^{\star}(\mathfrak M))=0$, $\Ker(\overline\varphi-1)$ computes $H^0(\mathcal C_{\varphi,\tau}^{\star}(\mathfrak M))$, and $\operatorname{Coker} (\overline\varphi-1)$ computes $H^1(\mathcal C_{\varphi,\tau}^{\star}(\mathfrak M))$.

In particular it implies when $T$ is crystalline with non-negative Hodge–Tate weights, $\mathcal C_{\varphi,\tau}^{\cris}(\mathfrak M)\hookrightarrow\mathcal C_{\varphi,\tau}^{\st}(\mathfrak M)$ is a quasi-isomorphism over $\mathbf Z_p$, and $H_f^1(K,T)=H_{\st}^1(K,T)$.
\end{thm}

\begin{proof}
Vanishing of $H^2(\mathcal C_{\varphi,\tau}^{\star}(\mathfrak M))$ is from Proposition~\ref{s6 prop-inverse}. Both $P$ and $Q$ are termwise surjective. The kernel of $P$ is $\bigl[N_\star(\mathfrak M)\xrightarrow{-(\varphi-1)}N_\star(\mathfrak M)\bigr]$ in degrees $1$ and $2$, which is acyclic by Proposition~\ref{s6 prop-inverse}, so $P$ is a quasi-isomorphism. Likewise, $\Ker(Q)=\bigl[u\mathfrak M\xrightarrow{\varphi_{\mathfrak M}-1}u\mathfrak M\bigr]$. Since $\varphi_{\mathfrak M}^j(u\mathfrak M)\subset u^{p^j}\mathfrak M$ and $\mathfrak M$ is $u$-adically complete, $-\sum_{j\geq0}\varphi_{\mathfrak M}^j$ converges on $u\mathfrak M$ and inverts $\varphi_{\mathfrak M}-1$ just as in Proposition~\ref{s6 prop-inverse}, so $Q$ is also a quasi-isomorphism. The crystalline-to-semistable inclusion is a quasi-isomorphism
by compatibility with the common two-term model in
\eqref{s6:eq-mod-u}.
By Theorem~\ref{thm H0 H1}, the inclusion
$H_f^1(K,T[1/p])\hookrightarrow H_{\st}^1(K,T[1/p])$
is therefore an isomorphism.
Taking inverse images in $H^1(K,T)$ gives the integral equality.
\end{proof}

\begin{rem}
The equality $H_f^1(K,T)=H_{\st}^1(K,T)$ in Theorem~\ref{s6 prop-deformation-retract} can also be seen directly from filtered $(\varphi,N)$-modules. Given a semistable extension $0\to V\to U\to\mathbf Q_p\to0$ where $V=T[1/p]$ is crystalline with weights in $[0,h]$ ($h\geq0$), lift $1\in D_{\st}(\mathbf Q_p)$ to $e\in D_{\st}(U)$. Since $\varphi(1)=1$ and $N(1)=0$, both $(\varphi-1)e$ and $N(e)$ belong to $D_{\cris}(V)$, on which $N=0$. Now applying $N\varphi=p\varphi N$ to $e$ yields $\varphi(N(e))=p^{-1}N(e)$. On the other hand, the effective Breuil--Kisin module of $T$ provides a $\varphi$-stable $W(k)$-lattice in $D_{\cris}(V)$, and scaling a nonzero $N(e)$ so that it lies in this lattice but not in $p$ times it would force $\varphi(N(e))=p^{-1}N(e)$ outside the lattice. Thus $N(e)=0$, so $N$ vanishes on all of $D_{\st}(U)$ and $U$ is crystalline.
\end{rem}

\begin{cor}
\label{s6 cor-congruences}
Let $\mathfrak M\in\operatorname{Mod}^{\star,\mathrm{eff}}_{\varphi,\tau}(\mathfrak S)$ and put $\overline{\mathfrak M}=\mathfrak M/u\mathfrak M$. For every $n\geq1$, there is a natural quasi-isomorphism
\[
\mathcal C_{\varphi,\tau}^{\star}(\mathfrak M)\otimes_{\mathbf Z_p}^{\mathbf L}\mathbf Z/p^n
\simeq\bigl[\overline{\mathfrak M}/p^n\xrightarrow{\overline\varphi-1}\overline{\mathfrak M}/p^n\bigr],
\]
whose degree-$1$ cohomology is $H^1(\mathcal C_{\varphi,\tau}^{\star}(\mathfrak M))/p^n$. 
\end{cor}

\begin{proof}
Since $\overline{\mathfrak M}$ is $p$-torsion-free, derived reduction modulo $p^n$ of the right-hand complex in \eqref{s6:eq-mod-u} is computed termwise, and the description of $H^1$ follows using $H^2(\mathcal C_{\varphi,\tau}^{\star}(\mathfrak M))=0$.
\end{proof}

\begin{rem}
\label{s6 cor-matrix}\label{s6 example-computations}
When $k$ is finite, $\bigl[\overline{\mathfrak M}\xrightarrow{\overline\varphi-1}\overline{\mathfrak M}\bigr]$ is a two-term complex of finite free $\mathbf Z_p$-modules of rank $[k:\mathbf F_p]\cdot\operatorname{rank}_{\mathfrak S}\mathfrak M$, its cohomology can be computed from the Smith normal form of the $\mathbf Z_p$-linear endomorphism $\overline\varphi-1$. For instance, if $k=\mathbf F_p$ and $T$ is an unramified rank-$1$ lattice with $\mathfrak M(T)=\mathfrak S e$ and $\varphi(e)=ae$ with $a\in\mathbf Z_p^\times$, the complex reduces to $[\mathbf Z_p\xrightarrow{a-1}\mathbf Z_p]$, with cohomology $(\mathbf Z_p,\mathbf Z_p)$ if $a=1$ and $(0,\mathbf Z/p^s)$ if $v_p(a-1)=s$ with $s\geq1$.
\end{rem}

\begin{eg}
\label{s6 example-missing-torsion}\label{s6 prop-integral-defect}
Assume $[K:\mathbf Q_p]<\infty$ and let $T\in\operatorname{Rep}^{\cris}_{\mathbf Z_p}(G_K)$ have Hodge--Tate weights in $[0,h]$ with $h\geq0$. By Theorems~\ref{thm H0 H1} and~\ref{s6 prop-mod-u} together with Remark~\ref{s6 cor-saturation}, the map on extension classes induced by $T_{\mathrm{cov}}$ fits into an exact sequence
\[
0\longrightarrow\operatorname{Coker}\bigl(\overline\varphi-1:
\overline{\mathfrak M}(T)\to\overline{\mathfrak M}(T)\bigr)
\longrightarrow H_f^1(K,T)\longrightarrow Q_\cris(T)\longrightarrow0,
\]
for some finite group $Q_\cris(T)$ that can be nonzero. 

Indeed, consider $T=\mathbf Z_p(-r)$ with $r\geq1$, we have $\mathfrak M(T)=\mathfrak S\{-r\}$ with Frobenius modulo $u$ equal to $p^r\sigma$ on $W(k)$, where $\sigma$ is the Witt Frobenius. Since $1-p^r\sigma$ has inverse $\sum_{j\geq0}p^{rj}\sigma^j$ on $W(k)$, Theorem~\ref{s6 prop-mod-u} gives $\mathcal C_{\varphi,\tau}^{\cris}(\mathbf Z_p(-r))\simeq0$, so $H^1(\mathcal C_{\varphi,\tau}^{\cris,\mathrm{fh}}(\mathbf Z_p(-r)))=0$ and $H_f^1(K,\mathbf Q_p(-r))=0$. As $H^0(K,\mathbf Q_p(-r))=0$, the cohomology sequence of $0\to\mathbf Z_p(-r)\to\mathbf Q_p(-r)\to(\mathbf Q_p/\mathbf Z_p)(-r)\to0$ yields
\[
H_f^1(K,\mathbf Z_p(-r))
\simeq H^0\bigl(K,(\mathbf Q_p/\mathbf Z_p)(-r)\bigr).
\]
When $K=\mathbf Q_p(\zeta_p)$ and $r=1$, the cyclotomic character has image $1+p\mathbf Z_p$, so $H_f^1(K,\mathbf Z_p(-1))=Q_\cris(\mathbf Z_p(-1))\simeq\mathbf Z/p$.

This reflects the failure of exactness also discussed in \cite[Example~3.1 and Remark~3.2]{CGM26}.
\end{eg}

\begin{rem}
\label{s6 comparison-summary}\label{s6 prop-principal-parts}
By Proposition~\ref{s6 prop-effective-extension}, for $\mathfrak M\in\operatorname{Mod}^{\star,\mathrm{eff}}_{\varphi,\tau}(\mathfrak S)$ the inclusion $\jmath_{\mathfrak M}$ in \eqref{eq:effective-to-ambient-Herr} is an equality after taking $\tau_{\leq1}$. For three-term complexes, the cokernel of $\jmath_{\mathfrak M}$ is the complex
\[
\mathcal P_\star(\mathfrak M)
=\big[0\longrightarrow\mathfrak M_E/\mathfrak M
\xrightarrow{\overline{\tau_{\mathfrak M}}-1}
N_\star^{\mathrm{fh}}(\mathfrak M)/N_\star(\mathfrak M)\big].
\]
From the short exact sequence of complexes
\[
0\longrightarrow\mathcal C_{\varphi,\tau}^{\star}(\mathfrak M)
\xrightarrow{\;\jmath_{\mathfrak M}\;}\mathcal C_{\varphi,\tau}^{\star,\mathrm{fh}}(\mathfrak M)
\longrightarrow\mathcal P_\star(\mathfrak M)\longrightarrow0,
\]
together with Proposition~\ref{s6 prop-effective-extension} and the vanishing $H^2(\mathcal C_{\varphi,\tau}^{\star}(\mathfrak M))=0$ (see Proposition~\ref{prop vanishing of H^2}), the induced long exact sequence gives $H^1(\mathcal P_\star(\mathfrak M))=0$ and an isomorphism
\[
H^2\!\left(\mathcal C_{\varphi,\tau}^{\star,\mathrm{fh}}(\mathfrak M)\right)
\simeq\operatorname{Coker}\!\left(\overline{\tau_{\mathfrak M}}-1:
\mathfrak M_E/\mathfrak M\longrightarrow N_\star^{\mathrm{fh}}(\mathfrak M)/N_\star(\mathfrak M)\right).
\]
In particular, $\overline{\tau_{\mathfrak M}}-1$ is injective (so any $a\in\mathfrak M_E$ with $(\tau_{\mathfrak M}-1)a\in N_\star(\mathfrak M)$ already belongs to $\mathfrak M$), and $\mathcal P_\star(\mathfrak M)$, hence also the cone of $\jmath_{\mathfrak M}$, is concentrated in degree $2$. We also expect that all cohomologies of $\mathcal P_\star(\mathfrak M)$ should vanish after inverting $p$.
\end{rem}

\subsection{Comparison with the classical Herr complex}
\label{s6:classical-comparison}

Put $C_p^\flat$ the tilt of $C_p$ and $\mathcal O_{\mathcal E}=\widehat{\mathfrak S[1/u]}_p$. Let $\widehat{\mathcal O_{\mathcal E}^{\mathrm{ur}}}$ be the $p$-adic completion of the maximal unramified extension of $\mathcal O_{\mathcal E}$ inside $W(C_p^\flat)$, and put $\mathcal O_{\mathcal E_\tau}=(W(C_p^\flat))^{G_L}$. Associated to a finite free $\mathbf Z_p$-module $T$ with continuous $G_K$-action its \'etale $(\varphi,\tau)$-module
\[
D(T)=\bigl(\widehat{\mathcal O_{\mathcal E}^{\mathrm{ur}}}\otimes_{\mathbf Z_p}T\bigr)^{G_{K_\infty}},\qquad
D_\tau(T)=\mathcal O_{\mathcal E_\tau}\otimes_{\mathcal O_{\mathcal E}}D(T),
\]
equipped with an additional operator $\tau_D$ (see \cite{Caruso_phitau}). With $\gamma$ as fixed in Section~\ref{sec notation}, \cite{Zha25} defines
\[
D_{\tau,0}(T)=\big\{z\in D_\tau(T) \,|\,(\gamma\otimes1)z
=\sum_{j=0}^{\chi(\gamma)-1}\tau_D^j(z)\big\}
\]
and proved that the integral Herr complex
\[
C_{\varphi,\tau}(T)
=\big[D(T)\xrightarrow{(\varphi-1,\tau_D-1)}
D(T)\oplus D_{\tau,0}(T)
\xrightarrow{(\tau_D-1)\ominus(\varphi-1)}D_{\tau,0}(T)\big]
\]
computes $H^i(K,T)$ in all degrees (see \cite[Theorem~2.12]{Zha25}). More precisely, we have quasi-isomorphism $C_{\varphi,\tau}(T)\simeq R\Gamma(K,T)$ (see \cite[Theorem~10.6]{GZ25}).

\begin{prop}
\label{prop embedding}
For $T\in\operatorname{Rep}^{\cris}_{\mathbf Z_p}(G_K)$, the embeddings of the coefficients induce a termwise injective morphism of complexes
\[\iota_{\mathrm{fh}}^\bullet:
\mathcal C_{\varphi,\tau}^{\cris,\mathrm{fh}}(T)\longrightarrow C_{\varphi,\tau}(T).
\]
Under $H^1(C_{\varphi,\tau}(T))\simeq H^1(K,T)$, the map
\[
H^1(\iota_{\mathrm{fh}}^\bullet):
H^1\!\left(\mathcal C_{\varphi,\tau}^{\cris,\mathrm{fh}}(T)\right)
\longrightarrow H^1(K,T)
\]
carries the class of an admissible module extension to that of the extension of $\mathbf Z_p$ by $T$ obtained via $T_{\mathrm{cov}}$. Consequently, the induced morphism
\[
\tau_{\leq1}\mathcal C_{\varphi,\tau}^{\cris,\mathrm{fh}}(T)
\longrightarrow R\Gamma(K,T)
\]
is an isomorphism on $H^0$ and an injection on $H^1$ with image $L_\cris(T)$, which equals $H_f^1(K,T)$ whenever $\operatorname{HT}(T)\subset[-r,0]$ for some $r\geq0$.

If $\mathfrak M(T)$ is effective, \textit{i.e.} $\operatorname{HT}(T)\subset[0,h]$ for some $h\geq0$, then $\iota_{\mathrm{eff}}^\bullet=\iota_{\mathrm{fh}}^\bullet\circ\jmath_{\mathfrak M(T)}$ fits into the commutative diagram
\[
\begin{tikzcd}[column sep=large,row sep=large]
\mathcal C_{\varphi,\tau}^{\cris}(T)
\arrow[r,hook,"\jmath_{\mathfrak M(T)}"]\arrow[dr,hook,"\iota_{\mathrm{eff}}^\bullet"']
&\mathcal C_{\varphi,\tau}^{\cris,\mathrm{fh}}(T)\arrow[d,hook,"\iota_{\mathrm{fh}}^\bullet"]\\
&C_{\varphi,\tau}(T),
\end{tikzcd}
\]
where $\jmath_{\mathfrak M(T)}$ induces isomorphisms on $H^0$ and $H^1$.
\end{prop}

\begin{proof}
Put $\mathfrak M=\mathfrak M(T)$, so scalar extension gives
\[
D(T)\simeq D(T_{\mathrm{cov}}(\mathfrak M))
\simeq\mathcal O_{\mathcal E}\otimes_{\mathfrak S}\mathfrak M.
\]
The element $E(u)$ is a unit in $\mathcal O_{\mathcal E}$ since it reduces modulo $p$ to $u^{\deg E}\in k((u))^\times$. Similarly, $A^{(2)}\hookrightarrow A_{\inf}\subset W(C_p^\flat)$ has $G_L$-fixed image and hence lands in $\mathcal O_{\mathcal E_\tau}$, where $E(u_0)$ is a unit because its reduction modulo $p$ is the nonzero element $\widetilde\pi^{\deg E}$ of the field $(C_p^\flat)^{G_L}$. The coefficient embeddings therefore extend to the localizations in \eqref{s6:eq-fh} and remain injective after tensoring with the finite free $\mathfrak S$-module $\mathfrak M$.

Under these embeddings, the subscript-zero equation for
$g=\gamma$ becomes the defining relation of $D_{\tau,0}(T)$.
Thus the maps in degrees $1$ and $2$ take values in
$D(T)\oplus D_{\tau,0}(T)$ and $D_{\tau,0}(T)$, respectively.
Compatibility with Frobenius and the $\tau$-operators then
shows that these maps commute with the differentials,
giving the termwise injective morphism of complexes
$\iota_{\mathrm{fh}}^\bullet$.

Given an admissible module extension $\mathfrak E$ with a lift $e$ of $1$ and cocycle $(\lambda,\mu)$ as in \eqref{eq:extension-lift}, scalar extension shows that the Frobenius and the identities of \eqref{eq:extension-lift} continue to hold in $D(T_{\mathrm{cov}}(\mathfrak E))$ and $D_\tau(T_{\mathrm{cov}}(\mathfrak E))$. The image of $(\lambda,\mu)$ in $D(T)\oplus D_{\tau,0}(T)$ therefore represents the $G_K$-extension $T_{\mathrm{cov}}(\mathfrak E)$, so $H^1(\iota_{\mathrm{fh}}^\bullet)$ sends $[(\lambda,\mu)]$ to $[T_{\mathrm{cov}}(\mathfrak E)]$. The remaining claims follow from Theorem~\ref{s6 thm-fh-main} and Proposition~\ref{s6 prop-effective-extension}.
\end{proof}

\section{\texorpdfstring{Relation with $F$-gauges and syntomic cohomology}{Relation with F-gauges and syntomic cohomology}}
\label{section F-gauge}

Let $X=\operatorname{Spf}(\mathcal O_K)$, with $K$ as in Section~\ref{sec notation}. Throughout this section, we work with the (non-logarithmic) prismatic site $X_{\Prism}$ and write $X^{\mathrm{Syn}}$ for the syntomification of $X$, from \cite[Definition~6.1.1]{bhatt-F-gauge}, recall that a prismatic $F$-gauge on $X$ is an object inside $D_{\mathrm{qc}}(X^{\mathrm{Syn}})$. For $\mathcal E\in D_{\mathrm{qc}}(X^{\mathrm{Syn}})$, put
\[
R\Gamma_{\mathrm{Syn}}(X,\mathcal E)
=R\Gamma(X^{\mathrm{Syn}},\mathcal E),\qquad
H^i_{\mathrm{Syn}}(X,\mathcal E)
=H^iR\Gamma_{\mathrm{Syn}}(X,\mathcal E).
\]
By \cite[Theorem~6.6.13]{bhatt-F-gauge}, a prismatic $F$-crystal $\mathfrak M_{\Prism}$ on $X_{\Prism}$ corresponds to a reflexive $F$-gauge $\mathcal E_{\mathfrak M}$ in the sense of \cite[Definitions~6.6.4 and~6.6.11]{bhatt-F-gauge}.

\begin{thm}
\label{thm main on coho of F-gauge}
\label{prop:effective-gauge-fibre}
Let $\mathfrak M_{\Prism}$ be an effective prismatic $F$-crystal on $X_{\Prism}$, and let $\mathfrak M$ be the crystalline prismatic $(\varphi,\tau)$-module obtained by evaluating $\mathfrak M_{\Prism}$ on the Breuil--Kisin prism with $\tau$ from the descent datum, as in Corollary~\ref{cor:crystal-tau-fh}. Then there are natural quasi-isomorphisms in $D(\mathbf Z_p)$:
\[
R\Gamma_{\mathrm{Syn}}(X,\mathcal E_{\mathfrak M})
\simeq\operatorname{fib}\!\left(
\varphi-1:
R\Gamma(X_{\Prism},\mathfrak M_{\Prism})
\longrightarrow R\Gamma(X_{\Prism},\mathfrak M_{\Prism})\right)
\simeq\mathcal C_{\varphi,\tau}^{\cris}(\mathfrak M).
\]
\end{thm}

\begin{proof}
Let $C_p$ be the completed algebraic closure fixed in Section~\ref{sec notation}, and let $R^\bullet$ be the $p$-completed \v{C}ech nerve of $\mathcal O_K\to\mathcal O_{C_p}$, with $n$-th term
\[
R^n=\left(\mathcal O_{C_p}^{\otimes_{\mathcal O_K}(n+1)}\right)^{\wedge}_p.
\]
As in the proof of \cite[Theorem~6.6.13]{bhatt-F-gauge}, each $R^n$ is a $p$-torsion-free quasiregular semiperfectoid $\mathcal O_{C_p}$-algebra. Let $(\Delta_n,I_n)$ be the initial prism inside $\Prism_{R^n}$, and put $M_n=\mathfrak M_{\Prism}(\Delta_n,I_n)$, which is a finite projective $\Delta_n$-module equipped with a semilinear Frobenius $\varphi_n$. Let $\mathcal E_{\mathfrak M,n}$ denote the pullback of $\mathcal E_{\mathfrak M}$ to $(R^n)^{\mathrm{Syn}}$.

By \cite[Definition~6.6.11]{bhatt-F-gauge}, the restriction of $\mathcal E_{\mathfrak M}$ to $(\mathcal O_{C_p})^{\mathrm{Syn}}$ is reflexive. Applying \cite[Lemma~6.6.10]{bhatt-F-gauge} along the map $\mathcal O_{C_p}\to R^n$ given by the first projection map identifies the filtration of $\mathcal E_{\mathfrak M,n}$ with the saturated Nygaard filtration:
\[
\operatorname{Fil}_N^iM_n
=\{m\in M_n:\varphi_n(m)\in I_n^iM_n \subset M_n[1/I_n]\}.
\]
Since $\mathfrak M_{\Prism}$ is assumed to be effective, we have $\varphi_n(M_n)\subset M_n$, and hence $\operatorname{Fil}_N^0M_n=M_n$.

Next, we compute $R\Gamma((R^n)^{\mathrm{Syn}},\mathcal E_{\mathfrak M,n})$. By \cite[Theorem~5.5.10 and Corollary~5.5.11]{bhatt-F-gauge}, the restriction of $\mathcal E_{\mathfrak M,n}$ to the filtered prismatization of $R^n$ has (derived) global sections $\operatorname{Fil}_N^0M_n$, since $M_n$ is derived $(p,I_n)$-complete and taking the degree-zero part of the Rees module is exact. Similarly, since the prismatization of $R^n$ is $\operatorname{Spf}(\Delta_n)$ by \cite[Theorem~5.5.7]{bhatt-F-gauge}, the restriction of $\mathcal E_{\mathfrak M,n}$ to the prismatization has (derived) global sections $M_n$. Using \cite[Example~6.1.7 and Remark~6.1.3]{bhatt-F-gauge}, the degree-zero filtered Frobenius and the filtration inclusion induce the two maps
\[
\begin{aligned}
\varphi_n:\operatorname{Fil}_N^0M_n&\longrightarrow M_n,
&m&\longmapsto\varphi_n(m),\\
1:\operatorname{Fil}_N^0M_n&\longrightarrow M_n,
&m&\longmapsto m.
\end{aligned}
\]
Because $\operatorname{Fil}_N^0M_n=M_n$, their difference is the $\mathbf Z_p$-linear endomorphism $\varphi_n-1:M_n\to M_n$, yielding
\[
R\Gamma((R^n)^{\mathrm{Syn}},\mathcal E_{\mathfrak M,n})
\simeq\operatorname{fib}(\varphi_n-1:M_n\longrightarrow M_n).
\]

By quasisyntomic descent for the prismatization and syntomification (see \cite[Remark~5.5.18 and the proof of Theorem~6.6.13]{bhatt-F-gauge}), taking mapping spectra from the structure sheaf to the respective coefficient complexes gives cosimplicial descent for derived global sections. For the prismatic side, \cite[Proposition~8.15]{BhattLuriePrismatization} identifies $R\Gamma(X_{\Prism},\mathfrak M_{\Prism})$ with the derived global sections on the prismatization of $X$, and evaluating on the prismatization $\operatorname{Spf}(\Delta_n)$ of $R^n$ recovers $M_n$. We thus obtain natural equivalences
\[
\begin{aligned}
R\Gamma_{\mathrm{Syn}}(X,\mathcal E_{\mathfrak M})
&\simeq\operatorname{Tot}_{[n]\in\Delta}
R\Gamma((R^n)^{\mathrm{Syn}},\mathcal E_{\mathfrak M,n}),\\
R\Gamma(X_{\Prism},\mathfrak M_{\Prism})
&\simeq\operatorname{Tot}_{[n]\in\Delta}M_n.
\end{aligned}
\]
Since the transition maps between the prisms $(\Delta_n,I_n)$ commute with Frobenius, the second equivalence fits into a commutative diagram
\[
\begin{tikzcd}[column sep=large,row sep=large]
R\Gamma(X_{\Prism},\mathfrak M_{\Prism})
\arrow[r,"\varphi-1"]\arrow[d,"\sim"']
&R\Gamma(X_{\Prism},\mathfrak M_{\Prism})\arrow[d,"\sim"]\\
\operatorname{Tot}_{[n]\in\Delta}M_n
\arrow[r,"\operatorname{Tot}(\varphi_n-1)"']
&\operatorname{Tot}_{[n]\in\Delta}M_n.
\end{tikzcd}
\]
Because fibers and totalizations are both limits in $D(\mathbf Z_p)$, so they commute. Taking horizontal fibers and substituting our description of $R\Gamma((R^n)^{\mathrm{Syn}},\mathcal E_{\mathfrak M,n})$ gives the first quasi-isomorphism. Finally, the comparison in Theorem~\ref{prop compute prismatic coho} commutes with Frobenius, and therefore identifies this fiber with $\mathcal C_{\varphi,\tau}^{\cris}(\mathfrak M)$ by Definition~\ref{def effective prismatic Herr complex}.
\end{proof}

\begin{cor}
\label{cor:effective-F-gauge-H2}
Under the hypotheses of Theorem~\ref{thm main on coho of F-gauge}, put $\overline{\mathfrak M}=\mathfrak M/u\mathfrak M$. There are natural quasi-isomorphisms in $D(\mathbf Z_p)$:
\[
R\Gamma_{\mathrm{Syn}}(X,\mathcal E_{\mathfrak M})
\simeq\bigl[\mathfrak M\xrightarrow{\varphi_{\mathfrak M}-1}\mathfrak M\bigr]
\simeq\bigl[\overline{\mathfrak M}\xrightarrow{\overline\varphi-1}\overline{\mathfrak M}\bigr],
\]
where both complexes are concentrated in degrees $0$ and $1$. In particular,
\[
H^i_{\mathrm{Syn}}(X,\mathcal E_{\mathfrak M})=0\qquad(i\geq2),
\]
and for every $n\geq1$ we have
\[
H^1_{\mathrm{Syn}}(X,\mathcal E_{\mathfrak M}/p^n)
\simeq H^1_{\mathrm{Syn}}(X,\mathcal E_{\mathfrak M})/p^n,
\]
where $\mathcal E_{\mathfrak M}/p^n$ denotes the cofiber of multiplication by $p^n$ on $\mathcal E_{\mathfrak M}$ in $D_{\mathrm{qc}}(X^{\mathrm{Syn}})$.
\end{cor}

\begin{proof}
Applying Theorem~\ref{s6 prop-mod-u} to the comparison in Theorem~\ref{thm main on coho of F-gauge} yields the displayed two-term complexes, from which the vanishing in degrees $i\geq2$ is immediate. Taking derived global sections of the cofiber sequence for multiplication by $p^n$ on $\mathcal E_{\mathfrak M}$ gives the exact sequence
\[
0\longrightarrow H^1_{\mathrm{Syn}}(X,\mathcal E_{\mathfrak M})/p^n
\longrightarrow H^1_{\mathrm{Syn}}(X,\mathcal E_{\mathfrak M}/p^n)
\longrightarrow H^2_{\mathrm{Syn}}(X,\mathcal E_{\mathfrak M})[p^n]
\longrightarrow0.
\]
Since $H^2_{\mathrm{Syn}}(X,\mathcal E_{\mathfrak M})=0$, the map $H^1_{\mathrm{Syn}}(X,\mathcal E_{\mathfrak M})/p^n\to H^1_{\mathrm{Syn}}(X,\mathcal E_{\mathfrak M}/p^n)$ is an isomorphism.
\end{proof}

We note that \cite[Lemma~4.11]{PentlandSyntomification} gives another explicit description of $R\Gamma_{\mathrm{Syn}}(X,\mathcal E_{\mathfrak M})$. Using this, he was able to prove \cite[Proposition~4.12]{PentlandSyntomification}, in particular, he was able to establishes the vanishing result in degree $2$ as stated in Corollary~\ref{cor:effective-F-gauge-H2}.

\begin{rem}
Let $T$ be a crystalline $G_K$-lattice with Hodge--Tate weights in $[0,h]$ for some $h\geq0$, and put $\mathfrak M=\mathfrak M(T)$. Combining Theorem~\ref{thm main on coho of F-gauge} with Theorem~\ref{thm H0 H1} identifies $H^0_{\mathrm{Syn}}(X,\mathcal E_{\mathfrak M})$ with $H^0(K,T)$ and yields in degree $1$ the composite
\[
H^1_{\mathrm{Syn}}(X,\mathcal E_{\mathfrak M})
\xrightarrow{\sim}
H^1\!\left(\mathcal C_{\varphi,\tau}^{\cris}(\mathfrak M)\right)
\lhook\joinrel\longrightarrow H_f^1(K,T),
\]
where the injection is induced by applying $T_{\mathrm{cov}}$ to the effective module extensions of Theorem~\ref{thm H0 H1}. This composite becomes an isomorphism after inverting $p$, though integral surjectivity can fail (see Example~\ref{s6 example-missing-torsion}).

For non-effective case, we note that the vanishing in degree $2$ fails in general, cf. \cite[Remark~6.7.4 and Footnote~86]{bhatt-F-gauge}). 

We also want to note that for crystalline lattices with Hodge--Tate weights in $[-r,0]$ ($r\geq0$), an integral comparison is proved in \cite[Corollary~A.9]{CGM26} (under the opposite Hodge--Tate sign convention).
\end{rem}

\appendix
\section{A convergence result of Frobenius iterates}
\label{s6:contraction}

To invert $\varphi-1$ via the geometric series $-\sum_{j\geq0}\varphi^j$, we need some estimation of $\varphi^j$ for the $(p,E(u_0))$-adic topology. This appendix aims to prove Proposition~\ref{s6 prop-inverse}.

\begin{lem}
\label{lem:small-ring-Frobenius-contraction}
Let $T\in\operatorname{Rep}^{\st}_{\mathbf Z_p}(G_K)$ have Hodge--Tate weights in $[0,h]$ for some $h\geq0$, and put $\mathfrak M=\mathfrak M(T)$. For every
\[
z\in\bigl(A_{\st}^{(2)}\cap W(\mathfrak m_{C_p}^{\flat})\bigr)
\otimes_{i_0,\mathfrak S}\mathfrak M,
\]
we have $\varphi^n(z)\to0$ in the $(p,E(u_0))$-adic topology as $n\to\infty$. When $T$ is crystalline, the same holds with $A^{(2)}$ and $p_0$ in place of $A_{\st}^{(2)}$ and $i_0$.
\end{lem}

\begin{proof}
Set $R=A_{\st}^{(2)}$, $J=\Ker(\operatorname{aug}:R\to\mathfrak S)$, and $I=(p,E(u_0))R$. Since $\operatorname{aug}$ is a morphism of $\delta$-rings splitting $i_0$, we have $R=\mathfrak S\oplus J$ as $\mathfrak S$-modules and $J\cap I=IJ$.

By \cite[Proposition~2.3.1(4)]{du-liu-prismaticphiGhatmodule}, $R/E(u_0)$ is a $p$-completed divided-power polynomial algebra over $\mathcal O_K$ with augmentation induced by $\operatorname{aug}$. Modulo $p$, every element of its divided-power augmentation ideal has $p$-th power $0$, since for $m>0$,
\[
\gamma_m(X)^p=\frac{(pm)!}{(m!)^p}\gamma_{pm}(X),
\qquad p\Big|\frac{(pm)!}{(m!)^p}.
\]
Thus $a^p\in J\cap I=IJ$ for all $a\in J$. Since $\operatorname{aug}$ is a $\delta$-map, $\delta(a)$ also lies in $J$, so $\varphi(a)=a^p+p\delta(a)\in IJ$, proving $\varphi(J)\subset IJ$.

By \cite[Lemma~3.3.1(2)]{du-liu-prismaticphiGhatmodule}, the map $R\to W(k)$ factors as $\operatorname{aug}$ followed by $\mathfrak S\to\mathfrak S/(u)=W(k)$, hence
\[
R\cap W(\mathfrak m_{C_p}^{\flat})=\operatorname{aug}^{-1}(u\mathfrak S)=u_0R+J.
\]
Put $e=\deg E(u)$ and pick $c\geq1$ with $p^c-1\geq e$. Because $u_0^e\in I$, we have $\varphi^c(u_0)=u_0^{p^c}\in Iu_0R$. Together with $\varphi(J)\subset IJ$ and $\varphi(I)\subset I$, this gives $\varphi^c(u_0R+J)\subset I(u_0R+J)$. Iterating yields
\[
\varphi^{cn}\!\left(R\cap W(\mathfrak m_{C_p}^{\flat})\right)
\subset I^n\left(R\cap W(\mathfrak m_{C_p}^{\flat})\right)
\qquad(n\geq0),
\]
and since $\varphi$ preserves both $I$ and $R\cap W(\mathfrak m_{C_p}^{\flat})$, the same inclusion holds with $\varphi^m$ on the left for every $m\geq cn$.

Because $\mathfrak M$ is effective, $\varphi_{\mathfrak M}^m(\mathfrak M)\subset\mathfrak M$ for all $m\geq0$. Writing $z=\sum_i r_i\otimes m_i$ with $r_i\in R\cap W(\mathfrak m_{C_p}^{\flat})$ and $m_i\in\mathfrak M$, we find $\varphi^m(z)=\sum_i\varphi^m(r_i)\otimes\varphi_{\mathfrak M}^m(m_i)\in I^n(R\cap W(\mathfrak m_{C_p}^{\flat}))\otimes_{i_0,\mathfrak S}\mathfrak M$ for $m\geq cn$, so $\varphi^m(z)\to0$ $I$-adically.

For $R=A^{(2)}$, the same argument applies using its diagonal augmentation, \cite[Proposition~2.2.8(2)]{du-liu-prismaticphiGhatmodule} for the divided-power quotient, and the proof of \cite[Lemma~3.3.1(3)]{du-liu-prismaticphiGhatmodule} for the reduction map to $W(k)$.
\end{proof}

\begin{prop}
\label{prop vanishing of H^2}\label{s6 prop-inverse}
For every $\mathfrak M\in\operatorname{Mod}^{\star,\mathrm{eff}}_{\varphi,\tau}(\mathfrak S)$, the endomorphism $\varphi-1$ of $N_\star(\mathfrak M)$ is bijective with continuous inverse $-\sum_{j\geq0}\varphi^j$, functorially in $\mathfrak M$. In particular, $H^2(\mathcal C_{\varphi,\tau}^{\star}(\mathfrak M))=0$.
\end{prop}

\begin{proof}
For each $g\in H_K$ with $\chi(g)\in\mathbf Z_{>0}$, the additive map
\[
\widehat{\mathfrak M}\longrightarrow\widehat{\mathfrak M},
\qquad z\longmapsto(g\otimes1)z-\sum_{j=0}^{\chi(g)-1}\tau_{\mathfrak M}^j(z)
\]
is $(p,E(u_0))$-adically continuous, so their common kernel $N_{\st}(\mathfrak M)$ is closed in $\widehat{\mathfrak M}$ and hence complete. For $z\in N_{\st}(\mathfrak M)$, Proposition~\ref{lem remove intersection with W} and Lemma~\ref{lem:small-ring-Frobenius-contraction} give $\varphi^j(z)\to0$. Since $\varphi$ preserves $N_{\st}(\mathfrak M)$, the series $b=-\sum_{j\geq0}\varphi^j(z)$ converges in $N_{\st}(\mathfrak M)$. Passing to the limit in
\[
(\varphi-1)\Bigl(-\sum_{j=0}^{n-1}\varphi^j(z)\Bigr)=z-\varphi^n(z)
\]
gives $(\varphi-1)b=z$, while any $z\in N_{\st}(\mathfrak M)$ with $(\varphi-1)z=0$ satisfies $z=\varphi^j(z)\to0$, so $z=0$. Since each $\varphi^j$ preserves $(p,E(u_0))^n\widehat{\mathfrak M}$ and commutes with morphisms in $\operatorname{Mod}^{\star,\mathrm{eff}}_{\varphi,\tau}(\mathfrak S)$, the inverse $-\sum_{j\geq0}\varphi^j$ is continuous and functorial.

In the crystalline case, Lemma~\ref{lem:small-ring-Frobenius-contraction} gives convergence of $-\sum_{j\geq0}\varphi^j(z)$ in $A^{(2)}\otimes_{p_0,\mathfrak S}\mathfrak M$ for its own $(p,E(u_0))$-adic topology, and continuity of the defining maps of $N_{\st}(\mathfrak M)$ ensures that the limit lies in $N_{\cris}(\mathfrak M)$, yielding bijectivity, continuity, and functoriality as above. Finally, $d^1(0,-b)=z$ shows that $d^1$ is surjective and hence $H^2(\mathcal C_{\varphi,\tau}^{\star}(\mathfrak M))=0$.
\end{proof}

\let\stdthebibliography\thebibliography
\let\stdendthebibliography\endthebibliography
\renewenvironment*{thebibliography}[1]{%
\stdthebibliography{GMW23}}
{\stdendthebibliography}

\bibliographystyle{amsalpha}

\bibliography{modified}

\end{document}